\documentclass[12pt]{article}
\usepackage{amsmath,amssymb,graphicx,latexsym}
\newcommand{\qed}{\hfill\vrule height6pt  width6pt
depth0pt \medskip}
\newcommand{\pr}{\mathbb{P}}
\newcommand{\arct}{\tan^{-1}}
\newcommand{\arcs}{\sin^{-1}}

\newcommand{\sgn}{\textrm{sgn}\,}
\newcommand{\sq}{\sqrt{2}/2}
\newcommand{\tsq}{\theta\sqrt{2}/2}
\newcommand{\bo}{\beta_1}
\newcommand{\bt}{\beta_2}
\newcommand{\gamo}{\gamma_1}
\newcommand{\gamt}{\gamma_2}
\newcommand{\nequiv}{\equiv\!\!\!\!\!\!/\ }
\title{Limit shapes for random square Young tableaux and plane
partitions}
\author{Dan Romik, Boris Pittel}
\begin{document}

\begin{center} \textsc{\large{Limit shapes for random square Young
tableaux and plane partitions}} \\\ \\
\footnotesize
\begin{tabular}{cc}
\textsc{\large{Boris Pittel}}
\footnote{Supported in part by the NSF Grant DMS-0104104.}
 & \textsc{\large{Dan Romik}} \\
 Department of Mathematics    & Department of Mathematics \\
 Ohio State University        & The Weizmann Institute of Science \\
 Columbus, Ohio 43210         & Rehovot 76100, Israel \\
 email: \tt{bgp@math.ohio-state.edu}
& email: \tt{romik@wisdom.weizmann.ac.il}
\end{tabular}
\end{center}

\footnotesize
\paragraph{Abstract.}Our main result is a limit shape theorem for the
two-dimensional surface defined by a uniform random $n\times n$ square
Young tableau. The analysis leads to a calculus of variations
minimization problem that resembles the minimization problems studied
by Logan-Shepp, Vershik-Kerov, and Cohn-Larsen-Propp. Our solution
involves methods from the theory of singular integral equations, and
sheds light on the somewhat mysterious derivations in these works. An
extension to rectangular diagrams, using the same ideas but involving
some nontrivial computations, is also given.

We give several applications of the main result. First, we show that
the location of a particular entry in the tableau is in the limit
governed by a semicircle distribution.

Next, we derive a result on the length of the longest increasing
subsequence in segments of a \emph{minimal Erd\"os-Szekeres
permutation}, namely a permutation of the numbers $1,2,\ldots,n^2$ whose
longest monotone subsequence is of length $n$ (and hence minimal by
the Erd\"os-Szekeres theorem).

Finally, we prove a limit shape theorem for the surface defined by a
random plane partition of a very large integer over a large square
(and more generally rectangular) diagram.

\vspace{40.0 pt}

\vspace{45.0 pt}
\noindent
\begin{tabular}{c} \qquad\qquad\qquad
 \qquad\qquad\qquad\qquad\qquad\qquad \\\hline \end{tabular}

$^1$Supported in part by the NSF Grant DMS-0104104.

\newpage
\tableofcontents

\newpage
\section{Introduction}

\subsection{Random square Young tableaux}

In this paper, we study the large-scale asymptotic behavior of uniform
random Young tableaux chosen from the set of tableaux of square
shape. Recall that a \emph{Young diagram} is a graphical
representation of a partition $\lambda:
\lambda(1)\ge\lambda(2)\ge\ldots\ge\lambda(k)$ of $n=\sum \lambda_i$ as an
array of \emph{cells}, where row $i$ has $\lambda_i$ cells. For a
Young diagram $\lambda$ (we will often identify a partition with its
Young diagram), a \emph{Young tableau} of shape $\lambda$ is a filling
of the cells of $\lambda$ with the numbers $1,2,\ldots,n$ such that the
numbers along every row and column are increasing.

A square Young tableau is a Young tableau whose shape is an $n\times
n$ square Young diagram. The number of such tableaux is
known by the hook formula of Frame-Thrall-Robinson (see
\eqref{eq:hook} below) to be
$$ \frac{(n^2)!}{[1\cdot (2n-1)][2\cdot(2n-2)]^2[3\cdot(2n-3)]^3
\ldots [(n-1)(n+1)]^{n-1} \ n^n }. $$ A square tableau
$T=(t_{i,j})_{i,j=1}^n$ can be depicted geometrically as a
three-dimensional stack of cubes over the two-dimensional square
$[0,n]\times[0,n]$, where $t_{i,j}$ cubes are stacked over the square
$[i-1,i]\times[j-1,j]\times\{0\}$. Alternatively, the function
$(i,j)\to t_{i,j}$ can be thought of as the graph of the
(non-continuous) surface of the upper envelope of this stack.  By
rescaling the $n\times n$ square to a square of unit sides, and
rescaling the heights of the columns of cubes so that they are all
between 0 and 1, one may consider the family of square tableaux as
$n\to\infty$. This raises the natural question, whether the shape of
the stack for a \emph{random} $n\times n$ square tableau exhibits some
asymptotic behavior as $n\to\infty$. The answer is given by the
following theorem, and is illustrated in Figure 1.

\paragraph{Theorem 1.} Let ${\cal T}_n$ be the set of $n\times n$
square Young tableaux, and let $\pr_n$ be the uniform probability
measure on ${\cal T}_n$. Then for the function
$L:[0,1]\times[0,1]\to[0,1]$ defined below, we have:

(i) \emph{Uniform convergence to the limit shape:} for all $\epsilon>0$,
$$ \pr_n\left( T\in{\cal T}_n : \max_{1\le i,j\le n}
\left|\frac{1}{n^2}t_{i,j} -
L\left(\frac{i}{n},\frac{j}{n}\right)\right| > \epsilon \right)
\xrightarrow[n\to\infty]{} 0.$$ 

(ii) \emph{Rate of convergence in the interior of the square:} for all
$\epsilon>0$,
$$ \pr_n\bigg( T\in{\cal T}_n : \max_{ \tiny{
\begin{array}{c} 1\le i,j\le n \\ 
  \min(i j,(n-i)(n-j)) > n^{3/2+\epsilon} \end{array}}}
\left|\frac{1}{n^2}t_{i,j} - L\left(\frac{i}{n},\frac{j}{n}\right)\right| > 
\frac{1}{n^{(1-\epsilon)/2}} \bigg) 
\xrightarrow[n\to\infty]{} 0. $$

\newpage

\begin{figure}[h!]
\hspace{-50.0pt}
\hbox{$
\begin{array}{c}
  \includegraphics{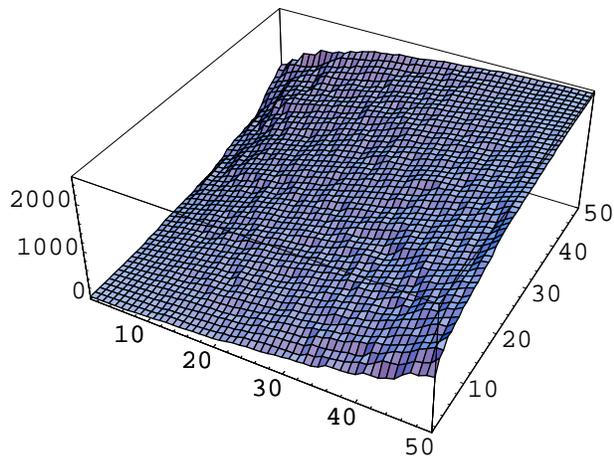} \\
  \textrm{\footnotesize{(a) 3D plot of simulated tableau}}
\end{array}$ $
\begin{array}{c}
  \includegraphics{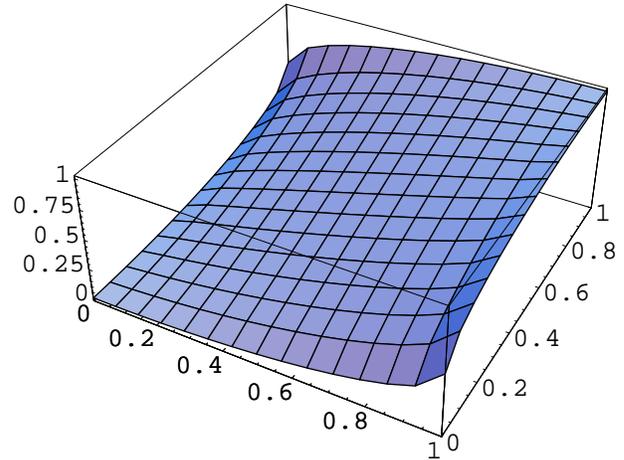} \\
  \textrm{\footnotesize{(b) The limit surface $L(x,y)$}}
\end{array}$
}
\medskip
\hspace{-50.0pt}
\hbox{$
\begin{array}{c}
  \includegraphics{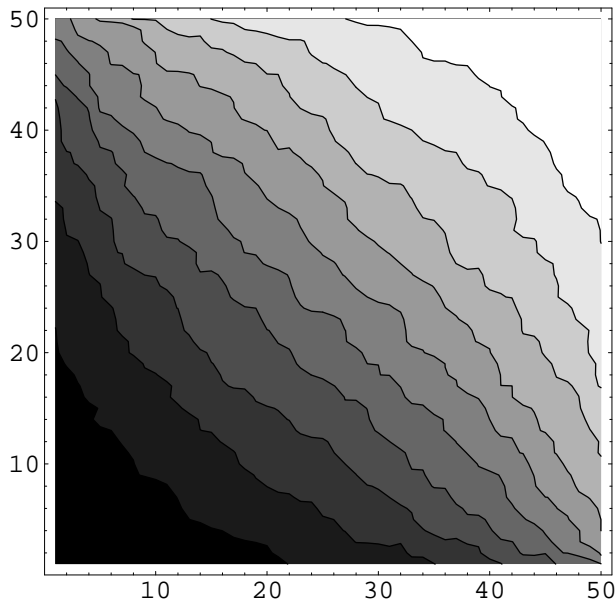} \\
  \textrm{\footnotesize{(c) Contour plot of simulated tableau}}
\end{array}$ $
\begin{array}{c}
  \includegraphics{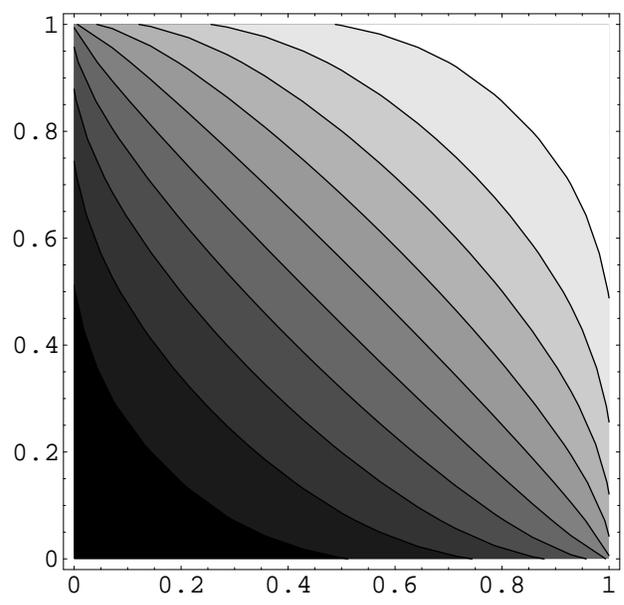} \\
  \textrm{\footnotesize{(d) Contour plot of $L$}}
\end{array}$
}
\caption{A simulated $50\times 50$ random tableau and the limit surface}
\end{figure}

\paragraph{Definition of $L$.} We call the function $L$ the
\emph{limit surface of square Young tableaux}. It is defined by the
implicit equation
$$ x+y =\frac{2}{\pi}(x-y)\arct\left(
\frac{ (1-2L(x,y))(x-y)}{\sqrt{4L(x,y)(1-L(x,y)) -
(x-y)^2}} \right) $$
$$ + \frac{2}{\pi}\arct \left(\frac{\sqrt{4L(x,y)(1-L(x,y))-
(x-y)^2}}{1-2L(x,y)}\right)     $$
for $0 \le y \le 1-x \le 1$, together with the reflection property
$$ L(x,y) = 1- L(1-x,1-y) $$
(where $\arct$ is the arctangent function).
It is more natural to describe $L$ in terms of its level
curves $\{L(x,y)=\alpha\}$. First, introduce the \emph{rotated coordinate
system}
\begin{equation}\label{eq:rotated}
u=\frac{x-y}{\sqrt{2}},\qquad v=\frac{x+y}{\sqrt{2}}.
\end{equation}
In the $u-v$ plane, the
square $[0,1]\times[0,1]$ transforms into the rotated square
$$ \Diamond = \{ (u,v)\in \mathbb{R}^2 : |u|\le \sqrt{2}/2,
|u|\le v \le \sqrt{2}-|u| \}. $$
Now define the one-parameter family of functions
$(g_\alpha)_{0\le\alpha\le 1}$ given by
$$g_\alpha : [-\sqrt{2\alpha(1-\alpha)},\sqrt{2\alpha(1-\alpha)}]\to
\mathbb{R}, $$
\begin{multline}\label{eq:minimizers} {\scriptsize
 g_\alpha(u) = \left\{ \begin{array}{ll}
\ \ \frac{2}{\pi}
u\arct\left(\frac{(1-2\alpha)u}{\sqrt{2\alpha(1-\alpha)-u^2}}
\right) +
\frac{\sqrt{2}}{\pi}
\arct
\left(\frac{\sqrt{2(2\alpha(1-\alpha)-u^2)}}{1-2\alpha}\right)
& 0 \le \alpha < \frac{1}{2}, \\
-\frac{2}{\pi}
u\arct\left(\frac{(1-2\alpha)u}{\sqrt{2\alpha(1-\alpha)-u^2}}
\right) -
\frac{\sqrt{2}}{\pi}
\arct
\left(\frac{\sqrt{2(2\alpha(1-\alpha)-u^2)}}{1-2\alpha}\right) +\sqrt{2}
& \frac{1}{2} < \alpha \le 1,  \\
\frac{\sqrt{2}}{2} & \alpha=\frac{1}{2}.
\end{array}\right.}
\end{multline}
Then in the rotated coordinate system, the surface $\bar{L}(u,v) =
L(x(u,v),y(u,v))$ can be described as the surface whose level curves
$\{\bar{L}(u,v)=\alpha\}$ are exactly the curves
$\{v=g_\alpha(u)\}$. That is,
$$ \{ (u,v)\in\Diamond : \bar{L}(u,v) = \alpha \} =
  \{ (u,v)\in\Diamond : |u|\le \sqrt{2\alpha(1-\alpha)},
v=g_\alpha(u) \}. $$
This is illustrated in Figure 2. It is straightforward to check that
the curves $v=g_\alpha(u)$ do not intersect, and so define a surface
\footnote{See equation \eqref{eq:partialdiff} in section 3.4.}.

\begin{figure}
\begin{center}
\includegraphics{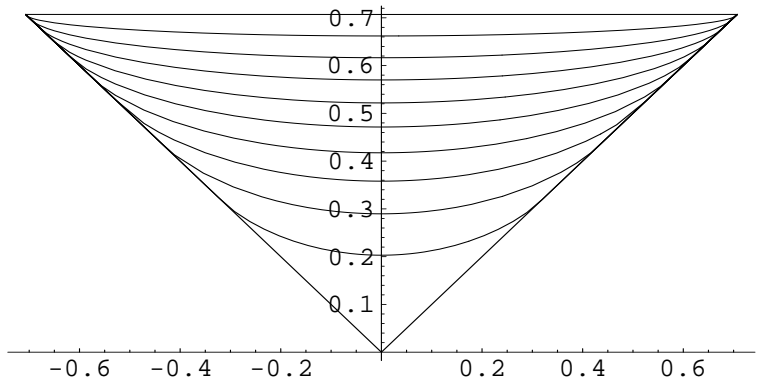}
\caption{The curves $v=g_\alpha(u)$ for
$\alpha = 0.05,\ 0.1,\ 0.15,\ 0.2,\ \ldots\ ,\ 0.5$}
\end{center}
\end{figure}

Note some special values of $L(x,y)$ which can be computed explicitly:
\begin{eqnarray*} L(t,0) &=& L(0,t)\ \, =\ \, \frac{1-\sqrt{1-t^2}}{2}, \\
L(t,1) &=& L(1,t)\ \,=\ \,\frac{1+\sqrt{2t-t^2}}{2}, \\
L(t,t) &=& \frac{1-\cos(\pi t)}{2}.
\end{eqnarray*}

\bigskip
The approach in proving Theorem 1 is the variational approach. Namely,
we identify the large-deviation rate functional of the level curves of
the random surface defined by the tableau, then analyze the functional
and find its minimizers. This will give Theorem 1(ii), with the rate
of convergence following from classical norm estimates for some
integral operators. The treatment of the boundary of the square,
required for Theorem 1(i), turns out to be more delicate, and will
require special arguments.

\subsection{Location of particular entries}

Theorem 1 identifies the approximate value of the entry of a typical square
tableau in a given location in the square. A dual outlook is to ask
where a given value $k$ will appear in the square tableau, since all
the values between 1 and $n^2$ appear exactly once. These questions
are almost equivalent. Indeed, if $k$ is approximately $\alpha\cdot
n^2$, then Theorem 1 predicts that with high probability the entry $k$
will appear in the vicinity of the level curve $\{L(x,y) = \alpha\}$
(the fact that this actually follows from Theorem 1 is a simple
consequence of the monotonicity property of the tableau along rows and
columns). However, one may ask a more detailed question about the
limiting distribution of the location of the entry $k$ on the level
curve. It turns out that its $u$-coordinate has approximately the semicircle
distribution. This is made precise in the following theorem.

\paragraph{Theorem 2.} For a tableau $T\in{\cal T}_n$ and $1\le k\le
n^2$, denote by $(i(T,k),j(T,k))$ the location of the entry $k$ in
$T$, and denote $X(T,k)=i(T,k)/n,\ Y(T,k)=j(T,k)/n$. Let $0<\alpha<1$,
let $k_n$ be a sequence of integers such that $k_n/n^2
\xrightarrow[n\to\infty]{} \alpha$, and for each $n$ let $T_n$ be a
uniform random tableau in ${\cal T}_n$. Then as $n\to\infty$, the
random vector $(X(T_n,k_n),Y(T_n,k_n))$ converges in distribution to
the random vector 
$$ (X_\alpha, Y_\alpha) := \left( \frac{V_\alpha+U_\alpha}{2},
 \frac{V_\alpha-U_\alpha}{2} \right), $$
where $U_\alpha$ is a random variable with density function
\begin{equation}\label{eq:semicircle}
f_{U_\alpha}(u) = \frac{\sqrt{2\alpha(1-\alpha)-u^2}}{\pi \alpha(1-\alpha)}
\mathbf{1}_{[-\sqrt{2\alpha(1-\alpha)},
\sqrt{2\alpha(1-\alpha)}]}(u)
\end{equation}
and $V_\alpha = g_\alpha(U_\alpha)$.

\bigskip
Theorem 2 is one of several aspects of our work which shows a deep
connection to the work of Logan-Shepp and Vershik-Kerov on the limit
shape of Plancherel-random partitions - see section 8 for discussion.

\subsection{Minimal Erd\"os-Szekeres permutations} The famous
Erd\"os-Szekeres theorem states that a permutation of $1,2,\ldots,n^2$
must have either an increasing subsequence of length $n$ or a
decreasing subsequence of length $n$. This can be proved using the
pigeon-hole principle, but also follows from the RSK correspondence
using the observation that a Young diagram of area $n^2$ must have
either width or height at least $n$.

For the width and height of a Young diagram of area $n^2$ to be
\emph{exactly} $n$, the diagram must be a square. From the RSK
correspondence it thus follows that to each permutation of
$1,2,\ldots,n^2$ whose longest increasing subsequence and longest
decreasing subsequence have length exactly $n$, there correspond a
pair of square $n\times n$ Young tableaux. Such a permutation has the
minimal possible length of a longest \emph{monotone} subsequence, and
it seems appropriate to term such permutations \emph{minimal
Erd\"os-Szekeres permutations} (we are not aware of any previous
references to these permutations, aside from a brief mention in
\cite{knuth}, exercise 5.1.4.9).

As an application of our limit shape result, we will prove the
following result on the length of the longest increasing subsequence
when just an initial segment of a random minimal Erd\"os-Szekeres
permutation is read.

\paragraph{Theorem 3.} For each $n$, let $\pi_n$ be a uniform random minimal
Erd\"os-Szekeres permutation of $1,2,\ldots,n^2$. For $1\le k\le n^2$, Let
$l_{n,k}$ be the length of the longest increasing subsequence in the
sequence $\pi_n(1), \pi_n(2),\ldots,\pi_n(k)$. Denote $\alpha=k/n^2$, and 
$\alpha_0=n^{-2/3+\epsilon}$. Then for any $\epsilon>0$, and 
$\omega(n)\to\infty$ however slowly,
$$ \max_{\alpha_0\le k/n^2\le 1/2}
\pr( |l_{n,k} -
2\sqrt{\alpha(1-\alpha)}n|> \alpha_0^{1/2}\omega(n)n) \xrightarrow[n\to\infty]{} 0.
$$

Thus the random fluctuations of $l_{n,k}$ around $2\sqrt{\alpha(1-\alpha)}n$ are not
likely to be of order substantially larger than $n^{2/3}$.

\subsection{A limit surface for random square plane partitions} 

Another probability model which was studied in the context of limit
shapes, is that of random plane partitions. If $\lambda$ is a Young
diagram, \emph{a plane partition of $n$ of shape $\lambda$} is an
array of positive integers $(p_{i,j})_{(i,j)\in\lambda}$ indexed by
the cells of $\lambda$, that sum to $n$ and which are weakly
decreasing along rows and columns, i.e., satisfy
$$ p_{i,j} \ge p_{i,j+1}, \qquad p_{i,j} \ge p_{i+1,j}. $$
Cerf and Kenyon \cite{cerfkenyon} proved a limit shape result for
random unrestricted plane partitions of an integer $n$. Cohn, Larsen
and Propp \cite{cohnlarspropp} proved a limit shape result for random plane
partitions whose three-dimensional graph is bounded inside a large box
of given relative proportions.

We apply Theorem 1 to prove a limit shape result for random plane
partitions of $m$ defined over an $n\times n$ square Young diagram,
when $m$ is much greater than $n^6$. This can be related to square
Young tableaux by the observation that when a plane partition does not
contain repeated entries (which in the asymptotic regime described
above happens with high probability), the order structure on the
entries of the plane partition is a Young tableau. The precise result is
the following.

\paragraph{Theorem 4.} For integers $n,m>0$, let ${\cal P}_{n,m}$ be the set of
plane partitions of $m$ of $n\times n$ square shape, and let
$\pr_{n,m}$ be the uniform probability measure on ${\cal
P}_{n,m}$. If $\pi=(p_{i,j})_{i,j=1}^n$ is an element of ${\cal
P}_{n,m}$, let its \emph{rescaled surface graph} be the function
$ \tilde{S}_\pi:[0,1)\times[0,1)\to[0,\infty)$ defined by
$$ \tilde{S}_\pi(x,y) = \frac{n^2}{m}p_{\lfloor n x\rfloor+1, \lfloor
n y\rfloor+1}. $$ 

Suppose $m$ and $n$ are sequences of integers that tend to infinity in
such a way that $m/n^6\to\infty$. Then for all
$\epsilon>0,\ x,y\in[0,1)$ we have
$$ \pr_{n,m}( \pi\in {\cal P}_{n,m} : |\tilde{S}_\pi(x,y)-M(x,y)|>
\epsilon ) \xrightarrow[\qquad]{} 0, $$
where $M:[0,1]\times[0,1]\to[0,\infty)$ is given by
$$ M(x,y) = -\log(L(x,y)). $$

\medskip
Theorem 4 may be related to a limiting case $\gamma\to\infty$ in the
limit shape result of Cohn-Larsen-Propp \cite{cohnlarspropp}. We have
not attempted to check this.

\subsection{Random rectangular Young tableaux and plane partitions} 

The methods which we will use to prove Theorems 1, 2, and 4 work
equally well for rectangular Young tableaux and plane partitions, in
the limit when the size of the rectangle grows and its relative
proportions tend to a limiting value $\theta>0$. For each possible value
$\theta$ of the ratio between the sides of the rectangle, there is a
limiting surface $L_\theta$ for random rectangular Young tableaux, and
a limiting surface $M_\theta$ for random rectangular plane
partitions. 
Analogously to the square tableaux, the rectangular $n_1\times n_2$
tableaux can be viewed as the result of applying the RSK algorithm to
a permutation of $\{1,\ldots,n_1 n_2\}$ with the property that the
lengths of the longest increasing and the longest decreasing
subsequences are exactly equal $n_1$ and $n_2$ (by the
Erd\"os-Szekeres theorem, the two lengths cannot be simultaneously
below $n_1$ and $n_2$, respectively). The proofs, which we include at
the end of the paper, require some nontrivial modifications, but the
final results are unexpectedly as elegant as for the square case.

Let $\theta>0$. We may assume that $\theta \le 1$, otherwise exchange
the two sides of the rectangle. Define
$L_\theta:[0,1]\times[0,\theta]\to[0,1]$, the \emph{limit surface of
rectangular tableaux with side ratio $\theta$}, as follows. For each
$0<\alpha<1$, the $\alpha$-level curve $\{(x,y):
L_\theta(x,y)=\alpha\}$ is given in rotated $u-v$ coordinates by
$$ \{ (u,h_{\theta,\alpha}(u)) : -\beta_1 \le u \le
\beta_2 \}, $$ where
%
%
%
%
%
%
\begin{eqnarray*}
\overline{\beta} &=& \sqrt{2\theta\alpha(1-\alpha)}, \\
\beta_1 &=& \overline{\beta} - \alpha(1-\theta)\sqrt{2}/2, \qquad
\beta_2\ \ =\ \  \overline{\beta} + \alpha(1-\theta)\sqrt{2}/2, \\
h_{\theta,\alpha}(u) &=& \theta\sqrt{2}/2 \pm
 (\beta_1 - \theta\sqrt{2}/2) +
\frac{2\overline{\beta}}{\pi} \bigg[ \pm (-\xi-\gamma_1) \arct
\sqrt{\frac{(1-\xi)(\gamma_1-1)}{(1+\xi)(\gamma_1+1)}} \\
& & +(\xi-\gamma_2)\arct
\sqrt{\frac{(1+\xi)(\gamma_2-1)}{(1-\xi)(\gamma_2+1)}} \\ & &
+ \frac{1}{2}\left( \arcs \xi+\frac{\pi}{2}\right)
\frac{1-\theta}{\sqrt{2}\, \overline{\beta}} \pm
\frac{\pi}{2}(\gamma_1-1) \bigg], \quad\qquad 0<\alpha\le \frac{1}{2},\\
\pm &=& \left\{ 
  \begin{array}{ll} + & \quad 0 < \alpha \le \theta/(1+\theta), \\
                    - & \quad \theta/(1+\theta) < \alpha \le 1/2,
  \end{array} \right. \\
\xi &=& \frac{u-\alpha(1-\theta)\sqrt{2}/2}{\overline{\beta}}, \qquad
        u \in [-\beta_1,\beta_2], \\
\gamma_1 &=& \frac{\alpha+\theta(1-\alpha)}{\sqrt{2}\,\overline{\beta}},
\qquad
\gamma_2\ \ =\ \  \frac{\theta\alpha+1-\alpha}{\sqrt{2}\,\overline{\beta}},
\\
h_{\theta,\alpha}(u) &=& (1+\theta)\sqrt{2}/2 -
h_{\theta,1-\alpha}((1-\theta)\sqrt{2}/2-u), \qquad \frac{1}{2} <
\alpha <1,
\end{eqnarray*}
see Figure 3. Set $$M_\theta(x,y) = -\log(L_\theta(x,y)).$$

\paragraph{Theorem 5.} For integers $n,m>0$, let ${\cal T}_{n,m}$ be
the set of tableaux whose shape is an $n\times m$ rectangular diagram,
and let $\pr_{n,m}$ be the uniform probability measure on
$\pr_{n,m}$. If
$T=(t_{i,j})_{i,j}\in{\cal T}_{n,m}$, define the rescaled tableau
surface of $T$ as the function
$\tilde{S}_T:[0,1)\times[0,m/n)\to[0,1]$ given by
$$ \tilde{S}_T(x,y) = \frac{1}{n m}t_{\lfloor n x\rfloor+1,
\lfloor n y \rfloor+1}. $$
Let $0<\theta\le 1$. If $m_n$ is a sequence of integers such
that $m_n/n \to\theta$ as $n\to\infty$, then for all $\epsilon>0,\ 
x\in[0,1), y\in [0,\theta)$,
$$ \pr_{n,m_n}(T\in{\cal T}_{n,m_n} : |\tilde{S}_T(x,y)-L_\theta(x,y)|
> \epsilon ) \xrightarrow[n\to\infty]{} 0. $$

\begin{figure}
\begin{center}
\includegraphics{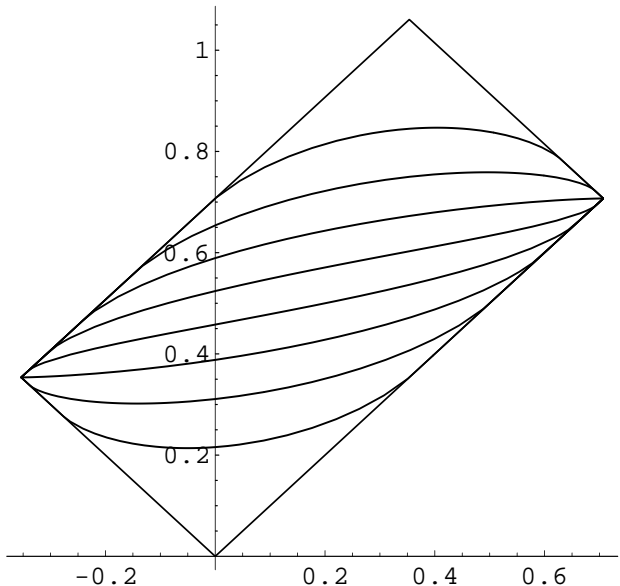}
\caption{The curves $h_{\theta,\alpha}$ for $\theta=0.5$,
$\alpha=k/9$, $k=1,2,\ldots,8$.}
\end{center}
\end{figure}

\paragraph{Theorem 6.} For each $\theta>0$, $M_\theta$ is the limit
surface of uniform random plane partitions of $m_n$ over a rectangular
diagram of sides $n$ and $k_n$, provided $k_n/n \to \theta$ and
$m_n/n^6\to\infty$ as $n\to\infty$. (The precise statement is by
analogy with Theorems 1, 4, and 5.)

\subsection{Organization of the paper}

The remainder of the paper is organized as follows: In the next
section, we present the variational approach to the limit surface of
random square Young tableaux, based on the hook formula of
Frame-Thrall-Robinson. The level curves of $L$ appear as minimizers of
a certain functional. This leads to a proof of Theorem 1 in the
interior of the square, except for the explicit identification of
$L$. Section 3 is dedicated to the derivation of the explicit formula
for the minimizer. Unlike the analogous results of Logan-Shepp,
Vershik-Kerov and Cohn-Larsen-Propp, we will show that there is no
need to guess the minimizer, by giving a technique for systematically
deriving it using an inversion formula for Hilbert transforms on a
finite interval. This may prove useful in similar problems.

In section 4, we complete the proof of Theorem 1, treating the more
delicate case of the boundary of the square, and prove Theorem 3. In
section 5, we discuss the \emph{hook walk} of Greene-Nijenhuis-Wilf
and the concept of the \emph{co-transition measure} of a Young
diagram. Using the explicit formulas for the co-transition measure
derived in \cite{romik}, we compute the co-transition measure of
the level curves $g_\alpha$, proving Theorem 2. In section 6, we prove
Theorem 4. In section 7 we give the computations necessary for
settling the rectangular case. In section 8, we discuss the
connections of our results to the theory of Plancherel-random
partitions, and some open problems.

\section{A variational problem for random square tableaux}

\subsection{A large-deviation principle}
One may consider a tableau $T\in{\cal T}_n$ as a path in the
\emph{Young graph} of all Young diagrams, starting with the empty
diagram, and leading up to the $n\times n$ square diagram, where each
step is of adding one box to the diagram. Identify $T$ with this
sequence $\lambda_T^0=\phi \subset \lambda_T^1 \subset \lambda_T^2
\subset \ldots \subset \lambda_T^{n^2} = \square_n$ of
diagrams. ($\lambda_T^k$ is simply the sub-diagram of the square
comprised of those boxes where the value of the entry of $T$ is $\le
k$.) Theorem 1 is then roughly equivalent, in a sense that will be
made precise later, to the statement that for each $1\le k\le n^2-1$, the
rescaled shape of $\lambda_T^k$ for a random
$T\in{\cal T}_n$ resembles the \emph{sub-level set}
$$ \{(x,y)\in[0,1]^2 : L(x,y) \le k/n^2 \} $$ of $L$, with
probability $1-o(1)$ as $n\to\infty$. It is this approach that leads
to the large-deviation principle. Namely, we can estimate the
probability that the sub-diagram $\lambda_T^k$ has a given shape:

\paragraph{Lemma 1.} For $T\in{\cal T}_n$, denote as before
$\lambda_T^0 \subset \ldots \subset \lambda_T^{n^2}$ the path in the
Young graph defined by $T$, and for each $0\le k\le n^2$, let
$\lambda_T^k : \lambda_T^k(1) \ge \lambda_T^k(2) \ge \ldots\ge
\lambda_T^k(n)$ be the lengths of the columns of $\lambda_T^k$ (some
of them may be $0$). For any Young diagram $\lambda:\lambda(1)\ge
\lambda(2)\ge \ldots\ge \lambda(n)$ whose graph lies within the $n\times n$
square, define the function $f_{\lambda}:[0,1]\to[0,1]$ by
\begin{equation}\label{eq:lambdagraph}
f_{\lambda}(x) = \frac{1}{n}\lambda(\lceil n x\rceil).
\end{equation}
(Note that this depends implicitly on $n$.) Let $0\le k \le n^2$, and
let $\alpha=k/n^2$.  Then for any given diagram $\lambda_0 \subseteq
\square_n$ with area $k$, we have
\begin{equation}\label{eq:lemma1}
\pr_n\left( T\in{\cal T}_n : \lambda_T^k =
\lambda_0 \right) =
\exp\bigg(-(1+o(1)) n^2 (I(f_{\lambda_0})+H(\alpha)+C)\bigg)
\end{equation}
as $n\to\infty$, where
\begin{eqnarray*}
 C &=& \frac{3}{2}-2\log 2, \\
 H(\alpha) &=& -\alpha\log(\alpha)-(1-\alpha)\log(1-\alpha), \\
 I(g) &=& \int_0^1 \int_0^1 \log|g(x)-y+g^{-1}(y)-x|dy\, dx, \\
 g^{-1}(y) &=& \inf\{x\in[0,1]: g(x)\le y\}.
\end{eqnarray*}
The $o(1)$ is uniform over all $\lambda_0$ and all $0\le k\le n^2$.

\paragraph{Proof.} For a Young diagram $\lambda:\lambda(1)\ge
 \lambda(2)\ge \ldots \ge \lambda(l)$ of area $m$, denote by $d(\lambda)$
the number of Young tableaux of shape $\lambda$ (also known as the
\emph{dimension} of $\lambda$, as it is known to be equal to the
dimension of a certain irreducible representation corresponding to
$\lambda$ of the symmetric group of order $m$). Recall the \emph{hook
formula} of Frame-Thrall-Robinson \cite{framehook}, which says that
$d(\lambda)$ is given by
\begin{equation}\label{eq:hook}
 d(\lambda) = \frac{m!}{\prod_{(i,j)\in \lambda} h_{i,j}},
\end{equation}
where the product is over all boxes $(i,j)$ in the diagram, and
$h_{i,j}$ is the \emph{hook number} of a box, given by
\begin{eqnarray*}
h_{i,j} &=& \lambda(i) - j + \lambda'(j) - i + 1\\ &=& 1+\textrm{number of
boxes either to the right of, or below $(i,j)$}
\end{eqnarray*}
(and where $\lambda'$ is the conjugate partition to $\lambda$.)
Then we have
\footnote{Note to the reader: this is probably the most important
formula in the paper!}
\begin{equation}\label{eq:dmeasure}
\pr_n\left( T\in{\cal T}_n : \lambda_T^k =
\lambda_0 \right) = \frac{d(\lambda_0)
d(\square_n\setminus \lambda_0)}{d(\square_n)},
\end{equation}
where $d(\square_n\setminus \lambda_0)$ means the number of fillings
of the numbers $1,\ldots,n^2-k$ in the cells of the skew-Young diagram
$\square_n\setminus \lambda_0$ that are monotonically
\emph{decreasing} along rows and columns. This is because
$\square_n\setminus\lambda_0$ can be thought of as an ordinary
diagram, when viewed from the opposite corner of the square. The
number of square tableaux whose $k$-th subtableau has shape
$\lambda_0$ is simply the number of tableaux of shape $\lambda_0$,
times the number of fillings of the numbers $k+1,k+2,\ldots,n^2$ in the
cells of $\square_n\setminus\lambda_0$ that are monotonically
increasing along rows and columns -- and these are of course
isomorphic to tableaux of shape $\square_n\setminus\lambda_0$, by
replacing each entry $i$ with $n^2+1-i$.

Take minus the logarithm of \eqref{eq:dmeasure} and divide by $n^2$,
using \eqref{eq:hook}. The right-hand side becomes
$$ a+b+c-d := \frac{1}{n^2}\log \left(\frac{(n^2)!}{k!(n^2-k)!}\right)
+ \frac{1}{n^2}\sum_{i=1}^n \sum_{j=1}^{\lambda(i)} \log(\lambda(i) -
j + \lambda'(j)-i+1) $$ $$ + \frac{1}{n^2}\sum_{i=1}^n
\sum_{j=\lambda(i)+1}^n \log(j-\lambda(i)+i-\lambda'(j)+1) -
\frac{1}{n^2} \sum_{i=1}^n \sum_{j=1}^n \log(2 n - i - j + 1).$$ 
By Stirling's formula, we have $ a = n^{-2} \log\binom{n^2}{k} =
H(\alpha)+o(1)$, with the required uniformity in $k$. The other
summands look like Riemann sums of double integrals. Indeed, we claim
that
\begin{eqnarray*}
b &=& \int_0^1 \int_0^{f_{\lambda_0}(x)} \log\bigg(
  f_{\lambda_0}(x)-y+f_{\lambda_0}^{-1}(y)-x \bigg) dy\,dx +
\frac{k}{n^2}\log n + o(1), \\
c &=& \int_0^1 \int_{f_{\lambda_0}(x)}^1 \log\bigg(
  y-f_{\lambda_0}(x)+x-f_{\lambda_0}^{-1}(y)\bigg) dy\,dx +
\frac{n^2-k}{n^2} \log n + o(1), \\
d &=& \int_0^1 \int_0^1 \log(2-x-y)dy\,dx + \log n + o(1)
= C + \log n +o(1),
\end{eqnarray*}
which on summing and exponentiating would give the lemma. Let us
prove, for example, the first of these equations. Write
\begin{eqnarray*}
b &=& \frac{1}{n^2}\sum_{i=1}^n
\sum_{j=1}^{\lambda(i)} \log(\lambda(i) - j + \lambda'(j)-i+1) \\ &=&
\frac{1}{n^2}\sum_{i=1}^n
\sum_{j=1}^{\lambda(i)} \log\left(\frac{\lambda(i) - j +
\lambda'(j)-i+1}{n}\right) + \frac{k}{n^2}\log n.
\end{eqnarray*}
Fix $1\le i\le n$ and $1\le j \le \lambda(i)$. Denote $h =
(\lambda(i) - j + \lambda'(j)-i+1)/n$. Approximate $n^{-2} \log h$ in
the above sum by the double integral
$$ Q := \int_{(i-1)/n}^{i/n} \int_{(j-1)/n}^{j/n} \log \left(
f_{\lambda_0}(x)-y + f_{\lambda_0}^{-1}(y)-x \right)dy\,dx.
$$
A change of variables transforms this (check the definition of
$f_{\lambda_0}$) into
$$ Q = \int_{-1/2n}^{1/2n} \int_{-1/2n}^{1/2n} \log (x+y+h) dx\, dy.
$$
%
%
%
%
%
Note that $h$ may take the values $1/n, 2/n, \ldots, (2n-1)/n$. If
$h=1/n$, then integrating we get
$$ Q = -\frac{\log n}{n^2} + n^{-2} \int_0^1 \int_0^1 \log(u+v)du\,dv
= \frac{\log h}{n^2} + O(n^{-2}). $$
If $h\ge 2/n$, by the integral mean value theorem, we have for some
$\eta \in [-1,1]$,
$$ Q = \frac{\log(h+\eta n^{-1})}{n^2} = \frac{\log h}{n^2} + O((n^3
h)^{-1}). $$
Clearly then the last estimate holds for $h=1/n$ as well. The sum of
the remainders over all $1 \le i \le n$, $1 \le j \le \lambda(i)$ is
of order
$$ n^{-2}\sum_{(i,j)\in \lambda_0} \frac{1}{h_{i,j}} \le n^{-2}
\sum_{m=1}^{2n-1} \frac{a(m)}{m}, $$
where
$$ a(m) := \#\{ (i,j)\in \lambda_0 : h_{i,j}=m \}. $$ 
Clearly $a(m) \le n$, since each row $i$ of $\lambda_0$ contains at
most one cell $(i,j)$ with $h_{i,j}=m$. This gives that the sum of the
remainders is of order
$$ n^{-2}\sum_{m=1}^{2n-1} \frac{n}{m} = O\left(\frac{\log
n}{n}\right),$$
which is indeed $o(1)$.
\qed

\subsection{Two formulations of the variational problem}

Lemma 1 says, roughly, that the exponential order of the probability
that a random square tableau $T$ has a given $k$-subtableau shape,
where $k$ is approximately $\alpha\cdot n^2$, is given by the value of
the functional $I$ on the boundary $g$ of the shape, plus some terms
depending only on $\alpha$.  Following the well-known methodology of
large deviation theory, the natural next step is to identify the
global minimum of $I$ over the appropriate class of functions, or in
other words to find the \emph{most likely} shape for the
$\alpha$-level set. If we can prove that there is a unique minimum,
and identify it, that will be a major step towards proving Theorem 1.
So we have arrived at the following variational problem.

\paragraph{Variational problem 1.} For each $0<\alpha<1$, any weakly
decreasing function $f:[0,1]\to[0,1]$ such that $\int_0^1
f(x)dx=\alpha$ is called \emph{$\alpha$-admissible}. Find the unique
$\alpha$-admissible function that minimizes the functional
$$ I(f) = \int_0^1 \int_0^1 \log|f(x)-y+f^{-1}(y)-x|dy\, dx. $$

\bigskip
We now simplify the form of the functional $I$, by first rotating the
coordinate axes by 45 degrees, and then reparametrizing the square by
the ``hook coordinates'' -- an idea used in \cite{vershikkerov1},
\cite{vershikkerov2}, \cite{loganshepp}. Let $u, v$ be the rotated
coordinates as in \eqref{eq:rotated}. Given an $\alpha$-admissible
function $f:[0,1]\to[0,1]$, there corresponds to it a function
$g:[-\sqrt{2}/2, \sqrt{2}/2]\to[0,\sqrt{2}]$, such that
$$ y = f(x) \iff v = g(u) $$
(see Figure 4). The class of $\alpha$-admissible functions translates
to those functions $g:[-\sqrt{2}/2,\sqrt{2}/2]\to[0,\sqrt{2}]$ that
are $1$-Lipschitz, and satisfy
$g(-\sqrt{2}/2)=g(\sqrt{2}/2)=\sqrt{2}/2 $ and
\begin{equation}\label{eq:alphacond}
\int_{-\sqrt{2}/2}^{\sqrt{2}/2} (g(u)-|u|)du = \alpha.
\end{equation}
We continue to call such functions $\alpha$-admissible. We call a
function admissible if it is $\alpha$-admissible for some $0\le \alpha
\le 1$.

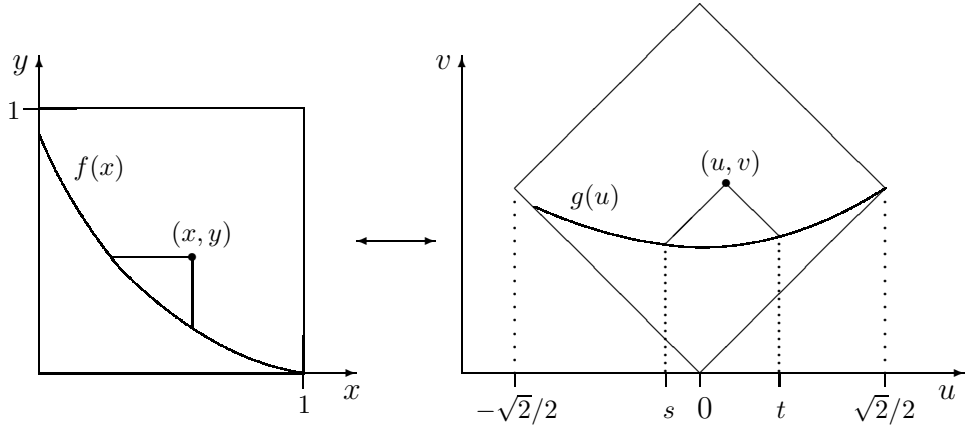
\begin{figure}[h!]
\begin{center}
\begin{picture}(500,160)(0,0)
\put(20,20){\vector(1,0){120}}
\put(20,20){\vector(0,1){120}}
\put(20,20){\framebox(100,100)}
\put(10,135){$y$}
\put(14,120){\line(1,0){20}}
\put(8,116){\footnotesize 1}
\put(135,10){$x$}
\put(120,14){\line(0,1){20}}
\put(118,5){\footnotesize 1}
\put(140,70){\vector(1,0){30}}
\put(170,70){\vector(-1,0){30}}
\put(180,20){\vector(1,0){190}}
\put(270,20){\line(-1,1){70}}
\put(270,20){\line(1,1){70}}
\put(200,90){\line(1,1){70}}
\put(340,90){\line(-1,1){70}}
\multiput(340,20)(0,5){14}{\circle*{1}}
\multiput(200,20)(0,5){14}{\circle*{1}}
\put(180,20){\vector(0,1){120}}
\put(360,10){$u$}
\put(170,135){$v$}
\put(270,20){\line(0,-1){5}}
\put(269,3){0}
\put(185,3){\footnotesize $-\sqrt{2}/2$}
\put(328,3){\footnotesize $\sqrt{2}/2$}
\put(200,20){\line(0,-1){5}}
\put(340,20){\line(0,-1){5}}
\qbezier(20,110)(30,85)(50,60)
\qbezier(50,60)(85,25)(120,20)
\qbezier(207,83)(232,72)(254,69)
\qbezier(254,69)(296,62)(340,90)
\put(280,92){\circle*{3}}
\put(270,97){\footnotesize $(u,v)$}
\put(78,64){\circle*{3}}
\put(70,70){\footnotesize $(x,y)$}
\put(280,92){\line(-1,-1){23}}
\put(280,92){\line(1,-1){20}}
\put(78,64){\line(-1,0){30}}
\put(78,64){\line(0,-1){27}}
\multiput(257,69)(0,-3){17}{\circle*{1}}
\multiput(300,72)(0,-3){18}{\circle*{1}}
\put(300,20){\line(0,-1){5}}
\put(257,20){\line(0,-1){5}}
\put(299,3){\footnotesize $t$}
\put(256,3){\footnotesize $s$}
\put(33,95){\footnotesize $f(x)$}
\put(221,84){\footnotesize $g(u)$}
\end{picture}
\caption{The rotated graph and the hook coordinates $s,t$}
\end{center}
\end{figure}

To derive the new form of the functional, write
$$ I(f) = I_1(f) + I_2(f) := \int_0^1 \int_0^{f(x)}
\log(h_f(x,y))dy\,dx + \int_0^1 \int_{f(x)}^1 \log(h_f(x,y))dy\,dx, $$
where $h_f(x,y)$ is the \emph{hook function} of $f$,
$$ h_f(x,y) = |f(x)-y + f^{-1}(y)-x|. $$
Now, set $$J(g) = J_1(g) + J_2(g) := I_1(f) + I_2(f),$$
where $f$ and $g$ are rotated versions of the same graph as in Figure
4. Then
$$ J_2(g) = \int_{-\sqrt{2}/2}^{\sqrt{2}/2} \int_{g(u)}^{\sqrt{2}-|u|}
  \log h_f(x,y) dv\,du. $$
Reparametrize this double integral by the \emph{hook coordinates} $s$
and $t$,
$$ s = \frac{f^{-1}(y)-y}{\sqrt{2}}, \qquad t =
\frac{x-f(x)}{\sqrt{2}} $$
(see Figure 4). The Lipschitz property ensures that this
transformation is one-to-one from the region
$$ \{ (u,v) : -\sqrt{2}/2 \le u \le \sqrt{2}/2,\ g(u)\le v \le 
\sqrt{2}-|u| \} $$
onto the region
$$ \Delta = \{(s,t) : -\sqrt{2}/2 \le s \le t \le \sqrt{2}/2 \}. $$
Therefore the integral transforms as
$$ J_2(f) = \int\!\!\int_\Delta \log\left(\sqrt{2}(t-s)\right)
 \left|\frac{\partial(u,v)}{\partial(s,t)}\right|ds\,dt.
$$
It remains to compute the Jacobian $\partial(u,v)/\partial(s,t)$. An
easy computation gives (see \cite{vershikkerov1},
\cite{vershikkerov2}, \cite{loganshepp})
$$ \frac{\partial(u,v)}{\partial(s,t)} =
\frac{1}{2}(1-g'(s))(1+g'(t)). $$
(This can be viewed geometrically as follows: draw on the $u$-axis in
Figure 4 the two intervals $[s,s+ds], [t,t+dt]$. The set of points in
the square for which the hook coordinates fall inside the two
intervals is approximately a parallelogram whose area is clearly seen
from the picture to be linear in $1-g'(s)$ and in $1+g'(t)$.)
So
$$ J_2(g) = \frac{1}{2}\int\!\!\int_\Delta
\log\left(\sqrt{2}(t-s)\right) (1-g'(s))(1+g'(t)) ds\,dt. $$
A similar computation for $J_1$, using ``lower'' instead of ``upper''
hook coordinates, shows that
$$ J_1(g) = \frac{1}{2}\int\!\!\int_\Delta
\log\left(\sqrt{2}(t-s)\right) (1+g'(s))(1-g'(t)) ds\,dt. $$
This gives
\begin{eqnarray*}
 J(g)&=& \frac{1}{2}\int\!\!\int_\Delta \log\left(\sqrt{2}(t-s)\right)
\big[(1-g'(s))(1+g'(t))+(1+g'(s))(1-g'(t))\big]ds\,dt  \\
&=& \frac{1}{2}
\int\!\!\int_\Delta
\log\left(\sqrt{2}(t-s)\right)\big(2-2g'(s)g'(t)\big)ds\,dt  \\
&=& -\frac{1}{2}
\int_{-\sqrt{2}/2}^{\sqrt{2}/2}\int_{-\sqrt{2}/2}^{\sqrt{2}/2} 
\log|t-s|\cdot g'(s)g'(t)ds\,dt + \log 2-\frac{3}{2}.
\end{eqnarray*}
We can now state a reformulation of the original variational problem.

\paragraph{Variational problem 2.} For each $0<\alpha<1$, a function
$g:[-\sqrt{2}/2,\sqrt{2}/2]\to[0,\sqrt{2}]$ is called
$\alpha$-admissible if: $g(-\sqrt{2}/2)=g(\sqrt{2}/2)=\sqrt{2}/2$; $g$
is 1-Lipschitz; and $\int_{-\sqrt{2}/2}^{\sqrt{2}/2}
(g(u)-|u|)du = \alpha$. Find the unique $\alpha$-admissible function
that minimizes the functional
\begin{equation}\label{eq:functionalk}
 K(g) = -\frac{1}{2}\int_{\sqrt{2}/2}^{\sqrt{2}/2}
\int_{-\sqrt{2}/2}^{\sqrt{2}/2} g'(s)g'(t) \log|s-t| ds dt.
\end{equation}

\subsection{Deduction of Theorem 1(ii)}

In the next section, we prove the following theorem.

\paragraph{Theorem 7.} For each $0<\alpha<1$, let $\tilde{g}_\alpha$
be the unique extension of $g_\alpha$ (defined in
\eqref{eq:minimizers}) to an $\alpha$-admissible function, namely
$$ \tilde{g}_\alpha(u) = \left\{ \begin{array}{ll} 
   g_\alpha(u) & |u|\le \sqrt{2\alpha(1-\alpha)} \\
   |u| & \sqrt{2\alpha(1-\alpha)} \le |u| \le \sqrt{2}/2 \end{array}
\right. $$
for $0<\alpha\le 1/2$, and
$$ \tilde{g}_\alpha(u) = \left\{ \begin{array}{ll} g_\alpha(u) &
|u|\le \sqrt{2\alpha(1-\alpha)} \\ \sqrt{2}-|u| &
\sqrt{2\alpha(1-\alpha)} \le |u| \le \sqrt{2}/2 \end{array} \right. $$
for $1/2 < \alpha < 1$.
Then:

(i) $\tilde{g}_\alpha$ is the unique solution to Variational
problem 2;

(ii) $K(\tilde{g}_\alpha) = -H(\alpha)+\log 2$;

(iii) For any $\alpha$-admissible function $g$ we have
$$ K(g) \ge K(\tilde{g}_\alpha) + K(g-\tilde{g}_\alpha). $$

\bigskip
Assuming this as proven, our goal is now to prove Theorem 1. At the
beginning of this section, we claimed that Theorem 1 was equivalent to
the statement that the subtableau $\lambda_T^k$ has shape
approximately described by the region bounded under the graph of the
level curve $\{L=k/n^2\}$ (which in rotated coordinates is given by
the curve $v=\tilde{g}_\alpha(u)$, where $\alpha=k/n^2$). We shall now
make precise the sense in which this is true, and see how this follows
from the fact that $\tilde{g}_\alpha$ is the minimizer.

For a continuous function $p:[-\sqrt{2}/2,\sqrt{2}/2]\to\mathbb{R}$,
define its supremum norm
$$ ||p||_\infty = \max_{u\in[-\sqrt{2}/2,\sqrt{2}/2]}
|p(u)|. $$

\paragraph{Lemma 2.} $K$ is continuous in the supremum norm on the
space of admissible functions.

\paragraph{Proof.} Consider the symmetric bilinear form
\begin{equation}\label{eq:bilinear}
\langle g,h \rangle = -\frac{1}{2}\int_{-\sqrt{2}/2}^{\sqrt{2}/2}
\int_{-\sqrt{2}/2}^{\sqrt{2}/2} g'(s)h'(t) \log|s-t|ds\,dt
\end{equation}
defined whenever $g$ and $h$ are almost everywhere differentiable
functions on $[-\sqrt{2}/2,\sqrt{2}/2]$ with bounded derivative. We
show that $\langle \cdot,\cdot\rangle$ is continuous in the supremum
norm with respect to any of its arguments, when restricted to the set
of $1$-Lipschitz functions; this will imply the lemma, since $K(g) =
\langle g,g\rangle$.
%
%
Write \eqref{eq:bilinear} more carefully as
\begin{equation*}
\langle g,h \rangle = -\frac{1}{2}\int_{-\sqrt{2}/2}^{\sqrt{2}/2}
g'(s) \cdot \lim_{\epsilon\searrow 0} \left[
\int_{-\sqrt{2}/2}^{s-\epsilon} h'(t) \log(s-t)dt +
\int_{s+\epsilon}^{\sqrt{2}/2} h'(t) \log(t-s)dt \right] ds.
\end{equation*}
For $s\in (-\sqrt{2}/2,\sqrt{2}/2)$ which is a point of
differentiability of $h$, integration by parts gives
\begin{multline*}
\int_{-\sqrt{2}/2}^{s-\epsilon} h'(t) \log(s-t)dt +
\int_{s+\epsilon}^{\sqrt{2}/2} h'(t) \log(t-s)dt = \\ =
h(t)\log(s-t)\bigg|_{t=-\sqrt{2}/2}^{t=s-\epsilon} - 
\int_{-\sqrt{2}/2}^{s-\epsilon}\frac{h(t)}{t-s}dt +
h(t)\log(t-s)\bigg|_{t=s+\epsilon}^{t=\sqrt{2}/2} -
\int_{s+\epsilon}^{\sqrt{2}/2} \frac{h(t)}{t-s}dt \\
= h\left(\frac{\sqrt{2}}{2}\right)\log\left(\frac{\sqrt{2}}{2}-s\right)
-h\left(-\frac{\sqrt{2}}{2}\right)\log\left(\frac{\sqrt{2}}{2}+s\right) 
\qquad\qquad\qquad\qquad\qquad\qquad\ \ 
\\ +(h(s-\epsilon)-h(s+\epsilon))\log\epsilon -
\int_{[-\sqrt{2}/2,s-\epsilon]\cup [s+\epsilon,\sqrt{2}/2]}
\frac{h(t)}{t-s}dt
\qquad\qquad\qquad\ \ \ 
\\ \xrightarrow[\ \ \ \epsilon \searrow 0\ \ \ ]{} 
h\left(\frac{\sqrt{2}}{2}\right)\log\left(\frac{\sqrt{2}}{2}-s\right)
-h\left(-\frac{\sqrt{2}}{2}\right)\log\left(\frac{\sqrt{2}}{2}+s\right)
- \pi\tilde{h}(s),
\qquad\quad\ \ 
\end{multline*}
where $\tilde{h}$ is the Hilbert transform of $h$, defined by the
principal value integral
$$ \tilde{h}(s) = \frac{1}{\pi}\int_{\mathbb{R}} \frac{h(t)}{t-s}dt $$
(think of $h$ as a function on $\mathbb{R}$ which is $0$ outside
$[-\sqrt{2}/2,\sqrt{2}/2]$.) Going back to \eqref{eq:bilinear}, this gives

\begin{multline}\label{eq:scalarprod}
\qquad \langle g,h\rangle = -\frac{1}{2}h\left(\frac{\sqrt{2}}{2}\right)
\int_{-\sqrt{2}/2}^{\sqrt{2}/2}
g'(s)\log\left(\frac{\sqrt{2}}{2}-s\right) ds \\ +
\frac{1}{2}h\left(\frac{-\sqrt{2}}{2}\right)
\int_{-\sqrt{2}/2}^{\sqrt{2}/2}
g'(s)\log\left(\frac{\sqrt{2}}{2}+s\right) ds
\qquad\qquad\qquad \  \\ +
\frac{\pi}{2} \int_{-\sqrt{2}/2}^{\sqrt{2}/2} g'(s)\tilde{h}(s)ds.
\qquad\qquad\qquad\qquad\qquad\qquad\qquad\qquad\qquad\quad\ 
\end{multline}
Now recalling that the Hilbert transform is an isometry on
$L_2(\mathbb{R})$ (see \cite{titchmarsh}, Theorem 90), and using the
fact that
$$ \bigg| \int_{-\sqrt{2}/2}^{\sqrt{2}/2} \log\left( \frac{\sqrt{2}}{2}\pm
s\right) ds \bigg| = \frac{2-\log 2}{\sqrt{2}} < 1,$$
this implies that for 1-Lipschitz functions $g,h_1,h_2$,
\begin{eqnarray*}
|\langle g,h_1-h_2\rangle| &\le& ||h_1-h_2||_\infty +
\frac{\pi}{2} \int_{-\sqrt{2}/2}^{\sqrt{2}/2}
|\tilde{h}_1(s)-\tilde{h}_2(s)|ds \\ &\le& ||h_1-h_2||_\infty +
2^{1/4}\frac{\pi}{2} \left(
\int_{-\sqrt{2}/2}^{\sqrt{2}/2}
\left(\tilde{h}_1(s)-\tilde{h}_2(s)\right)^2ds
\right)^{1/2} \\ &\le& ||h_1-h_2||_\infty +
 2^{1/4}\frac{\pi}{2} \left( \int_{\mathbb{R}}
\left(\tilde{h}_1(s)-\tilde{h}_2(s)\right)^2ds \right)^{1/2} \\ &=&
||h_1-h_2||_\infty + 2^{1/4}\frac{\pi}{2}
\left(\int_{-\sqrt{2}/2}^{\sqrt{2}/2}(h_1(s)-h_2(s))^2ds
\right)^{1/2} \\ &\le&
\left(1+2^{1/2} \frac{\pi}{2}\right) ||h_1-h_2||_\infty.
\end{eqnarray*}

\vspace{-26.0pt}
\qed

\bigskip\medskip
We have another use for \eqref{eq:scalarprod}. Let $f$ be a Lipschitz
function on $[-\sqrt{2}/2,\sqrt{2}/2]$ that satisfies $f(\pm
\sqrt{2}/2)=0$. Denote by
$$ F[f](x) = \int_{\mathbb{R}} f(t)e^{-ixt}dt $$
the Fourier transform of a function $f$. Recall the well-known
formulas
\begin{eqnarray*}
F[\tilde{f}](x) &=& i\cdot \sgn x\cdot F[f](x),\\
F[f'](x) &=& i\cdot x \cdot F[f](x),\\
\int_{\mathbb{R}} f_1(t)\overline{f_2(t)}dt &=& \frac{1}{2\pi}
\int_{\mathbb{R}}F[f_1](x) \overline{F[f_2](x)} dx.
\end{eqnarray*}
Then, by \eqref{eq:scalarprod}
\begin{multline}\label{eq:fourierk}
\qquad\qquad\ \,\,\,
K(f) = \langle f,f\rangle = 
\frac{\pi}{2} \int_{-\sqrt{2}/2}^{\sqrt{2}/2} f'(s)\tilde{f}(s)ds 
\\ =
\frac{1}{4}\int_{\mathbb{R}} F[f'](x)\overline{F[\tilde{f}](x)}dx =
\frac{1}{4}\int_{\mathbb{R}} |x|\cdot|F[f](x)|^2dx.
\qquad\qquad\qquad
\end{multline}
We note as a lemma an important consequence of this identity which we
shall need later on.

\paragraph{Lemma 3.} If $f$ is a Lipschitz function with $f(\pm
\sqrt{2}/2)=0$ as above, then $K(f)\ge 0$, and $K(f)=0$ only if
$f\equiv 0$.  \qed

\bigskip
Lemma 3 will be used in the next section to easily deduce uniqueness
of the minimizer. In fact, Theorem 7 gives all the necessary
information to prove a non-quantitative version of Theorem 1,
i.e. without the rate-of-convergence estimates. However, we can do
better, by noting that Theorem 7(iii), together with the
representation \eqref{eq:fourierk}, can be used to give
quantitative estimates for the rate of convergence in Theorem 1. We
prove the following strengthening of Lemma 3:

\paragraph{Lemma 4.} For every $r\in(2,3)$, 
there exists a constant $c=c(r) > 0$ such that for all 2-Lipschitz
functions $f:[-\sqrt{2}/2,\sqrt{2}/2]\to \mathbb{R}$ that satisfy
$f(\pm \sqrt{2}/2) = 0$, we have
$$ K(f) \ge c ||f||_\infty^r. $$

\paragraph{Proof.} Had the power of $|x|$ in \eqref{eq:fourierk} been
2, $K(f)$ would have been equal to $1/4$ times the squared $L_2$-norm
of $x F[f](x) = F[f'](x)$. Having $|x|$ in
\eqref{eq:fourierk} invites the conclusion that instead we are dealing
with the squared $L_2$-norm of $f^{(1/2)}(x)$, the \emph{fractional}
derivative of $f$ of order $1/2$.

To see that this is indeed the case, and to use the full power of such
an interpretation of $K(f)$, let us recall the corresponding
definitions. For $\alpha \in (0,1)$, the fractional derivative
$f^{(\alpha)}(x)$ of order $\alpha$ is defined by
\begin{equation}\label{eq:fractional}
 f^{(\alpha)}(x) = \frac{\alpha}{\Gamma(1-\alpha)} \int_0^\infty
 \frac{f(x)-f(x-t)}{t^{1+\alpha}} dt.
\end{equation}
The integral exists as $f(x)$ is Lipschitz and bounded. Clearly
$f^{(\alpha)}(x)\equiv 0$ for $x\le -\sqrt{2}/2$. Then
\begin{multline}\label{eq:b4}
\qquad\qquad
F[f^{(\alpha)}](x) = \int_\mathbb{R} e^{-ixt}f^{(\alpha)}(t)dt
\\  = \frac{\alpha}{\Gamma(1-\alpha)} \int_0^\infty
 \frac{1-e^{-ix\tau}}{\tau^{1+\alpha}} d\tau \cdot F[f](x)
= (ix)^{\alpha} F[f](x), \qquad\ \ 
\end{multline}
where
$$
(ix)^\alpha := \left\{ \begin{array}{ll}
  |x|^\alpha \exp(i\alpha\pi/2), & x>0, \\
  |x|^\alpha \exp(-i\alpha\pi/2),& x<0. \end{array}\right.
$$
Indeed, setting
$$ z^{1+\alpha} = |z|\exp(i(1+\alpha)\theta)),
\quad \textrm{if }z=|z|e^{i\theta},\ \ \theta\in(-\pi,\pi), $$
we have
\begin{eqnarray*}
\int_0^\infty \frac{1-e^{-ix\tau}}{\tau^{1+\alpha}}d\tau &=&
(ix)^\alpha \int_0^{i\infty} \frac{1-e^{-z}}{z^{1+\alpha}}dz =
(ix)^\alpha \int_0^\infty \frac{1-e^{-\tau}}{\tau^{1+\alpha}}d\tau \\
&=& (ix)^\alpha \frac{1}{\alpha} \int_0^\infty
\tau^{-\alpha}e^{-\tau}d\tau = (ix)^\alpha
\frac{\Gamma(1-\alpha)}{\alpha}.
\end{eqnarray*}
In particular, for $\alpha=1/2$, we get from \eqref{eq:b4} that
$$ \big|F[f^{(1/2)}](x)\big|^2 = |x|\cdot |F[f](x)|^2, $$
whence, by \eqref{eq:b4} and isometry of the Fourier transform,
\begin{equation} \label{eq:b5}
K(f) = \frac{1}{4}\int_\mathbb{R} |x|\cdot |F[f](x)|^2dx =
\frac{\pi}{2}|f^{(1/2)}(x)|^2.
\end{equation}
The fractional integration operator, inverse to that in
\eqref{eq:fractional}, is known to be given by
\begin{equation}\label{eq:b6}
f(x) = (I_\alpha f^{(\alpha)})(x), \qquad
(I_\alpha h)(x) := \frac{1}{\Gamma(\alpha)}\int_{-\infty}^x
(x-t)^{\alpha-1} h(t)dt.
\end{equation}
As a check, the Fourier transform of the RHS is
$$ \frac{1}{\Gamma(\alpha)} F[f^{(\alpha)}](x) \int_0^\infty
\tau^{\alpha-1} e^{-ix\tau}d\tau = (ix)^{-\alpha} F[f^{(\alpha)}](x) =
F[f](x). $$
By Theorem 383 in \cite{hardyetal}, for $p>1$ and
$$ 0 < \alpha < \frac{1}{p}, \qquad q=\frac{p}{1-\alpha p}, $$
$I_\alpha$ maps $L_p$ into $L_q$, and is bounded. That is, there
exists a constant $c(p)>0$ such that
\begin{equation}\label{eq:b7}
||I_\alpha h||_q \le c(p) ||h||_p.
\end{equation}
Introduce
$\psi(x)=f^{(\alpha)}(x)\mathbf{1}_{(-\infty,\sqrt{2}/2]}(x)$, so that
$\psi$ is supported by $[-\sqrt{2}/2,\sqrt{2}/2]$. According to
\eqref{eq:b6},
$$ (I_\alpha \psi)(x) = f(x), \qquad x\le \sqrt{2}/2. $$
So, using \eqref{eq:b7} and monotonicity of the $L_s$-averages, we
have
$$
||f||_q \le ||I_\alpha \psi ||_q \le c(p) ||\psi||_p \le c_1(p)
||\psi||_2 \le c_2(p) ||f^{(\alpha)}||_2, \quad c_1(p) :=
(\sqrt{2})^{1/p-1/2} c(p).
$$
In light of \eqref{eq:b5}, for $\alpha=1/2$ we obtain then
\begin{equation}\label{eq:b8}
||f||_q^2 \le c_2(p) K(f), \quad c_2(p):= \frac{2}{\pi}c_1(p)^2,
\quad \left(p\in (1,2),\ q=\frac{p}{1-p/2} \right).
\end{equation}
Let $x_0 \in (-\sqrt{2}/2,\sqrt{2}/2)$ be such that
$|f(x_0)|=||f||_\infty$. Since $f$ is $2$-Lipschitz,
$$ |f(x)| \ge ||f||_\infty - 2|x-x_0|, \quad |x-x_0| \le
\frac{||f||_\infty}{2}. $$
Then
$$
||f||_q^2 \ge \left(2 \int_0^{||f||_\infty/2} (||f||_\infty-2y)^qdy
  \right)^{2/q} = \frac{||f||_\infty^{2(q+1)/q}}{(q+1)^{2/q}},
$$
so, using \eqref{eq:b8}, we conclude that, for an absolute constant
$c^*(p,q)>0$,
$$ K(f) \ge c^*(p) ||f||_\infty^{2(q+1)/q}. $$
It remains to observe that
$$ \frac{2(q+1)}{q}=1+\frac{2}{p} $$
can be made arbitrarily close to 2 from above by selecting $p$
sufficiently close to 2 from below. This completes the proof.
\qed

\paragraph{Theorem 8.} For a Young diagram $\lambda$ whose graph lies
within the $n\times n$ square, let $g_{\lambda}(u)$ be the rotated
coordinate version of the function $f_{\lambda}(x)$ defined in
\eqref{eq:lambdagraph}. Denote $\alpha = k/n^2$.
Then for all $2<r<3$, there are constants $c=c(r)>0,
C=C(r)>0$ such that for any $\epsilon>0$ and for any $n$,
\begin{equation}\label{eq:theorem8}
\pr_n\bigg( T\in{\cal T}_n :
\max_{1\le k\le n^2-1} 
||g_{\lambda_T^{k}}-\tilde{g}_\alpha||_\infty >
\epsilon \bigg) \le C \exp(3n -c \,\epsilon^r n^2).
\end{equation}
Consequently, with probability subexponentially close to $1$, for all
$k$ the supnorm distance between $g_{\lambda_T^k} $ and
$\tilde{g}_{\alpha},\, (\alpha=k/n^2),$ does not exceed
$n^{-1/2+\delta}$, $(\delta>0)$.

\paragraph{Proof.} Let $p(m)$ be the number of partitions of an
integer $m$. It is known that for all $m$, $p(m) \le \exp(\pi
\sqrt{2m/3})$ (see \cite{apostol}, Theorem 14.5). Fix $n$, $1\le k\le
n^2-1$, $\epsilon>0$. Using \hbox{Lemma 1},
\begin{equation*}
\pr_n\bigg( T\in{\cal T}_n : ||g_{\lambda_T^{k}}
-\tilde{g}_\alpha||_\infty > \epsilon \bigg) =
\sum_{
\begin{array}{ll}\lambda_0\subseteq \square_n\textrm{ of area }k \\
||g_{\lambda_0}-\tilde{g}_\alpha||_\infty > \epsilon
\end{array}}
\pr_n\bigg( T\in{\cal T}_n : \lambda_T^{k} =
\lambda_0 \bigg)
\end{equation*}
\begin{equation}\label{eq:usinglemma}
\le p(k)\sup_{
\begin{array}{ll}\lambda_0\subseteq\square_n\textrm{ of area }k \\
||g_{\lambda_0}-\tilde{g}_\alpha||_\infty > \epsilon
\end{array}} \exp\bigg(-(1+o(1))n^2(K(g_{\lambda_0})+H(\alpha)-\log
2) \bigg).
\end{equation}
Let $\lambda_0$ be a diagram contained in $\square_n$ of area $k$,
such that $||g_{\lambda_0}-\tilde{g}_\alpha||_\infty > \epsilon$.
Since $g_{\lambda_0}$ is $\alpha$-admissible, using Theorem 7 and
Lemma 4 we have
$$ K(g_{\lambda_0})+H(\alpha)-\log 2 \ge 
K(g_{\lambda_0}-\tilde{g}_\alpha) >
c(r)||g_{\lambda_0}-\tilde{g}_\alpha||_\infty^r
\ge c(r)\epsilon^r. $$
Combining this with \eqref{eq:usinglemma} and with 
the above remark on the number of partitions of an integer gives
that for $n$ larger than some absolute initial bound,
\begin{equation*}
\pr_n\bigg( T\in{\cal T}_n : ||g_{\lambda_T^{k}}
-\tilde{g}_\alpha||_\infty > \epsilon \bigg) \le \exp(2.8 \sqrt{\alpha}n - c n^2
\epsilon^r ).
\end{equation*}
Taking the union bound over all $1\le k\le n^2-1$ gives
\eqref{eq:theorem8}.
\qed

\paragraph{Lemma 5.} For each $(x,y)\in (0,1)\times(0,1)$, let $(u,v)$
be their rotated coordinates as in \eqref{eq:rotated}. Let $\alpha_0 =
L(x,y)$, so that $|u|<\sqrt{2\alpha_0(1-\alpha_0)}$ and
$v=\tilde{g}_{\alpha_0}(u)$. There exist absolute constants $c_1,
c_2>0$ such that if we set
$$\sigma (x,y)=\min(xy,(1-x)(1-y)),$$
$$d(x,y)=c_1\sqrt{\sigma(x,y)},\quad \Delta(x,y) = c_2 \sigma^2(x,y), $$
we will have that for all $0<\alpha<1$ and $\delta < \Delta(x,y)$, if
$|\tilde{g}_\alpha(u)-\tilde{g}_{\alpha_0}(u)|<\delta\cdot d(x,y)$ then
$|\alpha-\alpha_0| < \delta$.

\paragraph{Proof.} Since $\tilde{g}_{\alpha}(u)$ increases with $\alpha$, it
suffices to prove existence of two absolute constants $\gamma_1,\gamma_2>0$ 
such that 
$$
|\tilde{g}_{\alpha}(u)-\tilde{g}_{\alpha_0}(u)|\ge \gamma_1\sigma^{1/2}(x,y)
|\alpha-\alpha_0|, \quad\textrm{if }|\alpha-\alpha_0|\le \gamma_2\sigma(x,y).
$$
Because of the symmetry property 
$\tilde{g}_{1-\alpha}(u)=\sqrt{2}-\tilde{g}_{\alpha}(u)$, we may
assume that $x+y\le 1$, or equivalently that $\alpha_0\le 1/2$. 


To prove the above claim, we note the following inequalities. Notice
first that
$$ \sqrt{2\alpha_0(1-\alpha_0)} \ge v \implies \alpha_0 \ge \frac{1-\sqrt{1-2
v^2}}{2}. $$
Likewise, $\alpha^{(-)}$ that corresponds to the lowest point $(u,u)$ is given
by
$$ \alpha^{(-)} = \frac{1-\sqrt{1-2 u^2}}{2}. $$
and we see that
\begin{equation}\label{eq:lowest}
\alpha_0 - \alpha^{(-)} \ge \frac{\sqrt{1-2 u^2}-\sqrt{1-2 v^2}}{2}
= \frac{v^2-u^2}{\sqrt{1-2u^2}+\sqrt{1-2 v^2}} \ge \frac{v^2-u^2}{2}
= xy.
\end{equation}
\eqref{eq:lowest} says that decreasing $\alpha_0$ by $x_0 y_0$ gives
us a feasible $\alpha$, for which $(u,\tilde{g}_\alpha(u))$ lies
between $(u,v)$ and the lowest point $(u,u)$, such that $u\le\sqrt{2\alpha
(1-\alpha)}$.

Let us estimate from above $\tilde{g}_\alpha(u)$ for
$\alpha\in[\alpha^{(-)},\alpha_0]$. From \eqref{eq:partialdiff} it
follows that
$$ \frac{\partial \tilde{g}_\alpha(u)/\partial
\alpha}{\sqrt{\tilde{g}_\alpha(u)^2-u^2}} \ge c $$
for some absolute constant $c>0$. (Indeed,
$2\alpha(1-\alpha)=\beta^2(\alpha) \ge \tilde{g}_\alpha(u)^2$.)
Integrating from $\alpha\in [\alpha^{(-)},\alpha_0]$ and
exponentiating, we obtain
$$ \frac{\tilde{g}_{\alpha_0}(u) + \sqrt{\tilde{g}_{\alpha_0}(u)^2-u^2}}
{\tilde{g}_{\alpha}(u) + \sqrt{\tilde{g}_{\alpha}(u)^2-u^2}} \ge
\exp(c(\alpha_0-\alpha)),$$
or equivalently
$$ \frac{\tilde{g}_{\alpha}(u) - \sqrt{\tilde{g}_{\alpha}(u)^2-u^2}}
{\tilde{g}_{\alpha_0}(u) - \sqrt{\tilde{g}_{\alpha_0}(u)^2-u^2}} \ge
\exp(c(\alpha_0-\alpha)).$$
Consequently
$$ \sqrt{\tilde{g}_\alpha(u)^2-\tilde{g}_{\alpha_0}(u)^2} \le
\cosh(c(\alpha_0-\alpha)) \sqrt{\tilde{g}_{\alpha_0}(u)^2-u^2} -
\sinh(c(\alpha_0-\alpha))\tilde{g}_{\alpha_0}(u), $$
or
$$ \tilde{g}_\alpha(u)^2 \le \bigg[\cosh(c(\alpha_0-\alpha))
\tilde{g}_{\alpha_0}(u) - \sinh(c(\alpha-\alpha_0))
\sqrt{\tilde{g}_{\alpha_0}(u)^2 - u^2}\bigg]^2, $$
so that
$$ \tilde{g}_{\alpha}(u) \le 
\cosh(c(\alpha_0-\alpha))
\tilde{g}_{\alpha_0}(u) - \sinh(c(\alpha-\alpha_0))
\sqrt{\tilde{g}_{\alpha_0}(u)^2 - u^2}. $$
Consequently, for some constants $c_i > 0$,
\begin{eqnarray*}
\tilde{g}_{\alpha}(u) - \tilde{g}_{\alpha_0}(u) &\le& -c_3
(\alpha_0-\alpha) [ (v^2-u^2)^{1/2} - c_4 (\alpha_0-\alpha)v ] \\
&=& - c_5(\alpha_0-\alpha)[ (x y)^{1/2} - c_6 (\alpha_0-\alpha)(x+y)]
\\ &\le& - c_7 (\alpha_0-\alpha) (x y)^{1/2},
\end{eqnarray*}
provided that
$$ \alpha_0-\alpha \le c_8 \frac{(x y)^{1/2}}{x+y}. $$
From \eqref{eq:lowest} we know that we can go below $\alpha_0$ by $xy$ at 
least.
Pick $\rho = \min(1,c_8/3)$; then the last inequality
holds for $\alpha_0-\alpha\le \rho x y$, and we have
\begin{equation}\label{eq: smallalpha}
\tilde{g}_{\alpha}(u) - \tilde{g}_{\alpha_0}(u) \le
- c_7 (\alpha_0-\alpha)(x y)^{1/2}, \qquad \alpha\in [\alpha_0 - \rho
x y, \alpha_0]. 
\end{equation}
Now for $\alpha_0\le \alpha\le 1/2$ we know that
$$ \frac{\partial \tilde{g}_{\alpha}(u)}{\partial \alpha} \ge c_9
\sqrt{v^2-u^2} = c_{10} (x y)^{1/2}, $$
so that 
\begin{equation}\label{eq: moderatealpha1}
 \tilde{g}_\alpha(u) - \tilde{g}_{\alpha_0}(u) \ge c_{10} (x
y)^{1/2} (\alpha-\alpha_0). 
\end{equation}
By symmetry, for $1/2\le\alpha\le 1-\alpha_0$,
\begin{equation}\label{eq: moderatealpha2}
\tilde{g}_{\alpha}(u)-\tilde{g}_{1-\alpha_0}\le -c_{10}((1-\alpha_0)-\alpha)
((1-x)(1-y))^{1/2},
\end{equation}
and, for $1-\alpha_0\le\alpha\le 1-\alpha_0+\rho(1-x)(1-y)$,
\begin{equation}\label{eq: largealpha}
\tilde{g}_{\alpha}(u)-\tilde{g}_{\alpha_0}(u)\ge c_7(\alpha-(1-\alpha_0))
((1-x)(1-y))^{1/2}.
\end{equation}
The inequalities \eqref{eq: smallalpha}, \eqref{eq: moderatealpha1},
\eqref{eq: moderatealpha2}, \eqref{eq: largealpha} prove the claim with
$\gamma_1=\min\{c_7,c_{10}\}$ and $\gamma_2=\rho$.
\qed

\paragraph{Proof of Theorem 1(ii).} We now prove Theorem 1(ii), the
part of Theorem 1 that deals with the interior of the square. The
treatment of the boundary of the square is more delicate and is
deferred to section 4, being essentially equivalent to Theorem 3.

Fix $1\le i, j\le n$ such that 
\begin{equation}\label{eq:interior}
\min(i j,(n-i)(n-j)) > n^{3/2+\epsilon}.
\end{equation}
Let $(u,v)$ be the rotated coordinates
corresponding to $(x,y)=(i/n,j/n)$. Let $\alpha_0 = L(i/n,j/n)$, so
that $v=\tilde{g}_{\alpha_0}(u)$ and
$|u|<\sqrt{2\alpha_0(1-\alpha_0)}$. For each tableau
$T=(t_{i,j})_{1\le i,j\le n}\in{\cal T}_n$ let $k_T = t_{i,j}$, and let
$\beta_T = k_T/n^2$. Note that $k_T$ is an integer representing the
smallest $s$ such that $\lambda_T^s$ contains the box $(i,j)$. This
implies that
$$ v \le g_{\lambda_T^{k_T}}(u) \le v+\frac{\sqrt{2}}{n}, $$ or
\begin{equation}\label{eq:sqrtn}
\left| g_{\lambda_T^{k_T}}(u)-\tilde{g}_{\alpha_0}(u) \right| \le
\frac{\sqrt{2}}{n}.
\end{equation}
Apply Lemma 5 with $(x,y)=(i/n,j/n)$ and
$\delta=n^{-(1-\epsilon)/2}$. Note that because of
\eqref{eq:interior}, for $n$ large we have $\delta < \Delta(x,y)$ as
required. Then, making use of Theorem 8 we get
\begin{eqnarray*}
\pr_n\left( T\in{\cal T}_n : \left| \frac{1}{n^2}t_{i,j}-
L\left(\frac{i}{n},\frac{j}{n}\right)
\right| > \frac{1}{n^{(1-\epsilon)/2}} \right)
&=& 
\pr_n\Big( T\in {\cal T}_n : |\beta_T - \alpha_0 | >
\frac{1}{n^{(1-\epsilon)/2}} \Big)
\end{eqnarray*}

\vspace{-20.0pt}
\begin{eqnarray*}
\textrm{(by Lemma 5)}
&\le& \pr_n\Big( T\in {\cal T}_n :
|\tilde{g}_{\beta_T}(u)-\tilde{g}_{\alpha_0}(u) | > 
\frac{d(i/n,j/n)}{n^{(1-\epsilon)/2}} \Big) \\
\textrm{(by \eqref{eq:sqrtn})}
&\le& \pr_n\Big( T\in{\cal T}_n : \left| g_{\lambda_T^{k_T
}}(u) - \tilde{g}_{\beta_T}(u) \right| > 
\frac{d(i/n,j/n)}{n^{(1-\epsilon)/2}}-\frac{\sqrt{2}}{n} \Big)
\\ \textrm{(for $n$ large enough, by \eqref{eq:interior})} &\le&
\pr_n\Big( T\in{\cal T}_n : \left| g_{\lambda_T^{k_T
}}(u) - \tilde{g}_{\beta_T}(u) \right| > 
\frac{d(i/n,j/n)}{2 n^{(1-\epsilon)/2}} \Big)
\\ \textrm{(by Theorem 8 with $r=2+\epsilon$)}
&\le& C \exp\left(3n - c n^2 \left(\frac{d(i/n,j/n)}{2
n^{(1-\epsilon)/2}}\right)^{2+\epsilon} \right) \\
\textrm{(for $n$ large, by \eqref{eq:interior})}
& \le & C' \exp(-c' n^{3/2}).
\end{eqnarray*}
Taking the union bound over all $1\le i,j\le n$ satisfying
\eqref{eq:interior} gives the result.
\qed

\section{Solution of the variational problem}

\subsection{Preliminaries}

In this section, we prove Theorem 7. We actually \emph{derive} the
explicit formula for the minimizer using methods of the calculus of
variations and the theory of singular (Cauchy-type) integral
equations. Our derivation makes only one a priori assumption (obtained
by educated guesswork and later verified by computer simulations) on
the graphical form that the minimizer would take, and so is in a sense
more systematic than the analogous treatments in the fundamental
papers \cite{loganshepp}, \cite{vershikkerov1}, \cite{vershikkerov2},
where the solutions are brilliantly guessed using the properties of
the Hilbert transform. We believe that our technique may prove useful
in the treatment of similar problems in the future.

First, observe that because of symmetry, we need only treat the
case $\alpha \le 1/2$; the mapping $g \to \sqrt{2}-g$ takes the set of
$\alpha$-admissible functions bijectively onto the set of
$(1-\alpha)$-admissible functions, and has the property that
$K(\sqrt{2}-g)=K(g)$.

Next, observe that for $\alpha=1/2$ the assertion is immediate,
because of Lemma 3.

We prove another fact that follows from general considerations, before
turning to the derivation of the minimizer.

\paragraph{Lemma 6.} For any $0<\alpha<1$, the functional $K$ has a unique
$\alpha$-admissible minimizer.

\paragraph{Proof.} The functional $K$ is continuous on the space of
$\alpha$-admissible functions, and is bounded below by Lemma 3. By the
Arzela-Ascoli theorem, the space of $\alpha$-admissible functions is
compact in the topology induced by the supremum norm (since the
admissible functions are uniformly bounded and
equicontinuous). Therefore $K$ has a minimizer. To prove that the
minimizer is unique, let $h_1$ and $h_2$ be two distinct
$\alpha$-admissible minimizers. Then $\tilde{h}=(h_1+h_2)/2$ is also
an $\alpha$-admissible function, and $g=(h_1-h_2)/2 \nequiv 0$, $g(\pm
\sqrt{2}/2) = 0$. So, using the parallelogram identity and Lemma 3,
$$ K(\tilde{h}) = \frac{1}{2}K(h_1)+\frac{1}{2}K(h_2)-K(g) < \min_{h
\textrm{ is }\alpha\textrm{-admissible}} K(h), $$
a contradiction.
\qed

\subsection{The derivation}

We now proceed with the derivation of the minimizer, which we shall
denote $h=h_\alpha$. The dependence on $\alpha$ will be suppressed
except where it is required. For the rest of this section, $\alpha$
will be a fixed value in $(0,1/2)$, unless stated otherwise.

First, note that, under the condition $h(\pm\sqrt{2}/2)=\sqrt{2}/2$,
the $\alpha$-condition $\int_{-\sqrt{2}/2}^{\sqrt{2}/2}
(h(u)-|u|)du=\alpha$ is equivalent to
\begin{equation}\label{eq:alphacondition}
-\int_{-\sqrt{2}/2}^{\sqrt{2}/2} u h'(u) du = \alpha - \frac{1}{2}.
\end{equation}
%
We now formulate a sufficient condition for $h$ to be a minimizer. It
is based on a standard recipe of the calculus of variations, the
Lagrange formalism. We form the Lagrange function
$$ {\cal L}(h,\lambda) = K(h) - \lambda
\int_{-\sqrt{2}/2}^{\sqrt{2}/2} u h'(u) du $$
and require that, for some $\lambda, h_\alpha$ be a local minimum
point of ${\cal L}(h,\lambda)$ in the convex set of functions $h$
subject to all the restrictions except the $\alpha$-condition
\eqref{eq:alphacondition}. To be sure, we ought to include into the
function a term $\lambda'$ times the integral of $h'$, since $h$ must
meet another constraint
\begin{equation}\label{eq:integral0}
\int_{-\sqrt{2}/2}^{\sqrt{2}/2} h'(u)du = 0.
\end{equation}
We chose not to, since -- in the square case -- even without this
constraint $h'(u)$ will turn out to be odd anyway. Since ${\cal
L}(h,\lambda)$ depends explicitly on $h'$ alone, we get the equations
for the sufficient condition in a simple-minded manner, by taking
partial derivatives of ${\cal L}$ with respect to $h'(s),\ s\in
(-\sqrt{2}/2,\sqrt{2}/2)$, and paying attention only to the constraint
$-1 \le h'(s) \le 1$. The resulting ``complementary slackness''
conditions are
\begin{equation}\label{eq:necessarycond}
w(s):= - \int_{-\sqrt{2}/2}^{\sqrt{2}/2} h'(t)\log|s-t|dt - \lambda s
\ \ \ \textrm{  is } \left\{
\begin{array}{ll}
= 0, & \textrm{if }-1 < h'(s) < 1, \\
\ge 0, & \textrm{if }h'(s)=-1, \\
\le 0, & \textrm{if }h'(s)=1.
\end{array}
\right.
\end{equation}

\paragraph{Lemma 7.} If $h$ is an $\alpha$-admissible function that,
for some $\lambda\in\mathbb{R}$, satisfies \eqref{eq:necessarycond}
for all $s\in(-\sqrt{2}/2,\sqrt{2}/2)$ for which $h'(s)$ is defined,
then $h$ is a minimizer.

\paragraph{Proof.} If $g$ is a $1$-Lipschitz function on
$[-\sqrt{2}/2,\sqrt{2}/2]$, then \eqref{eq:necessarycond} implies
that $(g'(s)-h'(s))w(s) \ge 0$ for all $s$ for which this is defined,
so
$$
\int_{-\sqrt{2}/2}^{\sqrt{2}/2} g'(s)w(s)ds \ge
\int_{-\sqrt{2}/2}^{\sqrt{2}/2} h'(s)w(s)ds. $$ 
If $g$ is $\alpha$-admissible, by \eqref{eq:alphacondition} this can
be written as
\begin{eqnarray*} 2 \langle h,g \rangle + \alpha - \frac{1}{2} &=& 
2 \langle h,g \rangle - \lambda \int_{-\sqrt{2}/2}^{\sqrt{2}/2} s
g'(s) ds \\ &\ge& 2 \langle h,h \rangle - \lambda
\int_{-\sqrt{2}/2}^{\sqrt{2}/2} s h'(s)ds = 2 \langle h,h \rangle +
\alpha - \frac{1}{2},
\end{eqnarray*}
which shows that
$$ \langle h,g\rangle \ge \langle h,h \rangle. $$
Therefore, by Lemma 3 applied to the function $g-h$,
$$
\langle g,g \rangle = \langle h,h \rangle + 2 \langle h, g-h \rangle +
\langle g-h, g-h \rangle \ge \langle h,h\rangle,
$$
so $h$ is a minimizer. \qed

\bigskip
We are about to prove part (i) of Theorem 7, namely that
$h=\tilde{g}_\alpha$ is the minimizer. Assuming this, note that in the
above proof we actually showed that
$$ \langle g,g \rangle \ge \langle \tilde{g}_\alpha, \tilde{g}_\alpha
\rangle + \langle g-\tilde{g}_\alpha, g-\tilde{g}_\alpha \rangle, $$
which is precisely the claim of part (iii) of Theorem 7. So it remains
to prove parts (i) and (ii).

Our challenge now is to determine an admissible $h$ that meets the
conditions \eqref{eq:necessarycond}. Now look at Figure 1(c) with your
head tilted 45 degrees to the right. Based on the shape of the level
curves, we make the following assumption: For some
$\beta=\beta(\alpha) \in (0,\sqrt{2}/2)$,
\begin{equation}\label{eq:guess}
h'(s)\ \ \ \textrm{  is }\left\{ \begin{array}{ll}
  = -1, & \textrm{if }-\sqrt{2}/2 < s < -\beta, \\
  \in (-1,1), & \textrm{if }-\beta<s<\beta, \\
  = +1, & \textrm{if }\beta<s<\sqrt{2}/2. \end{array}\right.
\end{equation}
Substituting this into \eqref{eq:necessarycond} gives that for
$-\beta<s<\beta$,
\begin{multline}\label{eq:multline}
-\int_{-\beta}^\beta h'(t)\log|s-t|dt = \lambda s -
 \int_{-\sqrt{2}/2}^\beta \log(s-t)dt + \int_\beta^{\sqrt{2}/2}
 \log(t-s)dt \\ =\lambda s +
(\sqrt{2}/2-s)\log(\sqrt{2}/2-s)-(\sqrt{2}/2+s)\log(\sqrt{2}/2+s)
\\ + (\beta+s)\log(\beta+s) - (\beta-s)\log(\beta-s)
\qquad\qquad\qquad\qquad\qquad
\end{multline}
Assume that $h'(s)$ is continuously differentiable on
$(-\beta,\beta)$. Differentiate \eqref{eq:multline}, to obtain
\begin{equation}\label{eq:diffmultline}
-\int_{-\beta}^\beta \frac{h'(t)}{s-t}dt = \lambda+ \log
 \frac{\beta^2-s^2}{\frac{1}{2}-s^2},
\end{equation}
where the left-hand side is a principal value integral. 

In the theory of integral equations this is known as an airfoil
equation. Solving it is tantamount to inverting a Hilbert transform on
a finite interval. Fortunately for us, it can be solved! The following
theorem appears in \cite{estradakanwal}, Section 3.2, p. 74. (See also
\cite{porterstirling}, Section 9.5.2.)

\paragraph{Theorem 9.} 
The general solution of the airfoil equation
$$
\frac{1}{\pi}\int_{-1}^1\frac{g(y)}{y-x}\,dx=f(x),\quad |x|<1,
$$
with the integral understood in the principal value sense, and $f(x)$ 
satisfying a H\"older condition, is given by
\begin{equation}\label{eq:airfoil}
g(x)=\frac{1}{\pi\sqrt{1-x^2}}\int_{-1}^1\frac{\sqrt{1-y^2}f(y)}{x-y}\,dy
+\frac{c}{\sqrt{1-x^2}}.
\end{equation}

\bigskip
Applying Theorem 9 to \eqref{eq:diffmultline}, we get the equation
\begin{equation}\label{eq:kanwal}
h'(s) = \frac{1}{\pi^2(\beta^2-s^2)^{1/2}} \int_{-\beta}^\beta
(\beta^2-t^2)^{1/2} \left( \lambda+\log
\frac{\beta^2-t^2}{\frac{1}{2}-t^2} \right) \frac{dt}{s-t} 
+ \frac{c}{(\beta^2-s^2)^{1/2}}.
\end{equation}
Here the integral is again in the sense of principal value, and the
equation must hold for some value of $c$. 

We evaluate the integral in \eqref{eq:kanwal}. Consider the
contribution of the $\lambda$-term first. Substituting $t=\beta \sin
x$ and later $u=\tan x/2$, we get
\begin{eqnarray*}
\int_{-\beta}^\beta \frac{(\beta^2-t^2)^{1/2}}{s-t}dt &=& \beta
\int_{-\pi/2}^{\pi/2} \frac{\cos^2 x}{s/\beta - \sin x}dx
\\ &=& \beta
\int_{-\pi/2}^{\pi/2} (s/\beta+\sin x)dx + \beta \left(1-(s/\beta)^2
\right) \int_{-\pi/2}^{\pi/2} \frac{dx}{s/\beta-\sin x}
\\ &=& \pi s + \frac{2\beta\left(1-(s/\beta)^2\right)}{s/\beta}
\int_{-1}^1 \frac{du}{u^2 - 2(\beta/s)u + 1}.
\end{eqnarray*}
For $|s|<\beta$, the denominator in the last integral has two real
roots, $u_1\in(-1,1)$ and $u_2\notin(-1,1)$. A simple computation
shows that the principal value of this integral at $u=u_1$ is zero. So
\begin{equation}\label{eq:lambdaterm}
\int_{-\beta}^\beta \frac{(\beta^2-t^2)^{1/2}}{s-t}dt = \pi s, \qquad
s\in(-\beta,\beta).
\end{equation}
Turn to the log-part of the integral in
\eqref{eq:kanwal}. Substituting $t=\tau \beta$, $s= v_1 \beta$,
$(2\beta^2)^{-1} = v_2^2$, we see that
\begin{equation}\label{eq:logterm}
\int_{-\beta}^\beta \frac{(\beta^2-t^2)^{1/2}}{s-t} \log
\frac{\beta^2-t^2}{\frac{1}{2}-t^2} dt = \beta [
I(s/\beta,\sqrt{2}/(2\beta)) - I(-s/\beta,\sqrt{2}/(2\beta)) ],
\end{equation}
where
$$ I(\xi,\gamma) = \int_{-1}^1 \frac{(1-\eta)^{1/2}}{\xi-\eta}
\log\frac{1+\eta}{\gamma+\eta} d\eta, \qquad \xi\in[-1,1], \gamma\ge 1,
$$
is evaluated in the following lemma.

\paragraph{Lemma 8.}
\begin{multline}\label{eq:evaluation}
I(\xi,\gamma) = \pi\bigg[ 1-\gamma+\sqrt{\gamma^2-1} -
\xi\log\left(\gamma+\sqrt{\gamma^2-1}\right) \\
- 2\sqrt{1-\xi^2} \arct \sqrt{
\frac{(\gamma-1)(1-\xi)}{(\gamma+1)(1+\xi)} } \ \bigg].
\qquad\qquad\qquad\qquad
\end{multline}

\paragraph{Proof.} Notice that $I(\xi,1)=0$, and, for $x>1$,
\begin{multline}\label{eq:multline2}
\frac{\partial I(\xi,x)}{\partial x} = - \int_{-1}^1
\frac{(1-\eta^2)^{1/2}}{(\xi-\eta)(x+\eta)}d\eta \\
= -\frac{1}{x+\xi}\left[ \int_{-1}^1
\frac{(1-\eta^2)^{1/2}}{\xi-\eta}d\eta + \int_{-1}^1
\frac{(1-\eta^2)^{1/2}}{x+\eta}d\eta \right] \qquad\qquad\qquad\quad\\
= -\frac{\pi\xi}{x+\xi}-\frac{1}{x+\xi}\int_{-1}^1
\frac{(1-\eta^2)^{1/2}}{x+\eta}d\eta,
\qquad\qquad\qquad\qquad\qquad\qquad\qquad\ \ \,
\end{multline}
see \eqref{eq:lambdaterm}. Substituting $\eta=\sin t$, $(t\in
[-\pi/2,\pi/2])$, and then $t=2 \arct u$, $(u\in[-1,1])$, we evaluate
\begin{multline}\label{eq:unnumbered}
\int_{-1}^1 \frac{(1-\eta^2)^{1/2}}{x+\eta}d\eta = (xt+\cos
t)|_{\pi/2}^{\pi/2} +(1-x^2) \int_{-\pi/2}^{\pi/2} \frac{dt}{x+\sin t}
\\  =  \pi x + 2(1-x^2)\int_{-1}^1 \frac{du}{x(1+u^2)+2u} =
\qquad\ \ \ \ \,\\
\qquad\quad
 =  \pi x - 2(x^2-1)^{1/2} \left[ \arct \sqrt{\frac{x+1}{x-1}} +
\arct \sqrt{\frac{x-1}{x+1}} \right] \\
 =  \pi (x-(x^2-1)^{1/2}).
\qquad\qquad\qquad\qquad\qquad\qquad\qquad\qquad\quad\ \,\,
\end{multline}
Combining this with \eqref{eq:multline2}, we obtain
$$
\frac{\partial I(\xi,x)}{\partial x} = -\pi +
\frac{\pi(x^2-1)^{1/2}}{x+\xi}.
$$
We integrate this equation from $x=1$ to $x=\gamma>1$, and use the
substitutions $x=\cosh t$, $t\in[0,t_0]$, with
$$ t_0 = \textrm{arccosh}\, \gamma = \log\left(
\gamma+(\gamma^2-1)^{1/2} \right), $$
and then $u=e^t$, $u\in[1,u_0]$, with
$$ u_0 = e^{t_0} = \gamma + (\gamma^2-1)^{1/2}. $$
We have
\begin{multline}\label{eq:multline3}
I(\xi,\gamma) = -\pi(\gamma-1) + \pi \int_0^{t_0} \frac{\sinh^2
t}{\cosh t + \xi}dt \\
= - \pi(\gamma-1)+\pi \left[ (\sinh t-\xi t)|_0^{t_0} + 2(\xi^2-1)
\int_0^{u_0} \frac{du}{u^2+2\xi u +1} \right]. \qquad\quad\ \ \,
\end{multline}
Here
\begin{equation}\label{eq:multline35}
(\sinh t-\xi t)|_0^{t_0} = (\gamma^2-1)^{1/2}-\xi
\log\left(\gamma+(\gamma^2-1)^{1/2} \right),
\end{equation}
and the last integral equals
\begin{multline}\label{eq:multline4}
\frac{1}{\sqrt{1-\xi^2}}\arct \frac{u+\xi}{(1-\xi^2)^{1/2}}
\bigg|_1^{u_0} = \frac{1}{\sqrt{1-\xi^2}} \arct
\frac{(u_0-1)(1-\xi^2)^{1/2}}{1-\xi^2+(u_0+\xi)(1+\xi)} \\
= \frac{1}{\sqrt{1-\xi^2}}
\arct\frac{u_0-1}{u_0+1}\sqrt{\frac{1+\xi}{1-\xi}} =
\frac{1}{\sqrt{1-\xi^2}} \arct
\sqrt{\frac{(\gamma-1)(1-\xi)}{(\gamma+1)(1+\xi)}}.
\end{multline}
Combining \eqref{eq:multline3}, \eqref{eq:multline35},
\eqref{eq:multline4} gives \eqref{eq:evaluation}.  \qed 

\bigskip
Now from \eqref{eq:lambdaterm}, \eqref{eq:logterm} and
\eqref{eq:evaluation} we get
\begin{multline}\label{eq:arctrad}
\qquad\qquad
h'(s) = \frac{c}{(\beta^2-s^2)^{1/2}} +
\frac{s}{\pi(\beta^2-s^2)^{1/2}}\left( \lambda- 2 \log
\frac{1+\sqrt{1-2 \beta^2}}{\sqrt{2} \beta} \right) \\
 + \frac{2}{\pi} \left( 
\arct \sqrt{\frac{(\gamma-1)(1+\xi)}{(\gamma+1)(1-\xi)}} -
\arct \sqrt{\frac{(\gamma-1)(1-\xi)}{(\gamma+1)(1+\xi)}} \right),
\qquad
\end{multline}
with $\xi=s/\beta$, $\gamma=\sqrt{2}/(2\beta)$, or, after some
simplification,
\begin{eqnarray*}
h'(s) &=& \frac{c}{(\beta^2-s^2)^{1/2}} +
\frac{s}{\pi(\beta^s-s^2)^{1/2}}\left( \lambda- 2 \log
\frac{1+\sqrt{1-2 \beta^2}}{\sqrt{2} \beta} \right) \\
& & + \frac{2}{\pi} \arct\frac{(1-2\beta^2)^{1/2}
s}{(\beta^2-s^2)^{1/2}}\,.
\end{eqnarray*}
We now observe that the only values of $c$ and $\lambda$ for which
the right-hand side is bounded as $s\nearrow \beta$, $s\searrow-\beta$,
and therefore has a chance of being the derivative of an
$\alpha$-admissible function, are
\begin{equation}\label{eq:lambdac}
c = 0, \qquad \lambda= 2 \log\frac{1+\sqrt{1-2 \beta^2}}{\sqrt{2}
\beta}.
\end{equation}
Therefore we get
\begin{equation}\label{eq:hprime}
h'(s) = \frac{2}{\pi} \arct\frac{(1-2\beta^2)^{1/2}
s}{(\beta^2-s^2)^{1/2}}.
\end{equation}
Note that $h'(s) \in (-1,1)$. We have determined $h'(s)$, except the
value of $\beta=\beta(\alpha)$ such that $h$ is $\alpha$-admissible,
i.e., satisfies \eqref{eq:alphacondition}. Rewrite
\eqref{eq:alphacondition} as
\begin{equation}\label{eq:alphabetacond}
\int_{-\beta}^\beta s h'(s)ds = \alpha - \beta^2.
\end{equation}
Besides evaluating this last integral, to compute $h(s)$ explicitly we
will need $\int_{-\beta}^s h'(u)du$. To this end, integrating the
first arctangent-of-radical function in \eqref{eq:arctrad} on the
interval $[-1,\xi]$, $(\xi\in(-1,1])$, we get
\begin{multline} \label{eq:eqx1}
\qquad
\int_{-1}^\xi \arct
\sqrt{\frac{(\gamma-1)(1+\eta)}{(\gamma+1)(1-\eta)}} d\eta \\ =
\xi \arct \sqrt{\frac{(\gamma-1)(1+\eta)}{(\gamma+1)(1-\eta)}} -
\frac{\sqrt{\gamma^2-1}}{2} \int_{-1}^\xi \frac{\eta\,
d\eta}{(\gamma-\eta)\sqrt{1-\eta^2}}. \qquad
\end{multline}
Substituting in the last integral $\eta=\sin t$, and then $u=\tan t$,
we transform it into
\begin{multline}\label{eq:eqx2}
\qquad
-t_0 - \frac{\pi}{2} + \gamma \int_{-\pi/2}^{t_0}
 \frac{dt}{\gamma-\sin t} \qquad
\qquad \qquad\qquad [t_0 = \arcs \xi] \\
= - t_0 - \frac{\pi}{2} + 2 \int_{-1}^{u_0}
 \frac{du}{1+u^2-2u/\gamma} \qquad \quad [u_0 = \tan t_0/2] 
\qquad\qquad \\
= -t_0 - \frac{\pi}{2} + \frac{2\gamma}{\sqrt{\gamma^2-1}}
\left( \arct \frac{u_0-\gamma^{-1}}{\sqrt{1-\gamma^{-2}}}+
\arct \frac{1+\gamma^{-1}}{\sqrt{1-\gamma^{-2}}} \right) 
\ \ \ 
\\ =
-t_0 - \frac{\pi}{2} + \frac{2\gamma}{\sqrt{\gamma^2-1}}\arct \left(
 \frac{1+u_0}{1-u_0} \sqrt{\frac{\gamma-1}{\gamma+1}} \right);
\qquad \qquad \qquad\qquad\quad\ \ 
\end{multline}
here
\begin{equation}\label{eq:eqx3}
\frac{1+u_0}{1-u_0} = \frac{1+\tan t_0/2}{1-\tan t_0/2} = \frac{1+\sin
t_0}{\cos t_0} = \frac{1+\xi}{\sqrt{1-\xi^2}} =
\sqrt{\frac{1+\xi}{1-\xi}}.
\end{equation}
From \eqref{eq:eqx1}, \eqref{eq:eqx2}, \eqref{eq:eqx3} we obtain
\begin{multline}\label{eq:eqx4} \qquad
 \int_{-1}^\xi \arct
 \sqrt{\frac{(\gamma-1)(1+\eta)}{(\gamma+1)(1-\eta)}}d\eta \\
 = (\xi-\gamma) \arct
 \sqrt{\frac{(\gamma-1)(1+\eta)}{(\gamma+1)(1-\eta)}} +
 \frac{\sqrt{\gamma^2-1}}{2}\left(\arcs\xi+\frac{\pi}{2}\right). \qquad
\end{multline}
In a similar fashion
\begin{equation}\label{eq:eqx5}
\int_{-1}^1 \eta \arct \sqrt{
\frac{(1+\eta)(\gamma-1)}{(1-\eta)(\gamma+1)}}d\eta =
\frac{\pi}{4}(1-\gamma^2+\gamma \sqrt{\gamma^2-1} ),
\end{equation}
and the integral in the negative arctangent in \eqref{eq:arctrad} is
obviously given by \eqref{eq:eqx5} as well. Using \eqref{eq:arctrad}
and \eqref{eq:eqx5}, we see that the $\alpha$-condition
\eqref{eq:alphabetacond} is equivalent to
$$ \beta^2(\gamma^2-\gamma \sqrt{\gamma^2-1}) = \alpha \iff
1-2\alpha = \sqrt{1-2\beta^2}, $$
the latter being possible only if $\alpha < 1/2$. In that case
\begin{equation}\label{eq:betaalpha}
\beta = \sqrt{2\alpha(1-\alpha)}.
\end{equation}
Consequently, see \eqref{eq:lambdac},
\begin{equation}\label{eq:lambdavalue}
\lambda = \log \frac{1-\alpha}{\alpha},
\end{equation}
and, see \eqref{eq:hprime},
\begin{equation}\label{eq:hprime2}
h'(s) = 
\frac{2}{\pi}\arct\left(\frac{(1-2\alpha)s}{\sqrt{2\alpha(1-\alpha)-s^2}}
\right), \qquad s\in
(-\sqrt{2\alpha(1-\alpha)},\sqrt{2\alpha(1-\alpha)}).
\end{equation}
Furthermore, denoting the integral in \eqref{eq:eqx4} by
$J(\xi,\gamma)$, we easily get
\begin{multline}\label{eq:hnotprime}
\qquad
h(s)= \beta + \int_{-\beta}^s h'(t)dt = \beta(1+J(\xi,\gamma) +
J(-\xi,\gamma) - J(1,\gamma)) \\ = 
\frac{2}{\pi}
s\arct\left(\frac{(1-2\alpha)s}{\sqrt{2\alpha(1-\alpha)-s^2}}
\right) +
\frac{\sqrt{2}}{\pi}
\arct
\left(\frac{\sqrt{2(2\alpha(1-\alpha)-s^2)}}{1-2\alpha}\right).
\end{multline}
We have derived a formula for a candidate minimizer, which we now
recognize as the function $\tilde{g}_\alpha$ that we defined in
section 2.  To be sure, this function was determined so as to meet the
ramifications of \emph{some} of the constraints. However, looking at
\eqref{eq:hprime2}, we see that $-1 < h'(s) < 1$ for $s\in
(-\beta,\beta)$, so $h$ is indeed $1$-Lipschitz, even though so far we
haven't paid attention to this constraint! Furthermore, since
$h^\prime(s)$ is odd, the constraint \eqref{eq:integral0} is met
automatically, and it is the reason why we were able to satisfy the
boundary constraints $h(-\sqrt{2}/2) = h(\sqrt{2}/2) =
\sqrt{2}/2$. Also, we determined $\beta$ from the requirement that $h$
should satisfy \eqref{eq:alphabetacond}, which under these boundary
conditions is equivalent to the $\alpha$-condition. We conclude that,
at the very least, $\tilde{g}_{\alpha}$ meets all the constraints,
thus is $\alpha$-admissible.

By Lemma 7, to prove that $\tilde{g}_\alpha$ is the minimizer, it only
remains to prove that $\tilde{g}_{\alpha}$ satisfies the conditions
\eqref{eq:necessarycond}. By \eqref{eq:diffmultline},
$w^\prime(s)\equiv 0$ for $|s|<\beta$. And $w(0) =0$ as $h^\prime(t)$
is odd. So $w(s)\equiv 0$ for $|s|<\beta$, hence the first condition
in \eqref{eq:necessarycond} is met. As for the remaining conditions,
by (anti)symmetry, it suffices to check, say, the third condition,
namely that
$$ F(s,\alpha) := -\int_{-\sqrt{2}/2}^{\sqrt{2}/2}
\tilde{g}_\alpha'(t)
\log|s-t| dt - \lambda(\alpha) s \le 0,\qquad \beta(\alpha) \le s \le
\sqrt{2}/2.
$$
Fix $0<s\le \sqrt{2}/2$, and let $\hat{\alpha}=(1-\sqrt{1-2s^2})/2$,
so that $\beta(\hat{\alpha})=s$. Clearly, because of the first
condition in \eqref{eq:necessarycond}, $F(s,\hat{\alpha})=0$. To
finish the proof, we will now show that $\partial F(s,\alpha)/\partial
\alpha > 0$ for $0<\alpha<\hat{\alpha}$. By \eqref{eq:multline},
\begin{equation}\label{eq:diffbeta}
\frac{\partial F(s,\alpha)}{\partial \beta} = - \int_{-\beta}^\beta
\frac{\partial \tilde{g}_\alpha'(t,\alpha)}{\partial \beta} \log |s-t|
dt - s \frac{\partial \lambda}{\partial \beta}.
\end{equation}
Using \eqref{eq:hprime} and simplifying gives
\begin{equation*}
\frac{\partial \tilde{g}_\alpha'(t)}{\partial \beta} = -\frac{2}{\pi
\beta(1-2 \beta^2)^{1/2}}\cdot \frac{t}{(\beta^2-t^2)^{1/2}}.
\end{equation*}
Since $\beta'(\alpha) = (1-2\beta^2)^{1/2}/\beta$, \eqref{eq:diffbeta}
becomes
\begin{equation*}
\frac{\partial F(s,\alpha)}{\partial \alpha} = \frac{2}{\pi \beta^2}
\int_{-\beta}^\beta \frac{t \log|s-t|}{(\beta^2-t^2)^{1/2}} dt +
\frac{s}{(1-\alpha)\alpha}.
\end{equation*}
Here the integral equals
$$
-(\beta^2-t^2)^{1/2} - \log|s-t|\bigg|_{-\beta}^\beta -
 \int_{-\beta}^\beta \frac{(\beta^2-t^2)^{1/2}}{s-t}dt = -\pi
 (s-(s^2-\beta^2)^{1/2}),
$$
see \eqref{eq:evaluation}. Therefore
\begin{eqnarray*}
\frac{\partial F(s,\alpha)}{\partial \alpha} &=& -\frac{2}{\beta^2}
\left(s-(s^2-\beta^2)^{1/2}\right) + \frac{s}{(1-\alpha)\alpha} \\
&=& s\left( \frac{1}{(1-\alpha)\alpha}-\frac{2}{\beta^2}\right) +
\frac{2}{\beta^2} (s^2-\beta^2)^{1/2} \\ & = &
\frac{2}{\beta^2} (s^2-\beta^2)^{1/2} > 0.
\end{eqnarray*}

\subsection{Direct computation of $K(\tilde{g}_\alpha)$}

Our next goal in this section is to show that
$K(\tilde{g}_\alpha)=-H(\alpha)+\log 2$. There are two ways to do
this. First, looking at the proof of Theorem 8, we see that we may
repeat the arguments of that proof (without assuming the value of
$K(\tilde{g}_\alpha)$ as in that proof) to deduce that the value
$M_\alpha$ of $K(\tilde{g}_\alpha)+H(\alpha)-\log 2$ \emph{must} be
$0$. For, if it were greater than $0$, then, denoting $k=\lfloor
\alpha n^2\rfloor$, we would have
\begin{eqnarray*}
1 = \pr_n\left( T \in {\cal T}_n \right) &=& 
\sum_{\lambda_0\textrm{ of area }k}
\pr_n\bigg( T\in{\cal T}_n : \lambda_T^k =
\lambda_0 \bigg) \\ &\le & p(n^2) \exp\left( - (1+o(1)) n^2 M_{k/n^2}
\right) \xrightarrow[n\to\infty]{} 0
\end{eqnarray*}
(since $M_\alpha$ is obviously continuous in $\alpha$.) On the other
hand, if $M_\alpha < 0$, then for some sufficiently large $n$, we
would have for some diagram $\lambda_0$ of area $\lfloor \alpha n^2
\rfloor$ contained in $\square_n$, that
$K(g_{\lambda_0})+H(\alpha)-\log 2 < 0$ (take a diagram for which
$g_{\lambda_0}$ approximates $\tilde{g}_\alpha$, and use Lemma 2). But
this again implies a contradiction:
$$ 1 \ge \pr_n\left( T \in {\cal T}_n : \lambda_T^{\lfloor \alpha n^2
\rfloor} = \lambda_0 \right) = \exp\left(-(1+o(1))n^2
(K(g_{\lambda_0})+H(\alpha)-\log 2)\right) > 1. $$

These last remarks notwithstanding, we find it worthwhile to compute
$K(\tilde{g}_\alpha)$ directly, if only to thoroughly test our
derivation of $\tilde{g}_\alpha$, and to show that all the integrals
involved can be evaluated explicitly.

For $h=\tilde{g}_\alpha$, rewrite \eqref{eq:alphacondition} as
$$ -\int_{-\sqrt{2}/2}^{\sqrt{2}/2} u(h'(u)-\sgn u)du = \alpha. $$
Using this, multiply both sides of \eqref{eq:necessarycond}
by $(h'(s)-\sgn s)$ and integrate, obtaining
\begin{multline}\label{eq:multline5}
K(h) = -\frac{\lambda \alpha}{2} -
\frac{1}{2}\int_{-\sqrt{2}/2}^{\sqrt{2}/2} h'(t) \bigg[ 2t \log|t| -
(t+\sqrt{2}/2)\log |t+\sqrt{2}/2| \\
- (t-\sqrt{2}/2)\log | t-\sqrt{2}/2 | \bigg] dt,
\qquad\qquad\ \ \,
\end{multline}
where we found before that $\lambda = \log((1-\alpha)/\alpha)$.
Denote
$$ S(t) = 2t \log|t| - (t+\sqrt{2}/2)\log |t+\sqrt{2}/2| \\ -
(t-\sqrt{2}/2)\log | t-\sqrt{2}/2 |, $$
and set 
$$K_1(h) = \int_{-\sqrt{2}/2}^{\sqrt{2}/2} h'(t)S(t)dt,$$
so that $K(h) = -\lambda\alpha/2 - K_1(h)/2$. Just like
\eqref{eq:diffbeta},
\begin{multline}\label{eq:diffbeta2}
\frac{\partial K_1(h_\alpha)}{\partial \beta} =
\int_{-\sqrt{2}/2}^{\sqrt{2}/2} \frac{\partial h_\alpha'(t)}{\partial
\beta} S(t) dt = \frac{2}{\pi\beta(1-2 \beta^2)^{1/2}}
\int_{-\beta}^\beta \frac{-t}{(\beta^2-t^2)^{1/2}}S(t)dt \\
= -\frac{2}{\pi\beta(1-2\beta^2)^{1/2}} \int_{-\beta}^\beta
(\beta^2-t^2)^{1/2}\bigg[ 2 \log|t|
-\log | t+ \sqrt{2}/2| - \log|t-\sqrt{2}/2| \bigg] dt.
\end{multline}
Denote
$$ E(s,\beta) = \int_{-\beta}^\beta (\beta^2-t^2)^{1/2}\log|t-s| dt,
$$
so that
\begin{equation}\label{eq:diffbeta3}
\frac{\partial K_1(h_\alpha)}{\partial \beta} =
2E(0,\beta)-E(-\sqrt{2}/2,\beta) - E(\sqrt{2}/2,\beta).
\end{equation}
By \eqref{eq:lambdaterm} and \eqref{eq:evaluation},
\begin{multline}\label{eq:diffs}
\frac{\partial E(s,\beta)}{\partial s} = \int_{-\beta}^\beta
\frac{(\beta^2-t^2)^{1/2}}{s-t}dt \\
= \left\{ \begin{array}{ll} \pi s & |s|< \beta, \\
   \pi (\sgn s)\left(|s|-(s^2-\beta^2)^{1/2}\right), &
 \beta < |s| < \sqrt{2}/2. \end{array}\right.
\qquad\qquad\qquad\quad\ \ 
\end{multline}
Therefore
\begin{equation}\label{eq:EEE}
2 E(0,\beta) = E(\beta,\beta)+E(-\beta,\beta)+ \pi\int_\beta^0 s ds +
\pi \int_{-\beta}^0 s ds = E(\beta,\beta)+ E(-\beta,\beta)-\pi
\beta^2.
\end{equation}
Likewise
\begin{equation}\label{eq:EEE2}
E(-\sqrt{2}/2,\beta)+E(\sqrt{2}/2,\beta) =
E(-\beta,\beta)+E(\beta,\beta) + 2\pi \int_\beta^{\sqrt{2}/2} \left(
s-(s^2-\beta^2)^{1/2}\right) ds,
\end{equation}
where
\begin{multline}\label{eq:multline6}
\int_\beta^{\sqrt{2}/2} (s^2-\beta^2)^{1/2}ds = \frac{1}{2}\left[
s(s^2-\beta^2)^{1/2}-\beta^2 \log\left( s+(s^2-\beta^2)^{1/2} \right)
\right]\bigg|_\beta^{\sqrt{2}/2} \\
= \frac{1}{2}\left( \frac{1-2\alpha}{2}-\alpha(1-\alpha)\log
\frac{1-\alpha}{\alpha} \right).
\end{multline}
So, using $\beta=(2\alpha(1-\alpha))^{1/2}$,
$$
E(-\sqrt{2}/2,\beta)+E(\sqrt{2}/2,\beta)=E(-\beta,\beta)+E(\beta,\beta)
+ \pi\left(-\beta^2+\alpha+\alpha(1-\alpha)\log\frac{1-\alpha}{\alpha}
\right),
$$
and, combining this relation with \eqref{eq:EEE}, we simpify
\eqref{eq:diffbeta3} to
$$
\frac{\partial K_1(h_\alpha)}{\partial \beta} = -\pi\left(
\alpha+\alpha(1-\alpha) \log\frac{1-\alpha}{\alpha} \right).
$$
So, by \eqref{eq:diffbeta2}
$$
\frac{\partial K_1(h_\alpha)}{\partial \alpha} =
\frac{\partial K_1(h_\alpha)}{\partial \beta}\cdot
\frac{(1-2\beta^2)^{1/2}}{\beta} = \frac{1}{1-\alpha} +
\log\frac{1-\alpha}{\alpha}. $$
Since $h_\alpha' \equiv 0$ at $\alpha=1/2$, we have $K_1(h)=0$ at
$\alpha=1/2$. Hence
\begin{multline}
K_1(h_\alpha) = \int_{1/2}^\alpha \left(\frac{1}{1-x}+\log
\frac{1-x}{x} \right) dx = -\log(1-\alpha)-2\log
2 \\ -(1-\alpha)\log(1-\alpha) -\alpha\log\alpha, \qquad\qquad\qquad
\end{multline}
which gives finally for $K(h_\alpha)$
$$ K(h) = \alpha\log \alpha+(1-\alpha)\log(1-\alpha)+\log 2 =
-H(\alpha) + \log 2.
$$
The proof of Theorem 7 is complete.
\qed

\subsection{The parametric family $\tilde{g}_\alpha$}

The minimality proof in section 3.2 relied on the possibility to
consider simultaneously the whole family of variational problems, and
thus to differentiate the minimizer $\tilde{g}_{\alpha}$ with respect
to $\alpha$. Moreover, to reveal a little secret, we anticipated the
formula \eqref{eq:lambdavalue} for the Lagrange multiplier
$\lambda$. According to a general (semiformal) recipe of the calculus
of variations (more specifically, mathematical programming), we
knew that this $\lambda$, dual to the $\alpha$-condition, should be
equal to $d K(\tilde{g}_\alpha)/d\alpha$, which we have proved to
be correct. The advantages of this approach of varying the parameter
$\alpha$ go even deeper than that. It will turn out that the partial
derivative of the minimizer $g_{\alpha}(\cdot)$ with respect to
$\alpha$ is the key to the distributional properties of the random
tableau. Using the formula for the minimizer, we compute easily that
\begin{equation}\label{eq:partialdiff}
\frac{\partial \tilde{g}_\alpha(u)}{\partial \alpha} = \left\{
\begin{array}{lll}
  0  & \qquad\qquad\quad\ 
& \sqrt{2\alpha(1-\alpha)} < |u| \le \sqrt{2}/2 \\
 \frac{\sqrt{2\alpha(1-\alpha)-u^2}}{\pi\alpha(1-\alpha)} & &
     |u| \le \sqrt{2\alpha(1-\alpha)} \end{array} \right.
\end{equation}
For each $\alpha$, direct integration reveals that
$\partial\tilde{g}_\alpha(u)/\partial\alpha$ is a probability density
function, i.e.
$$ \int_{-\sqrt{2}/2}^{\sqrt{2}/2}
\frac{\partial\tilde{g}_\alpha(u)}{\partial \alpha}\,du = 1. $$
(In fact, it is the density of the semicircle distribution, and it will play
a prominent role later. See sections 5, 8.1, 8.2.) This observation is in 
perfect harmony with the fact that $\tilde{g}_{\alpha}$ satisfies the 
$\alpha$-condition, thus providing a partial check of our computations. Indeed
\begin{eqnarray*}
\int_{-\sqrt{2}/2}^{\sqrt{2}/2} \left(\tilde{g}_\alpha(u)-|u|\right)du
&=& \int_{-\sqrt{2}/2}^{\sqrt{2}/2}
\left(\tilde{g}_\alpha(u)-\tilde{g}_0(u)\right)du 
\qquad\qquad\qquad\qquad\qquad\qquad\qquad
\\ & \hspace{-120.0pt}=& \hspace{-60.0pt}
\int_{-\sqrt{2}/2}^{\sqrt{2}/2} \left(\int_0^\alpha
\frac{\partial\tilde{g}_s(u)}{\partial s}ds\right)du =
\int_0^\alpha \left( \int_{-\sqrt{2}/2}^{\sqrt{2}/2}
 \frac{\partial\tilde{g}_s(u)}{\partial s} du\right) ds \\
& \hspace{-120.0pt}=& \hspace{-60.0pt}
\int_0^\alpha 1\, ds = \alpha.
\end{eqnarray*}
Had we been presented with the minimizer $\tilde{g}_{\alpha}$ ``out of
the blue'', this would have been the least computational way to prove
its $\alpha$-admissibility.

\section{The boundary of the square}

\subsection{Proof of Theorem 3}

In this section, we prove Theorem 3. As was remarked in section 1.3,
the RSK correspondence induces a correspondence between minimal
Erd\"os-Szekeres permutations $\pi$ of $1,2,\ldots,n^2$ and pairs
$T_1,T_2\in {\cal T}_n$ of square tableaux. By the well known result
of Schensted \cite{schensted}, in this correspondence the length
$l_{n,k}$ of the longest increasing subsequence in
$\pi(1),\pi(2),\ldots,\pi(k)$ is equal to the
length
$\lambda_{T_1}^k(1)$ of the first row of $\lambda_{T_1}^k$. So
the
distribution of $l_{n,k}$ under a uniform random choice of
minimal
Erd\"os-Szekeres permutation $\pi$ is equal to the
distribution of the
length of the first row of $\lambda_T^k$ in a
uniform random square
tableau $T\in{\cal T}_n$. Denoting for the
remainder of this section
$\alpha=\alpha(k)=k/n^2$, we can therefore
reformulate Theorem 3 as
stating that
\begin{equation}\label{eq:boundary}
\max_{\alpha_0\le k/n^2\le
1/2}
\mathbb{P}_n\left( T\in{\cal T}_n :
 \left|\lambda_T^{k}(1) -
2\sqrt{\alpha(1-\alpha)}n\right| > \alpha_0^{1/2}\omega(n)n \right)
\xrightarrow[n\to\infty]{} 0.
\end{equation}
Theorem 8 looks as if
it might imply \eqref{eq:boundary}. In fact, it
only implies a lower
bound on $\lambda_T^k(1)$. The reason is that
$g_{\lambda_T^k}$ can
be very close in the supremum norm to
$\tilde{g}_\alpha$ (as is known
to happen with high probability by
Theorem 8), while $n^{-1}
\lambda_T^k(1)$ might still be much larger
than $2
\sqrt{\alpha(1-\alpha)}$ (see \eqref{eq:firstrow}
below).

\paragraph{Lemma 9.} Let $\alpha_0=n^{-2/3+\epsilon}$,
$\delta=n^{-1/3(1-\epsilon)
}$, $\epsilon\in (0,2/3)$. Then

\begin{equation}\label{eq:boundaryweak}
\mathbb{P}_n\left(
T\in{\cal T}_n :
 \min_{\alpha_0\le\alpha\le 1/2} (\lambda_T^{k}(1)
-2\sqrt{\alpha(1-\alpha)}n)
\le -\delta n \right) = O(n^{-b})
\end{equation}
for every $b>0$.

\paragraph{Proof.} We use the
notation of Theorem 8. The length of the
first row $\lambda_T^k(1)$
can be extracted from the rotated
coordinate graph $g_{\lambda_T^k}$
using the following
relation:
\begin{equation}\label{eq:firstrow}
\frac{1}{n}\lambda_T^k(1)
= \sqrt{2}\, \inf\left\{ u \in [0,\sqrt{2}/2]
: g_{\lambda_T^k}(u) =
u \right\}.
\end{equation}
It follows from \eqref{eq:hprime2} that,
uniformly for $\alpha\in [\alpha_0,1/2]$
and $|u|<
\sqrt{2\alpha(1-\alpha)}$,
$$
|\partial\tilde{g}_{\alpha}(u)/\partial
u
-1|=\frac{2}{\pi}\tan^{-1}\frac{\sqrt
{2\alpha(1-\alpha)-u^2}}{(1-2\alpha)|u|}\ge
c(\sqrt{2\alpha(1-\alpha)}-|u|),
$$
$c>0$ being an absolute
constant. Consequently, for $\alpha\in
[\alpha_0,1/2]$,
$$
\tilde{g}_\alpha(\sqrt{2\alpha(1-\alpha)}-\delta)-(
\sqrt{2\alpha(1-\alpha)}-\delta)
\ge c\delta^{3/2}.$$ 
So if $T\in {\cal T}_n$ has the property that,
for some $k$ in question, 
$$ 
\lambda_T^{k} (1) -
2\sqrt{\alpha(1-\alpha)}n
< -\delta n, $$
then by
\eqref{eq:firstrow},
\begin{eqnarray*}
||g_{\lambda_T^{k
}}-\tilde{g}_\alpha||_\infty
&\ge& \sup \{
|g_{\lambda_T^{k}}(u)-\tilde{g}_\alpha(u)|
:
\sqrt{2\alpha(1-\alpha)}-\delta < u < \sqrt{2\alpha(1-\alpha)} \}
\\
& = & \sup \{ \tilde{g}_\alpha(u)-u
:
\sqrt{2\alpha(1-\alpha)}-\delta< u < \sqrt{2\alpha(1-\alpha)} \}
\ge
c\delta^{3/2}.
\end{eqnarray*}
So, by Theorem 8 with
$\epsilon:=c\delta^{3/2}$, 
\begin{multline*}
\qquad\qquad \mathbb{P}_n\left( T\in{\cal T}_n :
\min_{\alpha_0\le\alpha\le 1/2}(\lambda_T^{k} (1) -
2\sqrt{\alpha(1-\alpha)}n)
\le -\delta n \right) \\
\le \mathbb{P}_n\left( T\in{\cal T}_n : \max_{\alpha_0\le\alpha\le
1/2} || g_{\lambda_T^{k}}-\tilde{g}_\alpha||_\infty \ge
c\delta^{3/2}\right) \\ \le \exp(3n-\hat{c}n^2\delta^{3r/2})
\le \exp(3n-\hat{c}n^{2-r(1-\epsilon)/2}) \xrightarrow[n\to\infty]{}
0,\qquad\qquad\ 
\end{multline*}
provided that we choose a feasible $r$, i. e. $r\in (2,3)$, such
that $r<2(1-\epsilon)^{-1}$. 
\qed

\bigskip
To prove the upper bound and thus conclude the proof of Theorem 3, it
suffices to prove an upper bound for the \emph{expected value} of
$\lambda_T^k$, namely that, for $\alpha_0\le\alpha\le 1/2$,
\begin{equation}\label{eq:expected}
\mathbb{E}_n \left[\lambda_T^{k}(1) \right] \le 2\sqrt{\alpha(1-\alpha)}n +
O(\alpha_0^{1/2}n),
\end{equation}
where $\mathbb{E}_n$ denotes expectation with respect to the
probability measure $\mathbb{P}_n$. 
Indeed, choosing $\omega(n)\to\infty$ however slowly, we bound
$$
\mathbb{P}_n\left(T\in{\cal T}_n : \lambda_T^{k} (1)
\ge 2\sqrt{\alpha(1-\alpha)}n+ \alpha_0^{1/2}\omega(n)n\right)
\qquad\qquad\qquad\qquad\qquad\qquad\qquad $$

\vspace{-27.0pt}
\begin{eqnarray*}
\textrm{(by Markov's inequality)}
 & \le & (\alpha_0^{1/2}\omega(n)n)^{-1}\mathbb{E}_n
\left[\max(0,\lambda_T^k(1)-2\sqrt{\alpha(1-\alpha)}n)\right]\\
\textrm{(by Lemma 9, for any $b>0$)}
& \le & (\alpha_0^{1/2}\omega(n)n)^{-1}
\left(\mathbb{E}_n\left[\lambda_T^k(1)-
2\sqrt{\alpha(1-\alpha)}+\delta n\right]+O(n^{1-b}) \right)\\ 
& = & O((\alpha_0^{1/2}n+\delta n)/(\alpha_0^{1/2}\omega(n)n))=O(\omega(n)^{-1}).
\end{eqnarray*}
\bigskip  
Write
$$ \lambda_T^k(1) = \sum_{j=1}^k I_{n,j}, $$
where
$I_{n,j} = \lambda_T^j(1)-\lambda_T^{j-1}(1) = $
indicator of the event that $\lambda_T^j$ is obtained from
$\lambda_T^{j-1}$ by adding a box to the first row. Let $p_{n,j} =
\mathbb{E}_n(I_{n,j})$.

\paragraph{Lemma 10.} In the notation of Lemma 9, as $n\to\infty$,
$$ p_{n,j} \le \frac{n^2-2j}{n\sqrt{j(n^2-j)}}+
O(\delta n(n^2-2j+1)^{-1}), $$
uniformly for $\alpha_0\le j/n^2\le 1/2$.

\paragraph{Proof.} Let ${\cal Y}_{n,j}$ be the set of Young diagrams
of area $j$ contained in the $n\times n$ square. For a diagram
$\lambda \in {\cal Y}_{n,j}$, denote by $\textrm{next}(\lambda)$ the
diagram obtained from $\lambda$ by adding a box to the first
row. Then, conditioning $I_{n,j}$ on the shape $\lambda_T^{j-1}$, we
write
\begin{eqnarray*}
p_{n,j} &=& \mathbb{P}_n\left( \lambda_T^j =
\textrm{next}(\lambda_T^{j-1})\right) = \sum_{ \lambda\in {\cal Y}_{n,j-1}}
\frac{d(\lambda) d(\square_n\setminus
\textrm{next}(\lambda))}{d(\square_n)} \\ &=&
\sum_{\lambda \in {\cal Y}_{n,j-1}}
\frac{d(\textrm{next}(\lambda))d(\square_n\setminus
\textrm{next}(\lambda))}{d(\square_n)}\cdot
\frac{d(\lambda)}{d(\textrm{next}(\lambda))}
\end{eqnarray*}
This is nearly an average over ${\cal Y}_{n,j}$ with respect to the
measure \eqref{eq:dmeasure}; in fact, slightly less, since not any
$\lambda'\in {\cal Y}_{n,j}$ is of the form $\textrm{next}(\lambda)$
for some $\lambda\in{\cal Y}_{n,j-1}$. It follows from the convexity
of the function $x\to x^2$ that
\begin{multline}\label{eq:beforeamus}
\qquad\qquad
p_{n,j}^2 \le \sum_{\lambda \in {\cal Y}_{n,j-1}}
\frac{d(\textrm{next}(\lambda))d(\square_n\setminus
\textrm{next}(\lambda))}{d(\square_n)}\cdot \left(
\frac{d(\lambda)}{d(\textrm{next}(\lambda))} \right)^2 \\ =
\sum_{\lambda \in {\cal Y}_{n,j-1}}
\frac{d(\lambda)d(\square_n\setminus \lambda)}{d(\square_n)}\cdot
\frac{d(\lambda)d(\square_n\setminus \textrm{next}(\lambda))}
{d(\textrm{next}(\lambda))d(\square_n\setminus \lambda)}.
\qquad\qquad\qquad\qquad\ 
\end{multline}
We now note the amusing identity
\begin{equation}\label{eq:amus}
\frac{d(\lambda)d(\square_n\setminus \textrm{next}(\lambda))}
{d(\textrm{next}(\lambda))d(\square_n\setminus \lambda)}
= \frac{n^2-\lambda(1)^2}{j(n^2-j+1)}, \qquad (\lambda \in {\cal Y}_{n,j-1})
\end{equation}
which follows from writing out the hook products for $d(\cdot)$ in
\eqref{eq:hook} and observing cancellation of almost all the factors -
see Figure 5. Here is a proof of \eqref{eq:amus}. Clearly the only
hook lengths influenced by this operation are of the cells in the
first row and the $(\lambda(1)+1)$-th column. In particular,
$$ \frac{d(\lambda)}{d(\textrm{next}(\lambda))} = \frac{1}{j}
\prod_{i=1}^{\lambda(1)}
\frac{\lambda(1)-i+1+\lambda'(i)}{\lambda(1)-i+\lambda'(i)};
$$
here $\lambda'(i)$ is the number of cells in the $i$-th column of
$\lambda$. Clearly the fraction factors ``telescope'' on each
subinterval of $[1,\lambda(1)]$ where $\lambda'(i)$ is constant. Let
$[i_1,i_2]$ be such a (maximal) subinterval. Maximality implies that
$(i_2,\lambda'(i_2))$ is a corner of $\lambda$, and that
$(\lambda'(i_1)+1,i_1)$ is a corner of $\square_n\setminus
\textrm{next}(\lambda)$. Then
$$
\prod_{i=1}^{\lambda(1)}
\frac{\lambda(1)-i+1+\lambda'(i)}{\lambda(1)-i+\lambda'(i)} =
\frac{\lambda(1)-i_1+1+\lambda'(i_1)}{\lambda(1)-i_2+\lambda'(i_2)}
=
\frac{h_{\square_n\setminus\textrm{next}(\lambda)}(\lambda'(i_1)+1,\lambda(1)+1)}{h_\lambda(1,u_2)}
$$
where, say, $h_\lambda(u,v)$ denotes the hook length for a cell $(u,v)
\in \lambda$. Multiplying these fractions for all such subintervals
$[i_1,i_2]$, we get
\begin{equation} \label{eq:amus2}
\frac{d(\lambda)}{d(\textrm{next}(\lambda))} = \frac{1}{j}
\left(\prod_{(u,v)\in \textrm{corners}(\lambda)} f(u,v)\right)^{-1} \cdot
\left(\prod_{(u,v)\in \textrm{corners}(\square_n\setminus \lambda)}
  g(u,v) \right).
\end{equation}
Here $\textrm{corners}(\mu)$ is the corner set of a diagram $\mu$;
$f(u,v)$ is the hook length of a cell in the first row of $\lambda$
whose vertical leg ends at the corner $(u,v)\in
\textrm{corners}(\lambda)$; $g(u,v)$ is the hook length of a cell in 
$\square_n\setminus\textrm{next}(\lambda)$ from the
$(\lambda(1)+1)$-th column whose horizontal arm ends at the corner
$(u,v) \in \textrm{corners}(\square_n\setminus \lambda)$. Next,
considering separately the first row cells $(1,v)$, $v>\lambda(1)$,
the top $\lambda'(1)$ cells in the $(\lambda(1)+1)$-th column, and
finally the bottom $n-\lambda'(1)$ cells in that column, we obtain
\begin{equation}\label{eq:amus3}
\frac{d(\square_n\setminus
\textrm{next}(\lambda))}{d(\square_n\setminus \lambda)} =
\frac{n-\lambda(1)}{n^2-j+1} \cdot \prod_{k=2}^{\lambda'(1)}
\frac{\lambda(1)-\lambda(k)+k}{\lambda(1)-\lambda(k)+k-1} \cdot
\frac{\lambda(1)+n}{\lambda(1)+\lambda'(1)}.
\end{equation}
Here, analogously to the $d(\lambda)/d(\textrm{next}(\lambda))$ case,
\begin{multline}\label{eq:amus4}
\qquad
\frac{1}{\lambda(1)+\lambda'(1)} \prod_{k=2}^{\lambda'(1)}
\frac{\lambda(1)-\lambda(k)+k}{\lambda(1)-\lambda(k)+k-1} \\ =
\left(\prod_{(u,v)\in \textrm{corners}(\lambda)} f(u,v)\right) \cdot
\left(\prod_{(u,v)\in \textrm{corners}(\square_n\setminus \lambda)}
  g(u,v) \right)^{-1}. \qquad
\end{multline}
Putting \eqref{eq:amus2}, \eqref{eq:amus3}, \eqref{eq:amus4} together
gives \eqref{eq:amus}.

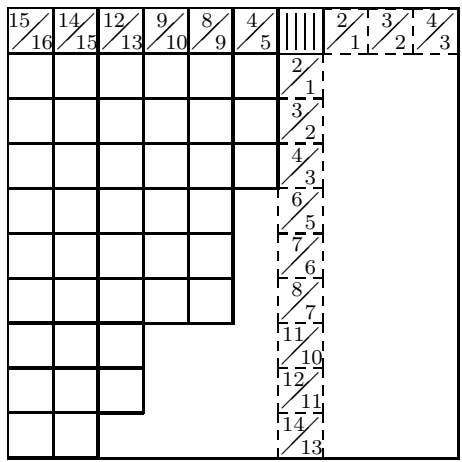
\begin{figure}[h!]
\begin{center}
\begin{picture}(200,200)(0,0)
\setlength{\unitlength}{1.7 pt}
\multiput(0,100)(0,-10){10}{\framebox(10,10)}
\multiput(10,100)(0,-10){10}{\framebox(10,10)}
\multiput(20,100)(0,-10){9}{\framebox(10,10)}
\multiput(30,100)(0,-10){7}{\framebox(10,10)}
\multiput(40,100)(0,-10){7}{\framebox(10,10)}
\multiput(50,100)(0,-10){4}{\framebox(10,10)}
\multiput(60,100)(0,-10){1}{\framebox(10,10)}
\multiput(62,101)(2,0){4}{\line(0,1){8}}
\put(0,10){\framebox(100,100)}
\multiput(60,90)(0,-10){9}{\dashbox{2}(10,10)}
\multiput(70,100)(10,0){3}{\dashbox{2}(10,10)}
\multiput(1,101)(10,0){6}{\line(1,1){8}}
\multiput(71,101)(10,0){3}{\line(1,1){8}}
\multiput(61,91)(0,-10){9}{\line(1,1){8}}
\put(0,106){\scriptsize{15}}
\put(5,101){\scriptsize{16}}
\put(11,106){\scriptsize{14}}
\put(15,101){\scriptsize{15}}
\put(21,106){\scriptsize{12}}
\put(25,101){\scriptsize{13}}
\put(33,106){\scriptsize{9}}
\put(35,101){\scriptsize{10}}
\put(43,106){\scriptsize{8}}
\put(46,101){\scriptsize{9}}
\put(53,106){\scriptsize{4}}
\put(56,101){\scriptsize{5}}
\put(73,106){\scriptsize{2}}
\put(76,101){\scriptsize{1}}
\put(83,106){\scriptsize{3}}
\put(86,101){\scriptsize{2}}
\put(93,106){\scriptsize{4}}
\put(96,101){\scriptsize{3}}
\put(63,96){\scriptsize{2}}
\put(66,91){\scriptsize{1}}
\put(63,86){\scriptsize{3}}
\put(66,81){\scriptsize{2}}
\put(63,76){\scriptsize{4}}
\put(66,71){\scriptsize{3}}
\put(63,66){\scriptsize{6}}
\put(66,61){\scriptsize{5}}
\put(63,56){\scriptsize{7}}
\put(66,51){\scriptsize{6}}
\put(63,46){\scriptsize{8}}
\put(66,41){\scriptsize{7}}
\put(61,36){\scriptsize{11}}
\put(65,31){\scriptsize{10}}
\put(61,26){\scriptsize{12}}
\put(65,21){\scriptsize{11}}
\put(61,16){\scriptsize{14}}
\put(65,11){\scriptsize{13}}
\thicklines
\put(0,10){\line(1,0){20}}
\put(20,10){\line(0,1){10}}
\put(20,20){\line(1,0){10}}
\put(30,20){\line(0,1){20}}
\put(30,40){\line(1,0){20}}
\put(50,40){\line(0,1){30}}
\put(50,70){\line(1,0){10}}
\put(60,70){\line(0,1){40}}
\end{picture}
\caption{Illustration of \eqref{eq:amus} for
$\lambda=(6,6,6,6,5,5,5,3,3,2)$: The numbers in the cells are the hook
lengths before and after the new cell is added.}
\end{center}
\end{figure}

Combining \eqref{eq:beforeamus} and \eqref{eq:amus} gives that
\begin{equation}\label{eq:boundpnj2}
p_{n,j}^2 \le \mathbb{E}_n\left[
\frac{n^2-\lambda_T^{j-1}(1)^2}{j(n^2-j+1)}\right]
\end{equation}
By Lemma 9, we may write
$$
\mathbb{E}_n(\lambda_T^{j-1}(1)) \ge 
\frac {2\sqrt{j(n^2-j)}}{n}-\delta n,$$
($\delta=n^{-(1-\epsilon)/3}$), for all $j/n^2\in [\alpha_0,1/2]$.
So,
using $\mathbb{E}_n^2[\lambda_T^{j-1}(1)]\le 
\mathbb{E}[(\lambda_T^{j-1}(1))^2]$,
$$ p_{n,j}^2 \le \frac{(n^2-2j)^2}{n^2\cdot j(n^2-j)}+\frac{4\delta}{\sqrt{j(n^2-j)}},
$$
or, using $(1+z)^{1/2}\le 1+z/2$ for $j<n^2/2$,
$$ p_{n,j}\le \frac{n^2-2j}{n\sqrt{j(n^2-j)}}+O(\delta n (n^2-2j+1)^{-1}).
$$
The estimate holds for $j=n^2/2$ as well, since 
$\delta^{1/2}n^2\to\infty$.
\qed

\bigskip
Note that \eqref{eq:boundpnj2} implies in particular the rough bound
$$
p_{n,j} \le \frac{n}{\sqrt{j (n^2-j+1)}},
$$
valid for all $j\le n^2$. Now, to complete the proof of Theorem 3, we use this bound for 
$j \le \alpha_0 n^2$ and Lemma 10 for $j>\alpha_0 n^2$. First
$$
\mathbb{E}_n\bigg[\lambda_T^{k}(1)\bigg] =
\sum_{j\le \alpha_0n^2} p_{n,j} +
\sum_{\alpha_0n^2<j\le k}p_{n,j}=\Sigma_1+\Sigma_2.$$
Here
$$
\Sigma_1\le 2\sum_{j\le\alpha_0n^2}j^{-1/2}=O(n\alpha_0^{1/2}),
$$
and
$$
\Sigma_2\le\sum_{\alpha_0n^2<j\le k}\frac{n^2-2j}{n\sqrt{j(n^2-j)}}+
O(\delta n\log n).
$$
The last sum is bounded above by
$$
n\int_{\alpha_0-n^{-2}}^{\alpha}\frac{1-2t}{\sqrt{t(1-t)}}\,dt=2n\sqrt{\alpha(1-\alpha)}
+O(n\alpha_0^{1/2}).
$$
Therefore, since $\alpha_0^{1/2}\gg \delta\log n$,
$$
\mathbb{E}_n[\lambda_T^k]\le 2n\sqrt{\alpha(1-\alpha)}+O(n\alpha_0^{1/2}).
$$
So \eqref{eq:expected} follows. Theorem 3 is proved.
\qed


\subsection{Proof of Theorem 1(i)}

With our enhanced understanding of the distribution of
$\lambda_T^k(1)$, we may now prove Theorem 1(i). First we show that
for individual boundary points, the tableau approaches the limit
surface. Fix $(x,y)$ on the boundary of the square. By symmetry, we
may assume that $y=0, 0<x<1$. Let $\epsilon>0$. Denote
$$ \alpha = L(x,0) = \frac{1-\sqrt{1-x^2}}{2}, $$ so that
$x=2\sqrt{\alpha(1-\alpha)}$. For any tableau $T \in {\cal T}_n$,
denote $k_T = t_{\lfloor nx\rfloor+1,1}$, and let $\beta_T =
k_T/n^2$. We want to show that with high probability,
$|\beta_T-\alpha|$ is small. Note that $k_T$ is an integer
representing the smallest $j$ for which $\lambda_T^j > nx$. Therefore
$nx\le \lambda_T^{k_T}(1) < nx+1$, or
\begin{equation}\label{eq:stam1}
\left| \lambda_T^{k_T}(1)-x \right| \le \frac{1}{n}
\end{equation}
The function $f(t):=L(t,0)=(1-\sqrt{1-t^2})/2$ is monotonically
increasing and uniformly continuous on $[0,1]$. Choose a $\delta>0$
such that $|t-t'|<\delta$ implies $|f(t)-f(t')|<\epsilon/3$. Choose
numbers $0=a_0<a_1<a_2<\ldots<a_N=1/2$ such that $a_{i+1}-a_i <
\epsilon/3$,\  $i=0,1,2,\ldots,N-1$. Denote
$x_i=f^{-1}(a_i)=2\sqrt{a_i(1-a_i)}$.

Let $T\in {\cal T}_n$ be a tableau that satisfies
\begin{equation}\label{eq:stam2}
\left| \frac{1}{n}\lambda_T^{\lfloor a_i n^2\rfloor}(1) - x_i\right| <
\frac{\delta}{2}, \qquad (i=1,2,\ldots,N)
\end{equation}
(this happens with high probability, by \eqref{eq:boundary}). Let
$0\le i<N$ be such that $a_i \le \beta_T < a_{i+1}$. Then clearly
\begin{equation}\label{eq:stam3}
x_i - \frac{\delta}{2} < \frac{1}{n}\lambda_T^{\lfloor a_i n^2 \rfloor}(1)
\le \frac{1}{n}\lambda_T^{k_T}(1) \le
\frac{1}{n}\lambda_T^{\lfloor a_{i+1} n^2\rfloor}(1) < x_{i+1}+
\frac{\delta}{2}
\end{equation}
Combining this with \eqref{eq:stam1} we get, for $n>2/\delta$,
$$ x_i - \delta < x < x_{i+1} + \delta.$$
Therefore
$$ a_i - \frac{\epsilon}{3} < \alpha=f(x) < a_{i+1}+\frac{\epsilon}{3},
$$
and, since also $a_i \le \beta_T < a_{i+1}$ and $a_{i+1}-a_i <
\epsilon/3$, we get
$$ |\beta_T-\alpha|<\epsilon. $$
Summarizing, we have shown that
\begin{multline}\label{eq:boundary2}
\mathbb{P}_n\left( T \in {\cal T}_n : |\beta_T-\alpha|<\epsilon
\right) \\ \ge
\mathbb{P}_n\left( T\in {\cal T}_n : \forall\  i=1,2,\ldots,N,\ \ \left|
\frac{1}{n}\lambda_T^{\lfloor a_i n^2\rfloor}(1)-x_i\right| <
\frac{\delta}{2} \right)
\xrightarrow[n\to\infty]{} 1.
\end{multline}
Theorem 1(i) now follows easily. It is enough to say that, because of
the monotonicity of the tableau $t_{i,j}$ as a function of $i$ and
$j$, and the monotonicity of the limit surface function $L$,
given $\epsilon>0$ we can find finitely many points
$(x_1,y_1),(x_2,y_2),\ldots,(x_N,y_N) \in [0,1]\times[0,1]\setminus
\{(0,0),(1,0),(0,1),(1,1)\} $ such that the event inclusion
\begin{multline}\label{eq:union}
\qquad
\bigg\{ T\in{\cal T}_n : \max_{1\le i,j\le n}
\left|\frac{1}{n^2}t_{i,j} -
L\left(\frac{i}{n},\frac{j}{n}\right)\right| > \epsilon \bigg
\} \subseteq \\ \bigcup_{i=1}^N
\bigg\{ T \in {\cal T}_n : \left| \frac{1}{n^2} 
t_{\lfloor n x_i\rfloor +1, \lfloor n y_i\rfloor+1} - L(x_i,y_i)
\right|  > \frac{\epsilon}{10}
\bigg\} \qquad\qquad\qquad
\end{multline}
holds. But now, the $\pr_n$-probability of each of the individual
events in this union tends to 0 as $n\to\infty$ -- because of Theorem
1(ii) for the points $(x_i,y_i)$ in the interior of the square (using
the continuity of the function $L$), and because of
\eqref{eq:boundary2} for the points on the boundary.  \qed

\section{The hook walk and the cotransition measure of a diagram}

In this section, we study the location of the $k$-th entry in the
random tableau $T \in {\cal T}_n$, when $k\approx \alpha\cdot
n^2$. The idea is to condition the distribution of the location of the
$k$-th entry on the shape $\lambda_T^k$ of the $k$-th subtableau of
$T$. Given the shape $\lambda_T^k$, the distribution of the location
of the $k$-th entry is exactly the so-called \emph{cotransition
measure} of $\lambda_T^k$ (see below). We know from Theorem 8 that
with high probability, the rescaled shape of $\lambda_T^k$ is
approximately described in rotated coordinates by the level curve
$v=\tilde{g}_\alpha(u)$. Romik \cite{romik} showed that the
cotransition measure is a continuous functional on the space of
continual Young diagrams, and derived an explicit formula for the
probability density of its $u$-coordinate. By substituting the level
curve $\tilde{g}_\alpha$ in the formula from \cite{romik}, we will get
exactly the semicircle density \eqref{eq:semicircle}, proving Theorem
2.

Let $\lambda:\lambda(1)\ge\lambda(2)\ge \ldots \ge\lambda(m)>0$ be a
Young diagram with $k=|\lambda|=\sum_i \lambda(i)$ cells. A cell
$c=(i,j)\in\lambda$ \ $(1\le i\le m, 1\le j \le \lambda(i))$ is called
a \emph{corner} cell if removing it leaves a Young diagram
$\lambda\setminus c$, or in other words if $j=\lambda(i)$ and ($i=m$
or $\lambda(i)>\lambda(i+1)$). If $T$ is a Young tableau of shape
$\lambda$, let $c_{\textrm{max}}(T)$ be the cell containing the
maximal entry $k$ in $T$. Obviously $c_{\textrm{max}}(T)$ is a corner
cell of $\lambda$.

The \emph{cotransition measure} of $\lambda$ is the probability
measure $\mu_\lambda$ on corner cells of $\lambda$, which assigns to a
corner cell $c$ measure
\begin{equation}\label{eq:cotransition}
\mu_\lambda(c) = \frac{d(\lambda \setminus c)}{d(\lambda)}
\end{equation}
(with $d(\lambda)$ as in \eqref{eq:hook}.) This is a probability
measure, since one may divide up the $d(\lambda)$ tableaux of shape
$\lambda$ according to the value of $c_{\textrm{max}}(T)$; for any
corner cell $c$, there are precisely $d(\lambda\setminus c)$ tableaux
for which $c_{\textrm{max}}(T)=c$. In other words $\mu_\lambda$ describes
the distribution of $c_{\textrm{max}}(T)$, for a uniform random choice
of a tableau $T$ of shape $\lambda$.

It is fascinating that there exists a simple algorithm to sample from
$\mu_\lambda$. This is known as the \emph{hook walk} algorithm of
Greene-Nijenhuis-Wilf, and it can be described as follows: Choose a
cell $c=(i,j)\in\lambda$ uniformly among all $k$ cells. Now execute a
random walk, replacing at each step the cell $c$ with a new cell $c'$,
where $c'$ is chosen uniformly among all cells which lie either to the
right of, or (exclusive or) below $c$. The walk terminates when a corner
cell is reached, and it can be shown \cite{greenenij} that the probability
of reaching $c$ is given by \eqref{eq:cotransition}. Figure 6 shows a
Young diagram, its corner cells and a sample hook walk path.

\begin{figure}[h!]
\begin{center}
\begin{picture}(100,110)(0,0)
\multiput(0,100)(0,-10){10}{\framebox(10,10)}
\multiput(10,100)(0,-10){10}{\framebox(10,10)}
\multiput(20,100)(0,-10){10}{\framebox(10,10)}
\multiput(30,100)(0,-10){8}{\framebox(10,10)}
\multiput(40,100)(0,-10){7}{\framebox(10,10)}
\multiput(50,100)(0,-10){5}{\framebox(10,10)}
\multiput(60,100)(0,-10){5}{\framebox(10,10)}
\multiput(70,100)(0,-10){5}{\framebox(10,10)}
\multiput(80,100)(0,-10){3}{\framebox(10,10)}
\multiput(90,100)(0,-10){1}{\framebox(10,10)}
\put(15,85){\circle{3}}
\put(45,85){\circle{3}}
\put(45,65){\circle{3}}
\put(65,65){\circle{3}}
\put(75,65){\circle{3}}
\multiput(18,85)(3,0){9}{\circle*{1}}
\multiput(45,81)(0,-3){5}{\circle*{1}}
\multiput(49,65)(3,0){5}{\circle*{1}}
\multiput(69,65)(3,0){2}{\circle*{1}}
\multiput(22,11)(2,0){4}{\line(0,1){8}}
\multiput(32,31)(2,0){4}{\line(0,1){8}}
\multiput(42,41)(2,0){4}{\line(0,1){8}}
\multiput(72,61)(6,0){2}{\line(0,1){8}}
\multiput(74,61)(2,0){2}{\line(0,1){2}}
\multiput(74,67)(2,0){2}{\line(0,1){2}}
\multiput(82,81)(2,0){4}{\line(0,1){8}}
\multiput(92,101)(2,0){4}{\line(0,1){8}}
\end{picture}
\caption{A Young diagram, its corners and a hook walk path}
\end{center}
\end{figure}
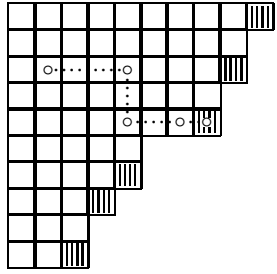

Now consider a sequence
$\lambda_n:\lambda_n(1)\ge\lambda_n(2)\ge\lambda_n(3)\ge \ldots $ of
Young diagrams for which, under suitable scaling, the shape converges
to some limiting shape described by a continuous function. More
precisely, let $f_{\lambda_n}(x)$ be as in \eqref{eq:lambdagraph},
and let $g_{\lambda_n}$ be its rotated coordinate version. Let
$f_\infty:[0,\infty)\to[0,\infty)$ be a weakly decreasing function,
and let $g_\infty$ be its rotated coordinate version, a 1-Lipschitz
function. In this more general setting, think of $g_{\lambda_n}$ and
$g_\infty$ as functions defined on all $\mathbb{R}$. Assume that:
there exists an $M>0$ such that $f_\infty(x) = 0$ for $x\ge M$, and on
$[0,M]$ $f$ is twice continuously differentiable, and its derivative
is bounded away from $0$ and $\infty$ (equivalently, for some $K<0<K'$,
$g_\infty(u)=|u|$ for $u \notin(K,K')$, and $g$ is twice continuously
differentiable in $[K,K']$ with derivative bounded awaw from -1 and
1). Finally, assume that
$$ || g_{\lambda_n} - g_\infty ||_\infty \xrightarrow[n\to\infty]{}
0. $$
For any $n$, let $(I_n,J_n)$ be a
$\mu_{\lambda_n}$-distributed random vector. Let $X_n = I_n/n,
Y_n=J_n/n$. We paraphrase results from \cite{romik}.

\paragraph{Theorem 10.} (Romik \cite{romik}, Theorems 1(b), 6) As
$n\to\infty$, $(X_n,Y_n)$ converges in distribution to the random
vector
$$ (X,Y) := \left( \frac{V+U}{2}, \frac{V-U}{2} \right), $$ where
$V=g_\infty(U)$ and $U$ is a random variable on $[K,K']$ with density function
\begin{equation}\label{eq:cotransition2}
\phi_U(x) = \frac{2}{\pi A} \cos\left( \frac{\pi
g_\infty'(x)}{2}\right) \sqrt{(x-K)(K'-x)} \exp\left(\frac{1}{2}
\int_K^{K'} \frac{g_\infty'(u)}{x-u}du \right),
\end{equation}
with
$$ A = \int_0^M f_\infty(x)dx = \int_K^{K'} (g_\infty(u)-|u|)du $$
and the integral in the exponential being a principal value integral.

\paragraph{Proof of Theorem 2.} We may assume that $0<\alpha<1/2$. The
proof of Theorem 2 now consists of an observation, a remark, and a
computation.

The observation is that the distribution of the location of the
$k_n$-th entry in a random tableau $T \in {\cal T}_n$ is the
distribution of the maximal entry in the shape $\lambda_T^{k_n}$ of
the $k_n$-th subtableau of $T$. Because by Theorem 8, this shape
(suitably rescaled and rotated) converges in probability to
$\tilde{g}_\alpha$ (Theorem 2 assumes $k_n/n^2 \to \alpha$), we may
apply Theorem 10 and conclude that Theorem 2 is true with a density for
$U_\alpha$ given by taking $g_\infty = \tilde{g}_\alpha$, $A=\alpha$,
$-K=K'=\sqrt{2\alpha(1-\alpha)}$ in \eqref{eq:cotransition2}.

The remark is that the above is not quite true, since
$\tilde{g}_\alpha$ does not satisfy the assumptions of Theorem 10! The
problem is that 
$$ -\lim_{\epsilon \searrow 0}
\tilde{g}_\alpha'(-\sqrt{2\alpha(1-\alpha)}+\epsilon) = \lim_{\epsilon
\searrow 0} \tilde{g}_\alpha'(\sqrt{2\alpha(1-\alpha)}-\epsilon) = 1,
$$ 
so the derivative is not bounded away from -1 and 1. However, since
this only happens near the two boundary points, going over the
computations in \cite{romik} shows that this is not a problem, and the
formula \eqref{eq:cotransition2} is still valid in this case
\footnote{ Alternatively, one may use the less explicit formula (8)
from \cite{romik}, which is valid even without the assumption that
$g_\infty'$ is bounded away from $\pm 1$, to verify directly that
\eqref{eq:semicircle} is the cotransition measure of
$\tilde{g}_\alpha$. }.

The computation is the verification that \eqref{eq:cotransition2}
gives the semicircle distribution \eqref{eq:semicircle} under the
above substitutions. We compute, using \eqref{eq:diffmultline} and the
identity $\cos(\arct v) = (1+v^2)^{-1/2}$:
\begin{eqnarray*}
\frac{2}{\pi A} &=& \frac{2}{\pi \alpha} \\
\sqrt{(x-K)(K'-x)} &=& \sqrt{2\alpha(1-\alpha)-x^2} \\
\exp\left(\frac{1}{2}
\int_K^{K'} \frac{\tilde{g}_\alpha(u)}{x-u}du \right) &=&
\sqrt{\frac{\alpha}{1-\alpha}} \cdot
\sqrt{\frac{\frac{1}{2}-x^2}{2\alpha(1-\alpha)-x^2}}, \\
\cos\left( \frac{\pi \tilde{g}_\alpha'(x)}{2}\right) &=& \cos
\left( \arct \frac{(1-2\alpha)x}{\sqrt{2\alpha(1-\alpha)-x^2}} \right)
\\ &=&
\left(1+\frac{(1-4\alpha(1-\alpha))x^2}{2\alpha(1-\alpha)-x^2}
\right)^{-1/2} =
\frac{\sqrt{2\alpha(1-\alpha)-x^2}}{2\sqrt{\alpha(1-\alpha)}
\sqrt{\frac{1}{2}-x^2}} 
\end{eqnarray*}
Multiplying the above expressions gives
$$ \phi_U(x) = \frac{1}{\pi\alpha(1-\alpha)}
\sqrt{2\alpha(1-\alpha)-x^2}, \qquad |x|\le \sqrt{2\alpha(1-\alpha)},
$$
as claimed.
\qed

\section{Plane partitions}

We prove Theorem 4 on the limit shape of plane partitions of an
integer $m$ over an $n\times n$ square diagram, when $n^6 = o(m)$. The
basic observation relating this to Young tableaux is that in this
asymptotic regime, almost all plane partitions have distinct parts. A
plane partition with distinct parts can be completely described by
separately giving the order structure on its parts -- a square Young
tableau -- and an unordered list of the parts, which is simply a
linear partition of $m$ into $n^2$ distinct parts. The structure of
these linear partitions is described by a limit shape theorem due to
Vershik and Yakubovich \cite{vershikyak1},
\cite{vershikyak2}. Combining these results will give us our proof of
Theorem 4.

We will use a result of Erd\"os and Lehner on partitions into a fixed
number of summands.

\paragraph{Theorem 11. (Erd\"os-Lehner \cite{erdoslehner})}
Let $p(m,k)$ denote the number of partitions of $m$ into $k$
parts. Let $q(n,k)$ denote the number of partitions of $m$ into $k$
\emph{distinct} parts. If $m$ and $k$ are sequences of integers that
tend to infinity in such a way that $k^3=o(m)$, then
$$ \frac{q(m,k)}{p(m,k)} \xrightarrow[\qquad]{} 1. $$
In words, if $k^3=o(m)$, almost all partitions of $m$ into $k$ parts
have no repeated parts.

\paragraph{Proof.} This is a combination of Corollary 4.3 and Lemma
4.4 in \cite{erdoslehner}.
\qed

\bigskip
For a Young diagram $\lambda$, denote by $p_\lambda(m)$ the number of
plane partitions of $m$ of shape $\lambda$. Denote by $q_\lambda(m)$
the number of plane partitions of $m$ of shape $\lambda$ with all
parts distinct.

\paragraph{Lemma 11.} Let $\lambda_m$ be a sequence of Young
diagrams. Let $k_m = |\lambda_m|$. If $k_m^3 = o(m)$ as $m\to\infty$,
then
$$ \frac{q_{\lambda_m}(m)}{p_{\lambda_m}(m)}
\xrightarrow[m\to\infty]{} 1. $$
In words, if $k_m^3 = o(m)$, almost all plane partitions of shape
$\lambda_m$ have no repeated parts.

\paragraph{Proof.} If $\lambda$ is a Young diagram of size
$k=|\lambda|$, a plane partition of $m$ of shape $\lambda$ is
described by the order structure on its parts, and the unordered set
of the parts. This gives the equation
$$ q_\lambda(m) = d(\lambda) q(m,k). $$
We claim that
\begin{equation} \label{eq:ppineq}
p_\lambda(m) \le d(\lambda) p(m,k).
\end{equation}
This will prove the claim, since then we will have
$$ \frac{q(m,k_m)}{p(m,k_m)} \le
\frac{q_{\lambda_m}(m)}{p_{\lambda_m}(m)} \le 1, $$ and the Lemma will
follow from Theorem 11. To prove \eqref{eq:ppineq}, we define a
mapping that assigns injectively to each plane partition
$\pi=(p_{i,j})_{(i,j)\in \lambda}$ a pair $(T,\mu)$, where
$T=(t_{i,j})_{(i,j)\in\lambda}$ is a Young tableau of shape $\lambda$,
and $\mu:\mu(1)\ge \mu(2)\ge \ldots \ge\mu(k)$ is a partition of $m$ into
$k$ parts. The mapping is defined as follows. Define a linear order
``$\prec$'' on the cells $(i,j)$ of $\lambda$, by stipulating that
$$ (i,j) \prec (i',j') \iff
p_{i,j}>p_{i',j'}\textrm{ or }
\bigg[p_{i,j}=p_{i',j'}\textrm{ and }
  \big(i<i'\textrm{ or }
    (i=i'\textrm{ and }j<j')\big)\bigg].
$$
Let $(i_1,j_1) \prec (i_2,j_2) \prec \ldots \prec (i_k,j_k)$ be the cells
of $\lambda$ sorted in this linear ordering, and set $t_{i_l,j_l}=l$
and $\mu(l) = p_{i_l,j_l}$, \ $l=1,2,\ldots,k$.

It is easy to verify that the mapping is injective and has the
required range. See Figure 7 for an illustration.
\qed

\begin{figure}[h!]
\begin{center}
$\pi =$\ 
\begin{tabular}{lllll}
7 & 7 & 6 & 5 & 2 \\
7 & 6 & 5 & 5 &   \\
7 & 5 & 2 &   &   \\
6 &   &   &   &   
\end{tabular}
$\xrightarrow[\qquad]{}
\begin{array}{l}
T =\ 
\begin{tabular}{lllll}
1 & 2 & 5 & 8 & 12 \\
3 & 6 & 9 & 10 &   \\
4 & 11 & 13 &   &   \\
7 &   &   &   &   
\end{tabular} \\ \ \\
\mu = 7,\,7,\, 7,\,7,\,6,\,6,\,6,\,5,\,5,\,5,\,2,\,2.
\end{array}$
\caption{Illustration of the proof of Lemma 11}
\end{center}
\end{figure}

\bigskip
Next, we recall the Vershik-Yakubovich limit shape theorem for
partitions of $m$ into $k$ distinct summands, when $k=o(\sqrt{m})$.

\paragraph{Theorem 12. (Vershik-Yakubovich \cite{vershikyak1},
\cite{vershikyak2})} Let $m=m_n$ and $k=k_n$ grow to infinity
as a function of some parameter $n$, in such a way that
$k=o(\sqrt{m})$. Let $\lambda_n: \lambda_n(1)> \lambda_n(2) >
\ldots >\lambda_n(k)$ be a sequence of uniform random partitions of $m$
into $k$ distinct parts. Then for any $t\ge 0,\ \epsilon>0$,
$$
\mathbb{P}\left( \left|\frac{1}{k}\# \{ 1 \le l \le k : \lambda_n(l) >
\frac{m}{k} t \} - e^{-t}\right|>\epsilon \right)
\xrightarrow[n\to\infty]{} 0.
$$
In words, the graph of the Young diagram of a uniform random partition
of $m$ into $k$ distinct parts, when $k=o(\sqrt{m})$, will with high
probability resemble the limit shape $e^{-t}$.\footnote{Actually, it
is more correct to say that this is the graph of the conjugate
partition $\lambda'$.}

\paragraph{Proof of Theorem 4.} Let $\pi=(p_{i,j})_{i,j=1}^n$ be the
random plane partition of $m$ over the square diagram
$\square_n$. Since $n^6=o(m)$ and $|\square_n|=n^2$, by Lemma 11 we may
assume that $\pi$ was chosen uniformly among all plane partitions of
$m$ of shape $\square_n$ \emph{with all parts distinct}, since this is
a set of probability close to $1$ in ${\cal P}_{n,m}$. Equivalently,
by the remarks at the beginning of this section, we may assume that
$\pi$ is selected by choosing independently a random Young tableau $T
\in {\cal T}_n$ and a random partition
$\mu:\mu(1)>\mu(2)>\ldots>\mu(n^2)$ of $m$ into $n^2$ distinct parts,
then setting $p_{i,j} = \mu(t_{i,j})$.

Fix $0\le x,y<1$. Let $\alpha=L(x,y)$. Let $i=\lfloor n x \rfloor+1$,
$j=\lfloor n y \rfloor + 1$, and $\beta = t_{i,j}/n^2$.

We need to show that $\tilde{S}_\pi(x,y) = (n^2/m)p_{i,j}$ is with
high probability very close to $-\log \alpha$. Let $\epsilon>0$ be
small. From Theorem 1(i), we know that with (asymptotically) high
probability
\begin{equation}\label{eq:betapp}
|\beta-\alpha| < \epsilon.
\end{equation}
Now, apply Theorem 12 for the random partition $\mu$ with
$t=-\log(\alpha-2\epsilon)$. This gives that with high probability
$$ \left| \frac{1}{n^2}\# \left\{ 1 \le l \le n^2 : \mu(l) >
\frac{m}{n^2}(-\log(\alpha-2\epsilon))\right\} - (\alpha-2\epsilon) \right| <
\epsilon, $$
or equivalently, since $\mu(1)>\mu(2)>\ldots>\mu(n^2)$,
$$ \left| \frac{1}{n^2} \max\left\{ 1 \le l \le n^2 : \mu(l) >
\frac{m}{n^2}(-\log(\alpha-2\epsilon))\right\} - (\alpha-2\epsilon) \right| <
\epsilon. $$
This implies in particular that
$$ \max\left\{ 1 \le l \le n^2 : \mu(l) >
\frac{m}{n^2}(-\log(\alpha-2\epsilon))\right\} <
n^2(\alpha-2\epsilon+\epsilon) = (\alpha-\epsilon)n^2, $$
hence, since by \eqref{eq:betapp}, $t_{i,j} > (\alpha-\epsilon)n^2$,
$$ p_{i,j} = \mu(t_{i,j}) \le \frac{m}{n^2}(-\log(\alpha-2\epsilon)). $$
Apply Theorem 12 again with
$t=-\log(\alpha+2\epsilon)$. This gives that with high probability
$$ \left| \frac{1}{n^2}\# \left\{ 1 \le l \le n^2 : \mu(l) >
\frac{m}{n^2}(-\log(\alpha+2\epsilon))\right\} - (\alpha+2\epsilon) \right| <
\epsilon, $$
or equivalently
$$ \left| \frac{1}{n^2} \max\left\{ 1 \le l \le n^2 : \mu(l) >
\frac{m}{n^2}(-\log(\alpha+2\epsilon))\right\} - (\alpha+2\epsilon) \right| <
\epsilon. $$
In particular this gives that
$$ \max\left\{ 1 \le l \le n^2 : \mu(l) >
\frac{m}{n^2}(-\log(\alpha+2\epsilon))\right\} >
n^2(\alpha+2\epsilon-\epsilon) = (\alpha+\epsilon)n^2, $$
hence, since by \eqref{eq:betapp}, $t_{i,j} < (\alpha+\epsilon)n^2$,
$$ p_{i,j} = \mu(t_{i,j}) > \frac{m}{n^2}(-\log(\alpha+2\epsilon)). $$
We have shown that the event
$$
-\log(\alpha+2\epsilon) < \frac{n^2}{m} p_{i,j} \le -\log(\alpha-2
\epsilon)
$$
holds with asymptotically high probability. Since $\epsilon$ was
arbitrary the result follows.
\qed

\section{Computations for the rectangular case}

The proof of Theorem 5 involves exactly the same ideas as the proof of
Theorem 1, with some more computations, which we include here for
completeness. The proof that Theorem 6 follows from Theorem 5 is
completely identical to the proof in section 6 that Theorem 4 follows
from Theorem 1.

Fix $0<\theta\le 1$ and $0<\alpha<1$. Our starting point is the
rotated-coordinate formulation of the variational problem whose
solution will yield the $\alpha$-level curve of the limit surface
$L_\theta$. The computations leading to this variational problem are
obvious generalizations of the corresponding computations for the
square case $\theta=1$, and are omitted.

\paragraph{Variational problem - the rectangular case.} A function
$h:[-\theta \sqrt{2}/2, \sqrt{2}/2]\to[0,\infty)$ is called
$\alpha$-admissible if $h$ is $1$-Lipschitz, and satisfies
\begin{eqnarray} \label{eq:admiss1}
 h(-\theta\sqrt{2}/2) &=& \theta\sqrt{2}/2, \\
\label{eq:admiss2}
h(\sqrt{2}/2) &=& \sqrt{2}/2, \\
\label{eq:admiss3}
\int_{-\theta\sqrt{2}/2}^{\sqrt{2}/2}(h(u)-|u|)du &=& \alpha \theta.
\end{eqnarray}
Find the unique $\alpha$-admissible $h$ that minimizes
$$ J(h) = -\frac{1}{2} \int_{-\theta\sqrt{2}/2}^{\sqrt{2}/2}
\int_{-\theta\sqrt{2}/2}^{\sqrt{2}/2} h'(s)h'(t)\log|s-t|ds\,dt.
$$

\bigskip
To derive the minimizer, first consider the case when $\alpha$ is
small. In that case, we make an assumption on the form of the
minimizer similar to \eqref{eq:guess}, but with a non-symmetric
interval $[-\beta_1(\alpha),\beta_2(\alpha)]$, where $\beta_1\in
(0,\theta\sqrt{2}/2)$, $\beta_2\in(0,\sqrt{2}/2)$. That is, we assume
that $h'$ has the form
\begin{equation}\label{eq:guessrect}
h'(s)\ \ \ \textrm{  is }\left\{ \begin{array}{ll}
  = -1, & \textrm{if }-\theta\sqrt{2}/2 < s < -\beta_1, \\
  \in (-1,1), & \textrm{if }-\beta_1<s<\beta_2, \\
  = +1, & \textrm{if }\beta_2<s<\sqrt{2}/2. \end{array}\right.
\end{equation}
Replace the conditions \eqref{eq:admiss1}, \eqref{eq:admiss2},
\eqref{eq:admiss3} with the equivalent set of conditions
\begin{eqnarray} \label{eq:adm1}
 h(-\theta\sqrt{2}/2) &=& \theta\sqrt{2}/2, \\
\label{eq:adm2}
\int_{-\theta\sqrt{2}/2}^{\sqrt{2}/2}h'(u)du &=& (1-\theta)\sqrt{2}/2,\\
\label{eq:adm3}
\int_{-\theta\sqrt{2}/2}^{\sqrt{2}/2}uh'(u)du +
\frac{1+\theta^2}{4}&=&\theta \alpha.
\end{eqnarray}
(In the square case $\theta=1$ we did not impose the restriction
\eqref{eq:adm2} on $h$ as we expected $h$ to be even, i.e. $h'$ to be
odd, so that the condition \eqref{eq:adm2} would be met
automatically. Not anymore in the rectangular case!) Then the
counterpart to \eqref{eq:multline} is: for $s \in (-\beta_1,\beta_2)$,
\begin{multline}\label{eq:rect1}
-\int_{-\beta_1}^{\beta_2} h'(t)\log|s-t|dt - \lambda s - \mu \\
+ (s+\tsq)\log|s+\tsq|+(s-\tsq)\log|s-\tsq| \\
- (s+\bo)\log|s+\bo|-(s-\bt)\log|s-\bt| \qquad\qquad\qquad\qquad \\
+\bo-\bt+(1-\theta)\sq
= 0. \qquad\qquad\qquad
\end{multline}
Here $\lambda$, $\mu$ are the Lagrangian multipliers dual to the
constraints \eqref{eq:adm3} and \eqref{eq:adm2}
respectively. Differentiating \eqref{eq:rect1} with respect to $s$
gives
\begin{equation}\label{eq:rect2}
-\int_{-\bo}^{\bt} \frac{h'(t)}{s-t}dt = \lambda + \log
 \frac{s+\bo}{s+\tsq} + \log \frac{\bt-s}{\sq-s}, \qquad
 s\in(-\bo,\bt).
\end{equation}
Introduce $a=(\bo+\bt)/2$, $b=(\bt-\bo)/2$, and substitute $s=a\xi+b,
t=a\eta+b$. The above equation becomes
\begin{equation}\label{eq:rect3}
-\int_{-1}^1 \frac{h'(a\eta+b)}{\xi-\eta}d\eta = \lambda +
 \log\frac{1+\xi}{\gamo+\xi} + \log\frac{1-\xi}{\gamt-\xi};
\end{equation}
here
\begin{equation}\label{eq:gammadef}
\gamo=\frac{\bt-\bo+\theta\sqrt{2}}{\bo+\bt}, \qquad
\gamt=\frac{\bo-\bt+\sqrt{2}}{\bo+\bt},
\end{equation}
and it is easy to check that $\gamo,\gamt>1$. Applying Theorem 9 to
\eqref{eq:rect3} and using Lemma 8, we obtain
\begin{multline}\label{eq:hprimerect}
h'(a\xi+b) =
\frac{1}{\pi^2(1-\xi^2)^{1/2}}(\pi\lambda\xi+I(\xi,\gamo)-I(-\xi,\gamt))
+ \frac{c'}{(1-\xi^2)^{1/2}}
\\ \quad
= \frac{\xi}{\pi(1-\xi^2)^{1/2}}\left(
\lambda-\log\left(\gamo+\sqrt{\gamo^2-1}\right)
       -\log\left(\gamt+\sqrt{\gamt^2-1}\right) \right) \\
+ \frac{2}{\pi}\left[
\arct\sqrt{\frac{(1+\xi)(\gamt-1)}{(1-\xi)(\gamt+1)}}
- \arct\sqrt{\frac{(1-\xi)(\gamo-1)}{(1+\xi)(\gamo+1)}} \right]
+ c(1-\xi^2)^{-1/2},
\end{multline}
$c', c$ being arbitrary constants. As in the symmetric case, if $h'(s)$
is to be bounded for $s\in(-\bo,\bt)$ (i.e. for $\xi\in(-1,1)$),
necessarily
\begin{equation}\label{eq:lambdarect}
c=0, \qquad
\lambda = \log\left(\gamo+\sqrt{\gamo^2-1}\right)
          +\log\left(\gamt+\sqrt{\gamt^2-1}\right).
\end{equation}
So we have
\begin{equation}\label{eq:hpr}
h'(a\xi+b) = 
\frac{2}{\pi}\left[
\arct\sqrt{\frac{(1+\xi)(\gamt-1)}{(1-\xi)(\gamt+1)}}
- \arct\sqrt{\frac{(1-\xi)(\gamo-1)}{(1+\xi)(\gamo+1)}} \right],
\end{equation}
for which indeed $|h'(a\xi+b)|\le 1$ holds.

We still need to find $\bo$ and $\bt$. Using \eqref{eq:guessrect},
rewrite \eqref{eq:adm2} and \eqref{eq:adm3} as, respectively,
\begin{eqnarray}
\label{eq:adm2new}
\int_{-\bo}^{\bt} h'(t)dt &=& \bt-\bo, \\
\label{eq:adm3new}
\int_{-\bo}^{\bt} th'(t)dt + \frac{1}{2}(\bo^2+\bt^2) &=&
\theta\alpha.
\end{eqnarray}
Now evaluating these integrals using \eqref{eq:hpr}, this gives the
equations
\begin{eqnarray*}
a\left[\left( \sqrt{\gamt^2-1}-\gamt\right)-
\left( \sqrt{\gamo^2-1}-\gamo\right)\right] &=& \bt-\bo, \\
-\frac{a^2}{2}\sum_{i=1}^2 \left[1-\gamma_i^2+\gamma_i
\sqrt{\gamma_i^2-1}\right] - b(\bt-\bo)+\frac{1}{2}(\bo^2+\bt^2) &=&
\theta\alpha,
\end{eqnarray*}
($a=(\bo+\bt)/2$, $b=(\bt-\bo)/2$). Excluding $\bo, \bt$ via
\eqref{eq:gammadef}, we obtain two equations for $\gamo, \gamt$,
namely
\begin{eqnarray}
\label{eq:eqgam1}
\gamo-\theta\gamt = \frac{1+\theta}{2}\left[
((\gamt^2-1)^{1/2}-\gamt)-((\gamo^2-1)^{1/2}-\gamo) \right], \\
\label{eq:eqgam2}
\theta\alpha =
\frac{(1+\theta)^2-(\gamo-\theta\gamt)^2}{2(\gamo+\gamt)^2} -
\frac{(1+\theta)^2}{4(\gamo+\gamt)^2}
\sum_{i=1}^2 \left[1-\gamma_i^2+\gamma_i
\sqrt{\gamma_i^2-1}\right].\ 
\end{eqnarray}
It seems a minor miracle that these equations can be solved
explicitly. Here is how. Isolating the difference of the radicals in
the first equation, multiplying both sides of the resulting equation
by the sum of radicals and cancelling the common factor $\gamo+\gamt$,
we obtain
$$
\sqrt{\gamo^2-1}+\sqrt{\gamt^2-1}=\frac{1+\theta}{1-\theta}(\gamt-\gamo).
$$
(In particular, $\gamt>\gamo$.) Combining this with the initial
equation, we express the radicals as linear combinations of $\gamo,
\gamt$:
\begin{eqnarray}\label{eq:radical1}
\sqrt{\gamo^2-1}&=&-\frac{1+\theta^2}{1-\theta^2}\gamo +
\frac{2\theta}{1-\theta^2}\gamt, \\
\label{eq:radical2}
\sqrt{\gamt^2-1}&=&-\frac{2\theta}{1-\theta^2}\gamo +
\frac{1+\theta^2}{1-\theta^2}\gamt.
\end{eqnarray}
Plugging these expressions for the radicals into \eqref{eq:eqgam2},
after collecting like terms, we obtain a \emph{quadratic} equation for
$x=\gamt/\gamo$:
$$
x^2(\theta(1-\alpha)+\alpha)-x(1-\theta)(1-2\alpha)-(1-\alpha+\theta\alpha)
= 0.$$
Consequently
$$ \frac{\gamt}{\gamo} = x =
\frac{\theta\alpha+1-\alpha}{\alpha+\theta(1-\alpha)}, $$
which, for $\alpha<1/2$, exceeds 1. Squaring both sides of
\eqref{eq:radical1} and substituting $\gamt=x\gamo$, we solve for
$\gamo$ to obtain
\begin{equation} \label{eq:eqgamnew}
\gamo =
\frac{\alpha+\theta(1-\alpha)}{2\sqrt{\theta\alpha(1-\alpha)}},
\qquad 
\gamt = \frac{\theta\alpha+1-\alpha}{2\sqrt{\theta\alpha(1-\alpha)}}.
\end{equation}
Direct checking reveals that these $\gamo, \gamt$ satisfy the
equations \eqref{eq:radical1}, \eqref{eq:radical2} themselves as long as
\begin{equation}\label{eq:condalphastar}
\alpha \le \alpha^* := \frac{\theta}{1+\theta}.
\end{equation}
For $\alpha>\alpha^*$, the gammas do not satisfy
\eqref{eq:radical1}. More precisely, $\gamo, \gamt$ would have satisfied
this equation, had we considered the negative value of
$\sqrt{\gamo^2-1}$. However, we need the positive value only. The
source of the trouble here is that $\gamo=1$ for
$\alpha=\alpha^*$. Tellingly, the boundary point $(-\bo,h(-\bo))$
reaches the corner $(-\tsq,\tsq)$ of the rotated rectangle at
$\alpha=\alpha^*$. Using \eqref{eq:gammadef}, we obtain: for $\alpha
\le \alpha^*$,
\begin{multline}\label{eq:betabeta}
\qquad\qquad\qquad\qquad\quad\ \ 
\bo = \sqrt{2\theta\alpha(1-\alpha)} - \alpha(1-\theta)\sq, \\
\bt = \sqrt{2\theta\alpha(1-\alpha)} + \alpha(1-\theta)\sq.
\qquad\qquad\qquad\qquad\qquad
\end{multline}
Using \eqref{eq:eqgamnew} and \eqref{eq:betabeta}, we simplify
\eqref{eq:hpr} to
\begin{equation}\label{eq:hprnew}
h'(t) = \frac{2}{\pi}\arct \left[
\frac{(1-\theta)\sqrt{\alpha(1-\alpha)} + \xi\sqrt{\theta}(1-2\alpha)}
{\sqrt{\theta(1-\xi^2)}}\right],
\end{equation}
where
$$ \xi = \frac{t-b}{a} =
\frac{t-\alpha(1-\theta)\sq}{\sqrt{2\theta\alpha(1-\alpha)}}, \qquad
t\in [-\beta_1,\beta_2].
$$
Furthermore, using
$$ h(s) = \bo + \int_{-\bo}^s h'(t)dt = \bo + a\int_{-1}^\xi
h'(a\eta+b) d\eta, $$
\eqref{eq:hpr}, and \eqref{eq:eqx4}, we obtain
\begin{multline}\label{eq:formulah}
h(s) = h(a\xi+b) = \bo \\ +\frac{2a}{\pi} \bigg[ -(\xi+\gamo) \arct
\sqrt{\frac{(1-\xi)(\gamo-1)}{(1+\xi)(\gamo+1)}} + (\xi-\gamt) \arct
\sqrt{\frac{(1+\xi)(\gamt-1)}{(1-\xi)(\gamt+1)}} \\
+ \frac{1}{2}\left(\arcs \xi+\frac{\pi}{2}\right)
\left(\sqrt{\gamt^2-1}-\sqrt{\gamo^2-1}\right) +
\frac{\pi}{2}(\gamo-1) \bigg].
\end{multline}
Finally, combining \eqref{eq:lambdarect} and \eqref{eq:eqgamnew}, we
compute
\begin{equation}\label{eq:lambdanew}
\lambda = \log \frac{\theta(1-\alpha)}{\sqrt{\theta\alpha(1-\alpha)}}
+ \log \frac{(1-\alpha)}{\sqrt{\theta\alpha(1-\alpha)}} =
\log\frac{1-\alpha}{\alpha},
\end{equation}
the same value as in the square case!

It remains to consider the range $\alpha^* < \alpha \le 1/2$. In this
case, it turns out that the formulas for the parameters $\bo, \bt$
remain the same, while the formula for the corresponding $h'$ changes
slightly. What is different is that now $h'(s)=1$ for $s\in [-\tsq,
-\bo]$ and
$$ h(-\bo) = \theta\sqrt{2} - \bo. $$
The latter condition means that now the boundary point $[-\bo,
h(-\bo)]$ lies on the longer side of the rectangle. The starting point
now is a modification of \eqref{eq:rect2}, stemming from
$h'(s)\equiv 1$, rather than $-1$, for $s\in [-\tsq,-\bo]$, namely
$$
-\int_{-1}^1 \frac{h'(a\eta+b)}{\xi-\eta}d\eta = \lambda - \log
 \frac{1+\xi}{\gamo+\xi} + \log \frac{1-\xi}{\gamt-\xi}.
$$
This leads to
$$ h'(a\xi+b) =
\frac{1}{\pi^2(1-\xi^2)^{1/2}}(\pi\lambda\xi-I(\xi,\gamo) -
I(-\xi,\gamt)) + \frac{c'}{(1-\xi^2)^{1/2}}, $$
($c'$ being a constant),
compare with the first line in \eqref{eq:hprimerect}, whence to
$$ 
\lambda = \log\left(\gamt+\sqrt{\gamt^2-1}\right)
          -\log\left(\gamo+\sqrt{\gamo^2-1}\right),
$$
compare with \eqref{eq:lambdarect}, and
\begin{equation*}\label{eq:hpr2}
h'(a\xi+b) = 
\frac{2}{\pi}\left[
\arct\sqrt{\frac{(1-\xi)(\gamo-1)}{(1+\xi)(\gamo+1)}}
+ \arct\sqrt{\frac{(1+\xi)(\gamt-1)}{(1-\xi)(\gamt+1)}} \right],
\end{equation*}
compare with \eqref{eq:hpr}. After integration, the final formula for
$h(s)$ is
\begin{multline}\label{eq:formulah2}
h(s) = \theta\sqrt{2} - \bo \\ +\frac{2a}{\pi} \bigg[ (\xi+\gamo) \arct
\sqrt{\frac{(1-\xi)(\gamo-1)}{(1+\xi)(\gamo+1)}} + (\xi-\gamt) \arct
\sqrt{\frac{(1+\xi)(\gamt-1)}{(1-\xi)(\gamt+1)}} \\
+ \frac{1}{2}\left(\arcs \xi+\frac{\pi}{2}\right)
\left(\sqrt{\gamo^2-1}+\sqrt{\gamt^2-1}\right) +
\frac{\pi}{2}(1-\gamo) \bigg].
\end{multline}
We add that,
despite the difference between the two formulas for $\lambda$ -- the
one above for $\alpha \ge \alpha^*$ and \eqref{eq:lambdarect} for
$\alpha\le \alpha^*$ -- the eventual expression is still that in
\eqref{eq:lambdanew}. The ``secret'' is that
$$ \sqrt{\gamo^2-1} =
\frac{|\sqrt{\alpha}-\sqrt{\theta(1-\alpha)}|}{2\sqrt{\theta\alpha(1-\alpha)}},
$$
with $\sqrt{\alpha}-\sqrt{\theta(1-\alpha)}$ changing its sign at
$\alpha^*$.

It remains to prove that $h$ is indeed a minimizer. Let $\alpha <
\alpha^*$. Consider $ s \in [-\tsq,-\bo]$. Since $h'(s)=-1$, we need
to show that $F(s,\alpha)\ge 0$, where $F(s,\alpha)$ is the left-hand
side expression in \eqref{eq:rect1}. The above computations show that
$F(s,\alpha) \equiv 0$ for $s\in (-\bo,\bt)$. As in the square case,
for fixed $s\in [-\tsq,\bo]$ let $\hat{\alpha}\in (\alpha,\alpha^*)$
be defined by $s=-\bo(\hat{\alpha})$. Then $F(s,\hat{\alpha})=0$,
and we will prove $F(s,\alpha)\ge 0$ if we show that $\partial
F(s,x)/\partial x < 0$ for all $x\in [\alpha,\hat{\alpha}]$. Since
$0\in (-\bo(x),\bt(x))$, $F(0,x)=0$, and we use the latter equation to
exclude the Lagrangian multiplier $\mu$ in the expression for
$F(s,x)$. Then, an easy computation shows that
\begin{equation}\label{eq:easycomp}
\frac{\partial F(s,x)}{\partial x} = -\int_{-\bo}^{\bt} \frac{\partial
h_x(t)}{\partial x} \log|s-t|dt + \int_{-\bo}^{\bt} \frac{\partial
h_x(t)}{\partial x} \log|t|dt -s \frac{d\lambda}{dx},
\end{equation}
where $\beta_i = \beta_i(x), \lambda=\lambda(x)$ are given by
\eqref{eq:betabeta} and \eqref{eq:lambdanew}. Let us evaluate
$\partial h_x(s)/\partial x$ for $s\in (-\bo,\bt)$. Differentiating
\eqref{eq:rect2} with respect to $x$ we obtain
$$ -\int_{-\bo}^{\bt} \frac{\partial h_x(t)/\partial x}{s-t}dt =
\frac{d\lambda}{dx} = - \frac{1}{x(1-x)}.
$$
Then, using Theorem 9 and \eqref{eq:lambdaterm},
\begin{equation}\label{eq:using}
\frac{\partial h_x(s)}{\partial x} = -\frac{\xi}{\pi
x(1-x)(1-\xi^2)^{1/2}} + \frac{c}{(1-\xi^2)^{1/2}}.
\end{equation}
Here we have to set $c=0$, as the equation \eqref{eq:adm2new} -- upon
differentiation with respect to $x$ -- leads to
\begin{equation}\label{eq:leadsto}
\int_{-1}^1 \frac{\partial h_x(t)}{\partial x} d\xi = 0.
\end{equation}
Hence
\begin{equation}\label{eq:hence}
\frac{\partial h_x(s)}{\partial x} = -\frac{\xi}{\pi
x(1-x)(1-\xi^2)^{1/2}}.
\end{equation}
Plugging this expression into \eqref{eq:easycomp}, integrating by
parts, and using \eqref{eq:unnumbered}, we transform
\eqref{eq:easycomp} into
\begin{equation}\label{eq:transform}
\frac{\partial F(s,x)}{\partial x} =
-\frac{\sqrt{(s-b)^2-a^2}}{x(1-x)} < 0.
\end{equation}
Let $\alpha \in (\alpha^*,1/2]$, and $s<-\bo(\alpha)$ again. Since now
$h'(s)=1$, we need to show that $F(s,\alpha)\le 0$. Let
$\tilde{\alpha} \in (\alpha^*,\alpha)$ be defined by $-s =
\bo(\tilde{\alpha})$. ($\tilde{\alpha}$ exists, uniquely, because
$\bo(x)$ is decreasing on $(\alpha,\alpha^*)$ and $-s < \bo(\alpha)$.)
Then $F(s,\tilde{\alpha})=0$, and so again we need to show that
$\partial F(s,x)/\partial x < 0$. The formula \eqref{eq:using}
continues to hold, and so does \eqref{eq:leadsto}, since now we have
$$ \int_{-\bo}^{\bt} h'(t)dt = \bt + \bo - \tsq, $$
and $h'(-\bo)=1$. Therefore \eqref{eq:hence} remains valid, which
implies \eqref{eq:transform}.

Analogously, $F(s,\alpha)\le 0$ for $s\ge \bt(\alpha)$ and $\alpha \in
(0,1/2]$. This finishes the proof that $h$ is the minimizer, the claim
which forms the core of the proof of Theorem 5.
\qed

\section{Discussion}

\subsection{Plancherel measure}

Let ${\cal Y}_k$ denote the set of Young diagrams of area $k$. The
Plancherel measures are the family of probability measures $\mu_k$ on
${\cal Y}_k$, defined by
\begin{equation}\label{eq:plancherel}
\mu_k(\lambda) = \frac{d(\lambda)^2}{k!}, \qquad (\lambda \in {\cal Y}_k).
\end{equation}
Alternatively, $\mu_k$ is sometimes defined as a measure on all Young
tableaux of size $k$, where
$$ \mu_k(T) = \frac{d(\textrm{shape}(T))}{k!}. $$ The measure on
diagrams is then the projection of the measure on tableaux under the
mapping that assigns to each tableau its shape. The measures $\mu_k$
are a \emph{projective family} of measures, in the following sense: If
$T$ is a $\mu_k$-random tableau, then the tableau $T'$ of size $k-1$
obtained by deleting the $k$-th entry from $T$ is a $\mu_{k-1}$-random
tableau. Therefore, all the $\mu_k$'s can be encompassed by a single
object $\mathbb{P}$, the \emph{infinite Plancherel measure}, which is
a measure on infinite tableaux -- i.e. fillings of the squares in the
positive quadrant of the plane with the positive integers that are
increasing along rows and columns -- for which the marginal
distribution of the shape of the $k$-th subtableau (the set of squares
where the entry of the infinite tableau is $\le k$) is given by
\eqref{eq:plancherel}. In other words, $\mathbb{P}$ can be thought of
as a measure on all sequences $\emptyset = \lambda_0 \subset \lambda_1
\subset \lambda_2 \subset ... $ of Young diagrams, where $\lambda_k$
has size $k$ and is obtained from $\lambda_{k-1}$ by the addition of a
box. So $\mathbb{P}$ is simply a natural Markovian coupling of the
measures \eqref{eq:plancherel}, known sometimes as the
\emph{Plancherel growth process}.

Much is known about Plancherel measure. It arises naturally in
representation theory, as a natural measure on the irreducible
representations of the symmetric group, and in combinatorics, as the
distribution of the output of the RSK algorithm applied to a uniform
random permutation in $S_k$. In particular, the length of the first
row $\lambda_k(1)$ of a $\mu_k$-random Young diagram has the same
distribution as the length $l_n(\pi)$ of the longest increasing
subsequence of a uniform random permutation $\pi$ in $S_k$, an important
permutation statistic.

Logan-Shepp \cite{loganshepp} and Vershik-Kerov \cite{vershikkerov1},
\cite{vershikkerov2} proved that the graph of a $\mu_k$-random Young
diagram, when rescaled by a factor of $\sqrt{k}$ along each axis and
drawn in rotated coordinates, with high probability resembles the
limit shape
$$ \Omega(u) = \left\{ \begin{array}{ll}
\frac{2}{\pi}(u\arcs(u/\sqrt{2})+\sqrt{2-u^2}) & |u| \le \sqrt{2}, \\
|u| & |u| > \sqrt{2}, \end{array} \right.$$
see Figure 8.

\begin{figure}[h!]
\begin{center}
\includegraphics{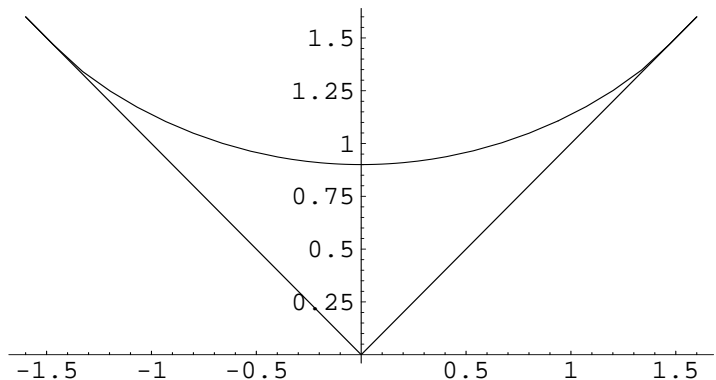}
\caption{The limit shape $v=\Omega(u)$}
\end{center}
\end{figure}
%
Gribov \cite{gribov} noted that this can be reinterpreted as a theorem
on the limit \emph{surface} of the Plancherel-random tableau, much in
the same spirit as Theorem 1. If $T$ is a $\mu_k$-random Young
tableau, then after rescaling the graph of $T$ is approximately
described in rotated coordinates by the surface $\Sigma:D\to[0,1]$,
where
$$ D = \{ (u,v): |u| \le \sqrt{2},\ \, |u| \le v \le \Omega(u) \} $$
is the two dimensional domain bounded between the graphs $|u|$ and
$\Omega(u)$, and for each $0 < \alpha < 1$, the $\alpha$-level curve
of $\Sigma$ is
\begin{equation*}\label{eq:time}
\{ (u,v)\in D : |u| \le \sqrt{2 \alpha},\ \, v =
\sqrt{\alpha}\Omega(u/\sqrt{\alpha}) \},
\end{equation*}
a shrunken copy of $\Omega$.

The approach in the papers of Logan-Shepp and Vershik-Kerov was the
variational approach, of analyzing the limiting integral functional
arising from \eqref{eq:plancherel}. Kerov \cite{kerov} considered the
following more dynamical approach: Assume that we have selected the
$\mu_k$-random diagram $\lambda_k$. Since under $\mathbb{P}$, the
sequence $\lambda_1 \subset \lambda_2 \subset ... $ is a
(nonhomogeneous) Markov chain with values in the Young graph, there is
a measure $\nu$ on the exterior corners of $\lambda_k$ (the boxes in
the complement of $\lambda_k$ that can be added to $\lambda_k$ to form
a Young diagram of size $k+1$), such that if we choose a $\nu$-random
corner of $\lambda_k$ and add the new box there, the resulting diagram
$\lambda_{k+1}$ will have distribution $\mu_{k+1}$. In other words,
$\nu$ is the probability transition measure of the Markov chain
$(\lambda_k)$. It is known as the \emph{transition measure of the
diagram $\lambda_k$}, and is in a sense dual to the co-transition
measure discussed in section 5.

Kerov showed that in the limit when the graph of the diagram
$\lambda_k$ becomes a smooth curve, the transition measure converges
to a limit. Imagine that in the limit, instead of attaching a new box
at a $\nu$-random corner, one attaches a $\nu$-fraction of a box at
each corner. So the curve grows in the ``tangent'' direction given by
$\nu$. Thus, the Plancherel growth process can be described in the
limit as a smooth flow on the (infinite-dimensional) space of
shapes. Kerov showed that $\Omega(u)$ is the unique shape which, after
rescaling, is invariant under this flow, and that this fixed point is
an attractor of the flow; this explains, in a way, (though does not
formally prove) its appearance as the limit shape for
Plancherel-random diagrams. Remarkably, the transition measure of
$\Omega$ (the limiting direction of the flow) is the semicircle
distribution.

Another interesting direction stemming from the study of Plancherel
measure is the connection to longest increasing subsequences of random
permutations. The limit shape result of Logan-Shepp and Vershik-Kerov
implies that the length $l_n(\pi)$ of the longest increasing
subsequence of a random permutation $\pi\in S_n$ is with high
probability at least $(1-o(1))2 \sqrt{n}$. Using additional arguments
(which were an inspiration for our proof of Theorem 3), Vershik and
Kerov showed also that $l_n(\pi)$ is with high probability at most
$(1+o(1))2\sqrt{n}$, solving the so-called Ulam's problem. More
recently, Baik, Deift and Johansson \cite{baiketal} showed that the
fluctuation of $l_n(\pi)$ around its asymptotic value $2\sqrt{n}$ has
a limiting distribution. More precisely,
$$ \frac{l_n(\pi) - 2\sqrt{n}}{n^{1/6}}
\xrightarrow[n\to\infty]{\textrm{ in distribution }} F. $$
Here $F$ is the Tracy-Widom distribution from random matrix theory,
defined as
$$ F(t) = \exp\left(-\int_t^\infty (x-t)u(x)^2 dx \right), $$
where $u(x)$ is the solution of the Painlev\'e II equation
$u''(x)=2u(x)^3 + xu(x)$ that is asymptotic to the Airy function
$\textrm{Ai}(x)$ as $x\to\infty$.
Other results along those lines can be found in \cite{baiketal2},
\cite{borodinetal}, \cite{johansson}, \cite{okounkov}; see also the
survey \cite{aldousdiaconis}

The distribution $F$ appears in random matrix theory as the limiting
distribution of the maximal eigenvalue of a GUE random
matrix. Following the Baik-Deift-Johansson result, it was found that
there are many striking parallels between Plancherel measure and
random matrix ensembles, see \cite{borodinetal}, \cite{oconnell},
\cite{ivanovolsh}. In particular, the transition measure of the
Plancherel-random diagram converges to the semicircle law, which is
also the limiting distribution of the empirical eigenvalue
distribution in the GUE and GOE random matrix ensembles. Ivanov and
Olshanski \cite{ivanovolsh} showed that this similarity is no mere
coincidence, but in fact appears also in the finer fluctuations of the
transition measure and eigenvalue distribution measure around the
semicircle distribution.

\subsection{The random square tableau as a deformation of Plancherel
measure}

The reader familiar with the works of Logan-Shepp and Vershik-Kerov
will undoubtedly have noticed the similarity between these results and
our analysis of the square tableau model. Define for each positive
integer $n$ and each $1 \le k \le n^2$, the probability measure
$\nu_{n,k}$ on ${\cal Y}_k$, by
\begin{equation}\label{eq:deform}
\nu_{n,k}(\lambda) = \frac{d(\lambda)d(\square_n\setminus
\lambda)}{d(\square_n)}, \qquad (\lambda \in {\cal Y}_k),
\end{equation}
where $d(\square_n\setminus\lambda)$ is taken as $0$ if $\lambda
\nsubseteq \square_n$. The measure $\nu_{n,k}$ is the distribution of
the $k$-th subtableau of a random $n\times n$ square tableau, and our
entire approach revolved around the analysis of its properties. It is
remarkable how many of the ideas used in the study of Plancherel
measure we have found useful in our study of square tableaux; first,
and most obviously, the variational problem that arises from
\eqref{eq:deform} resembles the variational problem studied by
Logan-Shepp and Vershik-Kerov. Although our approach in solving the
variational problem relied on the more methodical use of the inversion
formula for Hilbert transforms (an approach that could be applied the
Plancherel case as well!), we were greatly inspired by the methods
used in the Plancherel case. Secondly, our proof of Theorem 3 and the
treatment of the boundary of the square also follows closely the ideas
of Vershik and Kerov (with the notable difference, that our proof of
the upper bound uses the lower bound!). Finally, our Theorem 2 on the
location of particular entries, was inspired by Kerov's differential
model \cite{kerov} for Plancherel growth. By postulating the existence
of an analogous differential growth model for the $k$-subtableaux of
the square tableau, we were able to guess Theorem 2 from the formula
\eqref{eq:partialdiff}. This was later verified by a different method,
using the result from \cite{romik}.

Take another look at \eqref{eq:deform} and \eqref{eq:plancherel}. The
defining equations for $\mu_k$ and $\nu_{n,k}$ seem superficially
similar at best. In fact, they are closely related, and when $k$ is
very small these measures are quite close. To make this precise, we
first note the following curious identity. Define the \emph{falling
power} $a^{\downarrow b}$ of $a$ as $a^{\downarrow
b}=a(a-1)(a-2)\ldots (a-b+1)$. Then:

\paragraph{Lemma 12.} If $\lambda \in {\cal Y}_k$, $\lambda \subset
\square_n$, then
$$
 \frac{\nu_{n,k}(\lambda)}{\mu_k(\lambda)} = 
\frac{\prod_{j=1}^{\lambda'(1)} (n+j-1)^{\downarrow j}
\cdot \prod_{j=1}^{\lambda(1)} (n+j-1)^{\downarrow j}}
{(n^2)^{\downarrow k}}
$$

\paragraph{Proof.} Use the hook formula \eqref{eq:hook}. A computation
similar to the one in the proof of \eqref{eq:amus} shows that many of
the terms cancel. We omit the relatively simple details.
\qed

\bigskip
It follows using elementary estimates, which we again omit for the
sake of brevity, that

\paragraph{Theorem 13.} If $n\to \infty$ and $k=k(n)$ is such that
$k=o(n^{2/3})$, then
$$ \frac{\nu_{n,k}(\lambda)}{\mu_k(\lambda)}
\xrightarrow[n\to\infty]{} 1 $$ 
uniformly on the support ${\cal Y}_{n,k}$ of $\nu_{n,k}$ (the set of
diagrams of size $k$ contained in $\square_n$). In particular, the
total variation distance
$$ ||\nu_{n,k}-\mu_k||_1 := \sum_{\lambda \in {\cal Y}_k}
|\nu_{n,k}(\lambda) - \mu_k(\lambda)| \xrightarrow[n\to\infty]{} 0. $$

\vspace{-25.0 pt}
\qed

\bigskip
So in fact, when $k$ is small, $\nu_{n,k}$ is a kind of deformation of
the Plancherel measure $\mu_k$. In particular, for $k$ fixed and $n$
going to infinity, this implies the not-entirely-trivial fact that
$\mu_k$ is a probability measure. We remark that other deformations of
Plancherel measure have been used as a means to study Plancherel
measure itself -- see, e.g., \cite{johansson}. The phenomenon that a
small subtableau of a large random tableau has approximately the
Plancherel distribution was observed also in \cite{morsemckaywilf}
(see also \cite{stanley} for related results)
for 
a random tableau chosen uniformly among \emph{all} tableaux of size
$k$. Recently, it was shown \cite{pittel3} that the footprint of the
$k$ tallest stacks in a random unrestricted plane partition of high
volume also has in the limit the Plancherel distribution.

Another related observation is the easily checked fact that
$$ \sqrt{\alpha(1-\alpha)}\cdot
\tilde{g}_\alpha\left(\frac{u}{\sqrt{\alpha(1-\alpha)}}\right) 
\xrightarrow[\alpha \searrow 0]{} \Omega(u),
$$
i.e. the shape of the level curves of our limit surface $L$ converges
after rescaling to the Plancherel limit curve $\Omega$, as one
approaches the corner of the square. This is consistent with Theorem
13, although is not formally implied by it, as here $k$ is a small
constant times $n^2$. It seems likely that in the regime when $k$
grows like $o(n^2)$, but much faster than $n^{2/3}$, $\nu_{n,k}$ and
$\mu_k$ become mutually singular, even though the limit shapes
coincide.

\subsection{The probability of a square plane partition to have all
parts distinct}

Denote by $M_{nN}$ the total number of $n\times n$ square plane
partitions of $N$. Let ${\cal M}_{nN}$ be the total number of those
partitions with all parts distinct. From Lemma 11 it follows that if
\begin{equation}\label{eq:t1}
\lim_{n,N\to\infty} \frac{n^6}{N}=0,
\end{equation}
then
\begin{equation}\label{eq:t2}
\lim_{n,N\to\infty}\frac{{\cal M}_{nN}}{M_{nN}}=1.
\end{equation}
Our goal is to show that \eqref{eq:t1} is essentially necessary for
\eqref{eq:t2}. To motivate
the statement, notice that
the $k$-th largest part in a partition of $N$ into $n^2$ distinct parts is
$n^2-k+1$, at least. So
${\cal M}_{nN}=0$ unless $N\ge n^2(n^2+1)/2$. 

\paragraph{Theorem 14.} Suppose that $n^4/N\to 0$. (i) If $\lim n^6/N=\infty$,
then
\begin{equation}\label{eq:t3}
\lim_{n,N\to\infty}\frac{{\cal M}_{nN}}{M_{nN}}=0.
\end{equation}
(ii) If $\lim n^6/N=\alpha\in (0,\infty)$ then
\begin{equation}\label{eq:t4}
\lim_{n,N\to\infty}\frac{{\cal M}_{nN}}{M_{nN}}=e^{-\alpha/4}.
\end{equation}

\paragraph{Note.} Thus the reduction to the plane partitions with
distinct parts used in the proof of Theorem 4 is valid if and only if
$n^6=o(m)$.

\paragraph{Proof sketch of Theorem 14.} We prove \eqref{eq:t3}, \eqref{eq:t4} by
determining the asymptotic expressions of ${\cal M}_{nN}$ and $M_{nN}$.

\paragraph{Part 1.}
Begin with ${\cal M}_{nN}$. As in the proof of Lemma 10, we notice that
-- given a linear partition of $N$ into $n^2$ distinct parts -- the
number of the $n\times n$ square (descending) arrangements of these
parts equals the total number of $n\times n$ square Young
tableaux. So, denoting by $p_{nN}$ the total number of all such linear
partitions, and by $d(\square_n)$ the number of all such tableaux, we
obtain
\begin{equation}\label{eq:t7}
{\cal M}_{nN}=p_{nN}\cdot d(\square_n).
\end{equation}
Using the hook formula \eqref{eq:hook}, Euler's summation formula for $\sum_{s=1}^{n-1}
(n-s)\log (n-s)$, and two identities for the Gamma
function (see Bateman \cite{bateman}, Section 1.9)
\begin{eqnarray*}
\sum_{s=1}^ns\log s&=&\int_1^n\log \Gamma
(x)\,dx+\frac{n(n+1)}{2}+\frac{n}{2}\log 2\pi,\\
\log\Gamma (x)&=&
\left(x-\frac{1}{2}\right)\log x-x+\frac{1}{2}\log 2\pi \\ & &
+\int_0^{\infty}
\left[(e^t-1)^{-1}-t^{-1}+\frac{1}{2}\right]t^{-1}e^{-tx}\,dt,\quad x>0.
\end{eqnarray*}
we obtain 
\begin{equation}\label{eq:t11}
d(\square_n)\sim n^{11/12}\sqrt{2\pi}\exp\left(n^2\log n+n^2(-2\log
2+1/2)-\frac{1}{6}+\frac{\log 2}{12}
-C\right),
\end{equation}
where
\begin{equation}\label{eq:t10}
C:=\int_0^{\infty}\left[(e^t-1)^{-1}-
 t^{-1}+\frac{1}{2}-\frac{t}{12}\right]t^{-2}e^{-t}\,dt.
\end{equation}
(A cruder formula 
$$
d(\square_n)\sim \exp(n^2\log n+n^2(-2\log 2+1/2)+O(n\log n))
$$
was obtained, implicitly, in the proof of Lemma 1.)

As for $p_{nN}$, the total number of partitions of $N$ into $n^2$ distinct
parts, it is given by
\begin{equation}\label{eq:t12}
p_{nN}=[q^N t^{n^2}]\prod_{\ell=1}^{\infty}(1+q^{\ell}t).
\end{equation}
From a more general theorem of Vershik and Yakubovich \cite{vershikyak1},
based on \eqref{eq:t12}, Fristedt's conditioning defice \cite{fristedt}, and an
attendant local limit theorem result, it follows that
\begin{equation}\label{eq:t13}
p_{nN}\sim \frac{1}{2\pi N}
\exp\left(n^2\log\frac{N}{n^4}+2n^2-
\frac{n^6}{4N}(1+O(n^4/N))\right).
\end{equation}
Combining \eqref{eq:t11} and \eqref{eq:t13}, we arrive at
\begin{multline}\label{eq:t14} \qquad\qquad
{\cal M}_{nN}\sim\frac{n^{11/12}}{\sqrt{2\pi}N}\exp
\bigl(n^2\log(N/n^3)+n^2(-2\log 2+5/2)\\
-(n^6/4N)(1+O(n^4/N))+C^*\bigr), \qquad\qquad\\
C^*:=-\frac{1}{6}+\frac{\log 2}{12}-C, \qquad\qquad\qquad\qquad\qquad
\qquad\qquad\qquad\qquad\quad\ \ 
\end{multline}
with $C$ defined in \eqref{eq:t10}.

\paragraph{Part 2.} Turn now to $M_{nN}$, the number of all square $n\times n$
plane partitions of $N$. By the MacMahon formula for the number of
plane partitions with at most $n$ rows and $n$ columns,
\begin{equation}\label{eq:t15}
M_{nN}=[q^{N-n^2}]\prod_{\ell=1}^{\infty}(1-q^{\ell})^{-\ell}\cdot
\prod_{\ell>n}(1-q^{\ell})
^{2(\ell-n)}\cdot\prod_{\ell>2n}(1-q^{\ell})^{2n-\ell}.
\end{equation}
(Alternatively, this formula follows from the hook expression for the
generating functions of plane partitions with a given shape discovered
by Stanley \cite{stanleyhook}.)  We will use the techniques from
\cite{pittel3} inspired by Freiman's derivation of the main part of
Hardy-Ramanujan formula for the (linear) partition function, see
Postnikov \cite{postnikov}.

Let us take a close look at the generating function in \eqref{eq:t15}, which we
denote $p_n(q)$. Set $q=e^{-u}$, $\text{Re }u>0$. Taking logarithms,
using
\begin{equation}\label{eq:t16}
\log (1-e^{-mu})^{-1}=\sum_{j\ge 1}\frac{e^{-mju}}{j},
\end{equation}
and changing the summation order, we obtain
\begin{multline}\label{eq:t17}
\qquad\quad\ \ 
\log p_n(e^{-u})=u\sum_{j\ge 1}\frac{1}{uj}\frac{e^{uj}}{(e^{uj}-1)^2}
(1-e^{-unj})^2\\ =
n^2\sum_{j=1}^{\infty}\frac{e^{-unj}}{j}
+u n^3\sum_{j=1}^{\infty}\frac{\psi(unj)}{(unj)^3}
\qquad\qquad\quad
\\
\qquad\qquad\quad\ \ 
-\frac{1}{12}\sum_{j=1}^{\infty}\frac{e^{-uj}}{j}(1-e^{-unj})^2
+u\sum_{j=1}^{\infty}\phi(uj)(1-e^{-unj})^2;\\
\phi(z):=\frac{e^z}{z(e^z-1)^2}-\frac{1}{z^3}+\frac{e^{-z}}{12z};
\qquad\qquad\qquad\qquad\quad\ \ 
\\
\psi(z):=(1-e^{-z})^2-z^2e^{-z}.
\qquad\qquad\qquad\qquad\qquad\qquad\qquad\qquad
\end{multline}
Using \eqref{eq:t16}, read backward, for the first and the third sums,
and Euler's summation formula, with $m=1$, for the second and the
fourth sum, we obtain from \eqref{eq:t17}: if $|u|^2n^4\to 0$, then
\begin{multline}\label{eq:t19}
\qquad\qquad
p_n(e^{-u})\sim \frac{\exp\bigl(n^2(3/2-2\log 2)+(\log
2)/12+D\bigr)}{n^{1/12}}
\cdot (1-e^{-nu})^{-n^2};\\
D:=\int_0^{\infty}\phi(x)\,dx.
\qquad\qquad\qquad\qquad\qquad\qquad\qquad\qquad\qquad\qquad\qquad\ \ 
\end{multline}
(According to Maple the integral of $\psi(x)$ equals $3/2-2\log 2$.)
Now, by \eqref{eq:t15}, and the Cauchy formula
$$
M_{nN}=[q^{N-n^2}]p_n(q)=\frac{1}{2\pi i}
\oint_{z=\rho e^{i\theta}\atop \theta\in (-\pi,\pi]}\frac
{p_n(z)}{z^{N-n^2+1}}\,dz,
$$
where $\rho\in (0,1)$. In light of the last two equations, we set
$\rho=e^{-u_0}$, and select $u_0$ that 
minimizes $-n^2\log (1-e^{-nu})+(N-n^2)u$, that is 
$$
u_0=n^{-1}\log\left(1+\frac{n^3}{N-n^2}\right)\sim \frac{n^2}{N}.
$$ 
Clearly $u_0^2n^4\to 0$. 
Consider $|\theta|\le \theta_n=n^{-2}\epsilon_n$, $\epsilon_n=(n^4/N)^{1/2}$,
so that 
$|\theta|^2n^4\to 0$ as well. It can be shown without much difficulty that
\begin{multline}\label{eq:t21}
\frac{1}{2\pi
i}\int\limits_{|\theta|\le\theta_n}\frac{dz}{(1-e^{-nu})^{n^2}z^{N-n^2+1}}\,dz
=\frac{1}
{2\pi
i}\int\limits_{u_0-i\theta_n}^{u_0+i\theta_n}\frac{e^{v(N-n^2)}}{(1-e^{-nv})^{n^2}}\,dv\\
=\frac{1}{n^{n^2}}\cdot\frac{N^{n^2-1}}{(n^2-1)!}(1+O(n^4/N))
\sim\frac{n}{\sqrt{2\pi}N}\left(\frac{eN}{n^3}\right)^{n^2}.
\end{multline}
(The last integral, extended to the closed contour obtained by connecting the
points $u_0\pm i\theta_n$
with a circular arc centered at the origin, is exactly 
$$
\frac{1}{n^{n^2}}[t^{n^2-1}](1-nt)^{-(N-n^2+n+2)/n}=n^{-n^2}\binom{-n^{-1}(N-n^2+n+2)}{n^2-1}(-n)^{n^2-1},
$$
and the supplementary integral is less than this quantity by a factor
$(u_0/\theta_n)^{n^2}\sim
(n^4/N)^{n^2/2}$.) 
So, using \eqref{eq:t19}, \eqref{eq:t21},
\begin{multline}\label{eq:t22}
\frac{1}{2\pi i}
\oint\limits_{z=\rho e^{i\theta}\atop \theta\in (-\theta_n,\theta_n]}\frac
{p_n(z)}{z^{N-n^2+1}}\,dz\\
\sim\frac{n^{11/12}}{\sqrt{2\pi}N}\exp\biggl(n^2\log\frac{N}{n^3}+n^2(5/2-2\log
2)
+(\log 2)/12+D+o(1)\biggr).
\end{multline}
The proof that the contribution of $\theta\notin [-\theta_n,\theta_n]$ is negligible
compared with the last expression is based on cruder estimates, not unlike those in
\cite{pittel3}, and we omit it. Therefore
\begin{equation}\label{eq:t26}
M_{nN}\sim\frac{n^{11/12}}{\sqrt{2\pi}N}\exp\bigg(n^2\log\frac{N}{n^3}+
n^2(5/2-2\log 2)+(\log 2)/12+D\bigg).
\end{equation}
Comparing \eqref{eq:t26} and \eqref{eq:t14}, we see that
\begin{equation}\label{eq:t27}
\frac{{\cal M}_{nN}}{M_{nN}}=\exp\big(-(n^6/4N)(1+O(n^4/N))+A\big),
\end{equation}
where, recalling the definition of $C$ and $D$,
\begin{eqnarray*}
A&=&-\frac{1}{6}-C-D\\
&=&-\frac{1}{6}-\int\limits_0^{\infty}\left[\left(\frac{1}{e^t-1}-\frac{1}{t}+\frac{1}{2}\right)t^{-2}e^{-t}
+\frac{e^t}{t(e^t-1)^2}-\frac{1}{t^3}\right]\,dt\\
&=&-\frac{1}{6}-\int\limits_0^{\infty}\frac{d}{dt}\left(-\frac{1}{t(e^t-1)}+\frac{1}{2t^2}+\frac{e^{-t}}
{2t^2}\right)\,dt\\
&=&-\frac{1}{6}
 +\lim_{t\downarrow 0}\left(-\frac{1}{t(e^t-1)}+\frac{1}{2t^2}+
\frac{e^{-t}}{2t^2}\right)\\
&=&0.
\end{eqnarray*}
Thus we have
$$
\frac{{\cal M}_{nN}}{M_{nN}}=\exp\big(-(n^6/4N)(1+O(n^4/N))\big),
$$
which proves Theorem 14(i),(ii).
\qed

\subsection{Open problems} 

We conclude with some open problems.

\begin{itemize}
\item {\bf Gaussian fluctuations.} Prove a central limit theorem for
the fluctuations of $g_{\lambda_T^{\lfloor \alpha n^2\rfloor}}$
around $\tilde{g}_\alpha$, and for the fluctuations of the
cotransition measure of $\lambda_T^{\lfloor \alpha n^2\rfloor}$ around
the semicircle distribution, in the spirit of \cite{ivanovolsh}.

\item {\bf Limiting distribution of $l_{n,k}(\pi)$.} Find a scaling
sequence $a_n$ and a distribution function $F$ such that, in the
notation of Theorem 3,
$$
\frac{l_{n,\lfloor \alpha n^2 \rfloor} - 2\sqrt{\alpha(1-\alpha)} n}{a_n}
\xrightarrow[n\to\infty]{\textrm{ in distribution }} F. $$

\item {\bf Limit surface for random Young tableaux of given shape.}
Prove a limit surface theorem for random Young tableaux of other
shapes. In general, one can consider any decreasing function
$f:[0,\infty)\to[0,\infty)$ such that $\int_0^\infty f(x)dx = 1$ as a
\emph{continual Young diagram}, i.e. as a limit of the rescaled graphs
of a sequence of Young diagrams of increasing sizes. We conjecture
that for each such continual diagram $f$, there should exist a limit
surface $L_f$, defined on the domain
$$ D_f := \{ (x,y) : x\ge 0, \ \,0 \le y \le f(x) \} $$
bounded between the $x$-axis and the graph of $f$, that describes the
asymptotic behavior of almost all random Young tableaux of shape
approximated by $f$.

\end{itemize}

\section*{Acknowledgements}
We thank Nati Linial for his suggestion that we work together, which
eventually led to the conception of this paper.

\end{document}